\pdfoutput=1
\RequirePackage{ifpdf}
\ifpdf 
\documentclass[pdftex]{sigma}
\else
\documentclass{sigma}
\fi

\numberwithin{equation}{section}

\usepackage{mathtools} 
\usepackage{stmaryrd} 
\usepackage{physics} 
\usepackage{tensor} 
\usepackage{faktor} 
\usepackage{mathrsfs} 
\usepackage{bbm} 
\usepackage{tikz} 
\usetikzlibrary{cd}
\usepackage{enumitem} 
\setenumerate[1]{label={(\arabic*)}}
\usepackage[titletoc,title]{appendix} 
\newtheorem{Theorem}{Theorem}[section]
\newtheorem*{Theorem*}{Theorem}
\newtheorem{Corollary}[Theorem]{Corollary}
\newtheorem{Lemma}[Theorem]{Lemma}
\newtheorem{Proposition}[Theorem]{Proposition}
 { \theoremstyle{definition}
\newtheorem{Definition}[Theorem]{Definition}

\newtheorem{Example}[Theorem]{Example}
\newtheorem{Remark}[Theorem]{Remark} }
\newcommand{\fa}{\mathfrak{a}}

\newcommand{\fg}{\mathfrak{g}}
\newcommand{\fh}{\mathfrak{h}}
\newcommand{\fk}{\mathfrak{k}}

\newcommand{\fs}{\mathfrak{s}}
\newcommand{\fgl}{\mathfrak{gl}}

\newcommand{\fso}{\mathfrak{so}}
\newcommand{\fiso}{\mathfrak{iso}}

\newcommand{\fX}{\mathfrak{X}}
\newcommand{\fK}{\mathfrak{K}}

\newcommand{\fS}{\mathfrak{S}}
\newcommand{\fV}{\mathfrak{V}}
\newcommand{\fD}{\mathfrak{D}}

\newcommand{\RR}{\mathbb{R}}

\newcommand{\ZZ}{\mathbb{Z}}

\newcommand{\eD}{\mathscr{D}}
\newcommand{\eE}{\mathscr{E}}

\newcommand{\eL}{\mathscr{L}}

\newcommand{\eP}{\mathscr{P}}

\newcommand{\cH}{\mathcal{H}}
\newcommand{\cK}{\mathcal{K}}

\newcommand{\ssB}{\mathsf{B}}
\newcommand{\ssC}{\mathsf{C}}
\newcommand{\ssH}{\mathsf{H}}
\newcommand{\ssZ}{\mathsf{Z}}

\newcommand{\ssS}{\mathsf{S}}
\newcommand{\Sbundle}{\underline{S}}

\newcommand{\eLhat}{\widehat{\eL}}

\newcommand{\betahat}{\widehat{\beta}}
\newcommand{\gammahat}{\widehat{\gamma}}

\newcommand{\Htilde}{\widetilde{H}}

\newcommand{\betatilde}{\widetilde{\beta}}
\newcommand{\gammatilde}{\widetilde{\gamma}}

\newcommand{\thetatilde}{\widetilde{\theta}}
\newcommand{\fatilde}{\widetilde{\fa}}

\newcommand{\fgtilde}{\widetilde{\mathfrak{g}}}

\newcommand{\Wedge}{\mathchoice{{\textstyle\bigwedge}}{\bigwedge}{\bigwedge}{\bigwedge}}	
\newcommand{\Odot}{\mathchoice{{\textstyle\bigodot}}{\bigodot}{\bigodot}{\bigodot}}
\newcommand{\Otimes}{\mathchoice{{\textstyle\bigotimes}}{\bigotimes}{\bigotimes}{\bigotimes}}

\newcommand{\be}{\boldsymbol{e}}

\newcommand{\1}{\mathbbm{1}}

\newcommand{\pair}[2]{\left\langle #1,#2\right\rangle} 
\DeclareMathOperator{\Id}{Id}

\DeclareMathOperator{\End}{End}
\DeclareMathOperator{\Hom}{Hom}
\DeclareMathOperator{\GL}{GL}
\DeclareMathOperator{\Spin}{Spin}
\DeclareMathOperator{\Pin}{Pin}
\DeclareMathOperator{\Orth}{O}
\DeclareMathOperator{\SO}{SO}

\DeclareMathOperator{\ISO}{ISO}

\DeclareMathOperator{\Lie}{Lie}

\DeclareMathOperator{\Cl}{Cl}
\DeclareMathOperator{\Ad}{Ad}
\DeclareMathOperator{\ad}{ad}

\DeclareMathOperator{\Gr}{Gr}

\DeclareMathOperator{\evaluate}{ev}
\DeclareMathOperator{\pr}{pr}
\DeclareMathOperator{\Aut}{Aut}
\DeclareMathOperator{\GrAut}{GrAut}

\begin{document}
\allowdisplaybreaks

\newcommand{\arXivNumber}{2409.11306}

\renewcommand{\PaperNumber}{081}

\FirstPageHeading

\ShortArticleName{Killing (Super)Algebras Associated to Connections on Spinors}

\ArticleName{Killing (Super)Algebras Associated to Connections\\ on Spinors}

\Author{Andrew D.K. BECKETT}

\AuthorNameForHeading{A.D.K.~Beckett}

\Address{University of Edinburgh, UK}
\Email{\href{mailto:adkbeckett@proton.me}{adkbeckett@proton.me}}
\URLaddress{\url{https://adkbeckett.github.io/}}

\ArticleDates{Received November 28, 2024, in final form September 15, 2025; Published online September 30, 2025}

\Abstract{We generalise the notion of a Killing superalgebra, which arises in the physics literature on supergravity, to general dimension, signature and choice of spinor module and Dirac current. We also allow for Lie algebras as well as superalgebras, capturing a set of examples previously defined using geometric Killing spinors on higher-dimensional spheres. Our definition requires a connection on a spinor bundle -- provided by supersymmetry transformations in the supergravity examples and by the Killing spinor equation on the spheres -- and we obtain a set of sufficient conditions on such a connection for the Killing (super)algebra to exist. We show that these Lie (super)algebras are filtered deformations of graded subalgebras of (a generalisation of) the Poincar\'e superalgebra and then study such deformations abstractly using Spencer cohomology. In the highly supersymmetric Lorentzian case, we describe the filtered subdeformations which are of the appropriate form to arise as Killing superalgebras, lay out a classification scheme for their odd-generated subalgebras and prove that, under certain technical conditions, there exist homogeneous Lorentzian spin manifolds on which these deformations are realised as Killing superalgebras. Our results generalise previous work in the 11-dimensional supergravity literature.}


\Classification{17B66; 83E50; 17B56}

\tableofcontents

\section{Introduction}
\label{sec:intro}

{\bf Killing spinors in mathematics and physics.}
In both mathematics and physics, there are special sections of bundles of spinors known as \emph{Killing spinors}, although the term is used somewhat differently across disciplines. In Riemannian geometry, a Killing spinor on a spin manifold~$(M,g)$ is a non-zero section $\epsilon$ of a complex spinor bundle satisfying the equation~${\nabla_{\!X} \epsilon = \lambda X\!\cdot \epsilon
}$
for all vector fields $X$, where $\nabla$ is a lift of the Levi-Civita connection, $\cdot$ denotes the Clifford action and $\lambda$ is a complex constant. We will refer to such spinor fields as \emph{geometric Killing spinors}. Various generalisations have been considered, including allowing $\lambda$ to be a function \cite{Lichnerowicz1987} or replacing $\lambda X$ on the right-hand side of the equation with~$A(X)$ where~$A$ is an endomorphism field of~$TM$.\footnote{The term \emph{generalised Killing spinor} is most commonly used to refer to spinor fields $\epsilon$ in the special case where $A$ is $g$-symmetric \cite{Bar2005,Moroianu2014}, while the more general situation and in particular the $g$-skew case have also been studied~\cite{Habib2012}.}
The geometric context has also been generalised, for example by considering metrics with indefinite signature \cite{Bohle2003} as well as spinors in spin-$\mathbb{C}$ structures \cite{Moroianu1997} rather than spin structures.\looseness=-1

In physics, a not-dissimilar notion of Killing spinors arises in supergravity theories as a~condition for supersymmetry of ``bosonic'' backgrounds, i.e., those in which all fermions vanish. Such solutions can be considered as Lorentzian spin manifolds (or generalisations thereof, see~\cite{Lazaroiu2019,Lazaroiu2019_1}), typically equipped with some other ``background'' structure, such as higher gauge theory data, which induces differential forms on the base. Killing spinors are then spinor sections subject to a~differential equation parametrised by the other background data. For example, in 11-dimensional supergravity, the Killing spinor equation takes the form
\[
 \nabla_X \epsilon = \frac{1}{24}X \cdot F \cdot \epsilon - \frac{1}{8}F \cdot X \cdot \epsilon,
\]
where $F$ is a 4-form which arises as the field strength of a 3-form gauge potential \cite{Pakis2003}, and in this case we take the spinors to be real. This equation and its analogues in other supergravity theories arise by demanding that a supersymmetry variation of the gravitino field vanishes. Typically, there are further constraint equations (algebraic in the spinor) arising from vanishing variations of other fermionic fields which supplement the differential Killing spinor equation, although (outside of Section~\ref{sec:alg-kse}) in the present work our focus will be on the differential equation.

{\bf Killing superalgebras.}
The Killing spinors and a restricted set of Killing vectors which~preserve all the data of a supergravity background can be arranged into a \emph{Killing superalgebra}~\mbox{\cite{Figueroa-OFarrill2007_1,Figueroa-OFarrill2005}}, a Lie superalgebra which can be thought of as a supersymmetric generalisation of the algebra of infinitesimal isometries.\footnote{It should be noted that this notion is distinct from that of (geometric) Killing spinors as isometries of supermanifolds (or superisometries) found, e.g., in \cite{Alekseevsky1998}, although we will discuss cases where the two coincide.
	}
The brackets of this algebra are defined using Lie derivatives (of vector and spinor fields) and a map known as the \emph{Dirac current}, constructed via an invariant inner product on the (s)pinor module, which allows supergravity Killing spinors to be ``squared'' (paired symmetrically) to produce Killing vector fields of the restricted type.\footnote{A previous version of this work (as well as the author's Ph.D.~Thesis \cite{Beckett2024} from which much of its content is drawn) used the term \emph{squaring map} for this construction, following \cite{Figueroa-OFarrill2012,Figueroa-OFarrill2014,Hustler2016}, but this usage conflicts with other works, e.g., \cite{Cortes2021,Shahbazi2024_1} where it denotes a polyvector-valued map of which the vector-valued map here called \emph{Dirac current} is a component. The updated terminology is also consistent with other related works \cite{Beckett2021,deMedeiros2016,deMedeiros2018,Figueroa-OFarrill2017_1,Figueroa-OFarrill2017}.}

One might wonder whether other types of spinor sections, such as geometric Killing spinors, can be used to define Killing superalgebras or some analogue thereof. In Riemannian signature, a~symmetric Dirac current map does not always exist; instead, an inner product on spinors can be used to construct a \emph{skew}-symmetric pairing of spinors into vector fields. Inspired by the Killing superalgebras of supergravity, Figueroa-O'Farrill considered the problem of whether, using this pairing, one can organise geometric Killing spinors into a Killing \emph{algebra} (rather than a superalgebra) \cite{Figueroa-OFarrill2008_except}. He demonstrated that this is indeed the case for the standard \mbox{7-,}~8- and 15-spheres and that the resulting algebras are, respectively, a geometric realisation of~$\mathfrak{so}(9)$ which exhibits the triality of an $\mathfrak{so}(8)$ subalgebra and realisations of the exceptional Lie algebras~$E_8$ and~$F_4$. On the other hand, he showed that it is not possible in general to construct such an algebra on an arbitrary Riemannian spin manifold: the vector fields formed by pairing Killing spinors may not be Killing vectors (although they are \emph{conformal} Killing), and the Jacobi identities provide further obstructions.

Although, as mentioned above, there are conventional choices of Dirac current in Lorentzian or Riemannian signature, in any signature there may be many choices of such maps which may be either symmetric or skew-symmetric -- these choices were classified in~\cite{Alekseevsky1997}. With any such choice, it is natural to consider whether one might be able to define a notion of Killing algebra or superalgebra on a spin manifold using some special type of spinor, and to ask what structure such algebras might have and what they can tell us about the geometry. Our goal in this work will be to investigate these questions by considering equations of the form
$\nabla_X \epsilon = \beta(X) \epsilon
$
where $\beta$ is a one-form with values in endomorphisms of the spinor bundle, $\epsilon$ is a section of the spinor bundle and $X$ is an arbitrary vector field. All types of Killing spinors discussed above are solutions to such an equation. If we define a connection on spinors by $D=\nabla-\beta$, then solutions to this equation are precisely those spinor sections which are parallel with respect to $D$;\footnote{Spinors of this type have been referred to as \emph{generalised Killing spinors} \cite{Cortes2021,Lazaroiu2016}, but the term \emph{differential spinors} is used in more recent work in this direction \cite{Shahbazi2024_2,Shahbazi2024_1} in order to avoid confusion with the use of that term, e.g., in \cite{Bar2005,Moroianu2014}.}
we will find that the existence of a Killing (super)algebra imposes conditions on the connection $D$. The question of associating such an algebra to an \textit{a priori} unspecified connection on spinors was briefly considered in the standard Lorentzian context \cite{Hustler2016}, but our approach will be more systematic: we allow for arbitrary dimension, signature, $N$-extension of spinors and choice of Dirac current (including its symmetry); moreover, we consider more carefully the global geometric structures required and present the conditions on $D$ in a different manner. We will also go much further in our analysis of the general structure of these algebras.

{\bf Classification of supergravity backgrounds and homogeneity.}
Much of the interest in Killing superalgebras in the literature is driven by their application to the classification of supersymmetric supergravity backgrounds, an important open problem closely related to the study of string theory vacua; see the recent review \cite{Gran2019} and the references therein. An important organising principle in this effort is the \emph{fraction of supersymmetry}, denoted $\nu$. The \emph{number of supersymmetries}, i.e., the dimension of the space of Killing spinors $\fS_D$ of a background, is bounded above by the number of supercharges of the theory, and $\nu$ is the fraction of this bound which is attained. For example, 11-dimensional supergravity has 32 supercharges, so \smash{$\nu=\frac{\dim\fS_D}{32}$}. Killing superalgebras are especially relevant for analysis of the \emph{highly supersymmetric} case $\nu>\frac{1}{2}$ due to a result known as the \emph{homogeneity theorem} (our Theorem~\ref{thm:homogeneity} is a gloss) which implies that if $\nu>\frac{1}{2}$, the tangent space at every point of the background is spanned by Killing vector fields of the restricted kind; in particular, the background is (up to local isomorphism) a homogeneous space for the even part of the Killing superalgebra \cite{Figueroa-OFarrill2012,Figueroa-OFarrill2014,Hustler2016}. Unfortunately, the theorem is only valid in Lorentzian signature (and for a particular class of Dirac currents), since its proof relies on the causality of the Dirac current and properties of subspaces of Lorentzian inner product spaces. Nonetheless, it can happen that local homogeneity still holds in examples in other signatures, such as the Killing algebras of geometric Killing spinors on higher-dimensional spheres \cite{Figueroa-OFarrill2008_except}.

{\bf Spencer cohomology and filtered deformations.}
The homogeneity theorem suggests that the classification of highly supersymmetric backgrounds might be approached from a new direction; namely, via a classification of the possible Killing superalgebras. This was first considered for 11-dimensional backgrounds by Figueroa-O'Farrill and Santi \cite{Figueroa-OFarrill2017_1,Figueroa-OFarrill2017,Santi2022} who showed that these superalgebras are \emph{filtered subdeformations} of the Poincar\'e superalgebra; that is, they are filtered Lie superalgebras whose associated graded superalgebras are graded subalgebras of the Poincar\'e superalgebra.\footnote{See Section~\ref{sec:filtered-def} and Appendix~\ref{sec:grading-filtration-superspace} for definitions and further background.} Filtered deformations of graded Lie superalgebras were studied by Cheng and Kac \cite{Cheng1998}, who showed that they are governed by the (generalised) Spencer cohomology of the graded superalgebra being deformed. It was shown by exploiting the homogeneity theorem that in the highly supersymmetric case, the subdeformations corresponding to Killing superalgebras can be parametrised in terms of the Spencer cohomology of the Poincar\'e superalgebra itself, rather than that of the subalgebra being deformed, a significant simplification~\cite{Figueroa-OFarrill2017_1}. Moreover, it is possible to explicitly reconstruct (up to local isometry) a highly supersymmetric background from the Spencer deformation data associated to its Killing superalgebra. This point of view allows the classification problem for highly supersymmetric 11-dimensional backgrounds to be transformed into an algebro-geometric problem; namely, one of determining the admissible Spencer deformation data \cite{Santi2022}.

{\bf Aim and structure.}
Our primary goal in this work is to extend the results of \cite{Figueroa-OFarrill2017_1} to a more general setting. Many of the same techniques are employed, but care must be taken with some technical points which are significantly simplified in that more narrow context, including in the definition of the Killing superalgebra itself, as we will discuss at length. The paper is organised as follows.

In Section~\ref{sec:background}, we provide some necessary background on our \emph{flat model (super)algebras} -- Poincar\'e superalgebras and their generalisations (Section~\ref{sec:flat-model-algs}) whose filtered subdeformations we will later study -- and on the filtered deformations and Spencer cohomology (Section~\ref{sec:spencer-thy}) of graded Lie superalgebras. Additional background, including conventions for (pseudo-)inner product spaces and spinors (Appendix~\ref{sec:inner-products-spinors}), ordinary and spinorial Lie derivatives (Appen\-dix~\ref{sec:lie-der}), gradings and filtrations (Appendix~\ref{sec:grading-filtration-superspace}) and homogeneous spaces (Appendix~\ref{sec:homog-spaces}) can be found in Appendix~\ref{app:further-background}.

In Section~\ref{sec:killing-spinors-superalgebras}, we first tackle the question of whether, given a connection $D$ on a bundle of spinors, one can define a Lie (super)algebra using the $D$-parallel spinors and some restricted set of Killing vectors. The answer is given in the form of Theorem~\ref{thm:killing-algebra-exist}, our first main result. We then formally define a class of \emph{admissible connections} and their associated Killing (super)algebras (see Definition~\ref{def:killing-spinor}). Our second main result, Theorem~\ref{thm:killing-algebra-filtered}, shows that each such algebra is a~filtered subdeformation of one of the generalisations of the ($N$-extended) Poincar\'e superalgebra discussed in \cite{Alekseevsky1997}, called here \emph{flat model $($super$)$algebras}. This generalises the result known in the 11-dimensional Lorentzian case \cite[Theorem~12]{Figueroa-OFarrill2017_1} and implicitly deduced in other dimensions~\cite{Beckett2021,deMedeiros2016,deMedeiros2018}. Some comments are made in Section~\ref{sec:alg-kse} on the extension of these results to the case of Killing spinors under algebraic constraints of the type imposed in supergravity by the vanishing variations of fermions other than the gravitino.

In Section~\ref{sec:filtered-def-poincare}, we temporarily forgo geometry and approach the issue from the algebraic angle, considering filtered subdeformations of flat model superalgebras through the lens of Spencer cohomology. For the sake of clarity, we work with Lie superalgebras throughout this section, although the Lie algebra case is entirely analogous, at least for the first two subsections. From Section~\ref{sec:homogeneity} onward, we work in Lorentzian signature and with \emph{causal} Dirac currents in order to take advantage of the homogeneity theorem (see Theorem~\ref{thm:homogeneity}), which we exploit to develop the notion of a \emph{geometrically realisable} highly supersymmetric subdeformation of a Poincar\'e superalgebra (see Definition~\ref{def:geom-real}) as well as related notions of \emph{admissible} and \emph{integrable} Spencer cocycles (see Definitions~\ref{def:admissible-class} and~\ref{def:integrable-cocycle}). Theorem~\ref{thm:admiss-coycle-integr} is the main result of this section, describing the ``integration'' of such a cocycle into a realisable subdeformation. In Section~\ref{sec:classification-odd-gen}, we discuss a scheme for classifying realisable subdeformations using Spencer cocycle data, summarised in Proposition~\ref{prop:lie-pair-subdef-corresp}. The approach of this section is similar to existing work on the 11-dimensional Lorentzian case \cite{Figueroa-OFarrill2017_1, Santi2022}, but our more general context introduces new complications, particularly in the first two subsections (where there is no homogeneity assumption) and on the issues of integrability and classification. We make some comments in Remark~\ref{rem:relaxing-homogeneity-assumption} on the extent to which the assumptions of Lorentzian signature and symmetric, causal Dirac current can be relaxed.

Finally, in Section~\ref{sec:high-susy-spin-mfld}, we again exploit homogeneity to work our way back to the geometric perspective from the algebraic by discussing the ``reconstruction'' of a supersymmetric background (as a homogeneous space) from a highly supersymmetric filtered subdeformation, realising the latter as (a subalgebra of) a Killing superalgebra. This culminates in our final result, Theorem~\ref{thm:homog-spin-mfld}, generalising the 11-dimensional case \cite[Theorem~13]{Figueroa-OFarrill2017_1}.

In a companion paper to this work \cite{Beckett2024_3}, we provide some simple 2-dimensional examples to demonstrate some of the theory presented here.

\section{Background}
\label{sec:background}

In this section, we provide some necessary mathematical preliminaries and establish nomenclature and notation. Further background is provided in Appendix~\ref{app:further-background}, which readers unfamiliar with spinors, filtrations, $\ZZ$-gradings on superalgebras, Lie superalgebra cohomology or homogeneous spaces, or those who require a review, may wish to read first.

\subsection{Flat model algebras}
\label{sec:flat-model-algs}

The Euclidean and Poincar\'e algebras are the isometry algebras of Euclidean and Minkowski space respectively. They can be considered as ``flat models" for isometry algebras in Riemannian and Lorentzian geometry in the informal sense that they are the isometry algebras of the maximally symmetric flat spaces in their respective settings. A more precise meaning can be given to this using Cartan geometry, in which a geometric structure on a manifold is viewed as a local equivalence to some homogeneous space -- see \cite{Wise2010} for a pedagogical introduction emphasising links to pseudo-Riemannian geometry and application to general relativity in particular and \cite{Cap2009_1} for a more systematic treatment. Work by \v{C}ap and Neusser \cite{Cap2009} shows that, under certain regularity conditions, one can assign a Lie algebra of symmetries to a Cartan geometry, and that this algebra is a \emph{filtered deformation}\footnote{Filtrations are discussed in Appendix~\ref{sec:filtered-alg} and filtered deformations in Section~\ref{sec:filtered-def}. The term ``filtered deformation'' is not used in \cite{Cap2009} but was previously introduced by Cheng and Kac \cite{Cheng1998}.}
of a subalgebra of the flat model algebra. Pseudo-Riemannian manifolds with signature $(p,q)$ can be considered as Cartan geometries for the local model $\ISO(p,q)/\Orth(p,q)$ with a Cartan connection encoding both the spin connection and a soldering form. If this spin connection is Levi-Civita, the symmetry algebra of the Cartan connection is precisely the Lie algebra of infinitesimal isometries (Killing vectors), and the result discussed above says that it is filtered, with its associated graded Lie algebra being a subalgebra of $\fiso(p,q)$.

Although we do not use the formalism of Cartan geometry in this work, we can view the result on the structure of Killing superalgebras already mentioned in the introduction as a ``super'' generalisation of this phenomenon in pseudo-Riemannian geometry, and it would be subsumed by a generalisation of the \v{C}ap--Neusser result \cite{Cap2009} to super-Cartan geometry. Our (much more explicit) version of this result is Theorem~\ref{thm:killing-algebra-filtered}, which tells us that Killing (super)algebras are filtered subdeformations of the $\ZZ$-graded (super)algebras we will introduce in this section. Motivated by the connection to Cartan geometry, we will refer to these $\ZZ$-graded (super)algebras as \emph{flat model $($super$)$algebras}.

See Appendix~\ref{sec:inner-products-spinors} for further background on inner product spaces and related structures including Clifford algebras and spinors.

\subsubsection{Isometry algebra of an inner product space}

Let $(V,\eta)$ be a (pseudo-)inner product space. We denote by $\fiso(V,\eta)$ the isometry algebra of (the affine space modelled on) $(V,\eta)$ and by $\fso(V,\eta)$ the Lie algebra of $\eta$-skew-symmetric endomorphisms of $V$, although we usually suppress $\eta$ in the notation. As a vector space $\fiso(V)=V\oplus\fso(V)$, and it is equipped with the Lie bracket
\begin{align}\label{eq:flat-bracket}
 \comm{A}{B} = AB - BA,\qquad
 \comm{A}{v} = Av,\qquad
\comm{v}{w} = 0,
\end{align}
where $A,B\in\fso(V)$ and $v,w\in V$. In the Euclidean case $\fiso(V)$ is known as the \emph{Euclidean algebra}; in the Lorentzian case it is the \emph{Poincar\'e algebra} and $\fso(V)$ is the \emph{Lorentz algebra}.

We write $\fiso(p,q):=\fiso(\RR^{p,q})$ where $\RR^{p,q}$ is the standard inner product space with signature~$(p,q)$; if $\eta$ has signature $(p,q)$, the isomorphism $(V,\eta)\cong\RR^{p,q}$ induces a Lie algebra isomorphism $\fiso(V)\cong\fiso(p,q)$.

\subsubsection{Isometry algebra with spinors}
\label{sec:isom-alg-spinors}

We now construct a $\ZZ$-graded Lie algebra or Lie superalgebra\footnote{See Appendix~\ref{sec:graded-superalg} for a definition.} $\fs=\fs_{-2}\oplus\fs_{-1}\oplus\fs_{0}$ whose even part is $\fiso(V)$ as follows. Let
$\fs_0 = \fso(V)$, $
 \fs_{-1} = S$, $
 \fs_{-2} = V$, $
 \fs_i = 0$ otherwise,
where $(V,\eta)$ has signature $(p,q)$ and $S$ is a (possibly $N$-extended) spinor module; that is, we take
\begin{itemize}\itemsep=0pt
	\item $S = \ssS_N := N\ssS_1$ for some integer $N\geq 1$ for $p-q\neq 0,4\mod 8$, where $S_1$ is the unique irreducible real spinor module of $\Cl_{\overline 0}(p,q)$,
	\item $S = \ssS_{(N_+,N_-)} = N_+\ssS_+\oplus N_-\ssS_-$ for integers $N_+, N_-\geq 0$ not both zero if $p-q= 0,4\mod 8$ where $S_\pm$ are the two (non-isomorphic) irreducible real spinor modules of $\Cl_{\overline 0}(p,q)$~-- without loss of generality, we can take $N_+\geq N_-$.
\end{itemize}
We have used conventions for spinor modules explained in Appendix~\ref{sec:pinor-spinor}. We let the bracket on $\fs_{\overline 0}=\fiso(V)$ be the one defined by \eqref{eq:flat-bracket} and note that this respects the $\ZZ$-grading. In order to extend this to a Lie (super)bracket $\fs\otimes\fs\to\fs$ in a way which respects the $\ZZ$-grading, we must specify maps $\fs_{-1}\otimes\fs_0\to\fs_{-1}$ and $\fs_{-1}\otimes\fs_{-1}\to\fs_{-2}$ -- all other components must be zero. We take the first of these maps to be the natural action of $\fso(V)$ on $S$, so the brackets are given by
\begin{align}\label{eq:flat-spinor-bracket}
 \comm{\epsilon}{\zeta} = \kappa(\epsilon,\zeta),\qquad
 \comm{A}{\epsilon} = A\cdot\epsilon,\qquad
 \comm{v}{\epsilon} = 0,
\end{align}
for $A\in\fso(V)$, $v\in V$ and $\epsilon,\zeta\in S$, where we leave the map $\kappa\colon \Wedge^2 S\to V$ (for a Lie algebra) or~${\kappa\colon \Odot^2 S\to V}$ (for a Lie superalgebra) to be determined. The action of $A\in\fso(V)$ on $\epsilon\in S$ is defined in Appendix~\ref{sec:pin-spin-group} (see equation~\eqref{eq:sov-spinor-action}).

It remains to check the Jacobi identities involving elements of $\fs_{-1}$. The $[\fs_0,\fs_0,\fs_{-1}]$ identity is satisfied precisely because the $\fs_0\otimes\fs_{-1}\to\fs_{-1}$ component of the bracket is a module action. The $[\fs_0,\fs_{-1},\fs_{-2}]$ and $[\fs_{-2},\fs_{-1},\fs_{-1}]$ identities are trivial since $\comm{\fs_{-1}}{\fs_{-2}}=\comm{\fs_{-2}}{\fs_{-2}}=0$. The remaining component is $[\fs_0,\fs_{-1},\fs_{-1}]$, which is satisfied if and only if $\kappa$ is $\fso(V)$-equivariant. Thus we make the following definition.

\begin{Definition}[flat model (super)algebra, Dirac current]\label{def:flat-model-alg}
	Let $(V,\eta)$ be an inner product space and $S$ an (extended) real spinor module of $\fso(V)$. Then
	\begin{itemize}\itemsep=0pt
		\item A \emph{flat model algebra} is a $\ZZ$-graded Lie algebra $\fs$ with underlying $\ZZ$-graded vector space $V\oplus S\oplus \fso(V)$ and bracket defined by \eqref{eq:flat-bracket} and \eqref{eq:flat-spinor-bracket} for some $\fso(V)$-equivariant map $\kappa\colon \Wedge^2 S\to V$.
		\item A \emph{flat model superalgebra} is a $\ZZ$-graded Lie superalgebra $\fs$ with underlying $\ZZ$-graded vector space $V\oplus S\oplus \fso(V)$ and bracket defined by \eqref{eq:flat-bracket} and \eqref{eq:flat-spinor-bracket} for some $\fso(V)$-equivariant map $\kappa\colon \Odot^2 S\to V$.
	\end{itemize}
	In either case, $\kappa$ is referred to as the \emph{Dirac current}.
	
	Furthermore, we say that $\fs$ is \emph{minimal} if $S=\ssS_1=\ssS$ or $\ssS_{(1,0)}=\ssS_+$ and that it is \emph{$N$-extended} if $S=\ssS_N$ for $N>1$ or $S=\ssS_{(N_+,N_-)}$ for $N=N_++N_->1$; in the latter case we also say it is \emph{$(N_+,N_-)$-extended}.
\end{Definition}

We note that all flat model (super)algebras are referred to as \emph{$N$-$($super$)$-extended Poincar\'e algebras}, \emph{$\pm N$-extended Poincar\'e algebras} \cite{Alekseevsky1997} or \emph{Poincar\'e Lie superalgebras} \cite{Cortes2020,Gall2021}, but we reserve the name \emph{Poincar\'e} for Lorentzian signature. Our $\ZZ$-grading convention follows that used in previous work on Spencer cohomology and supergravity \cite{Beckett2021,deMedeiros2016,deMedeiros2018,Figueroa-OFarrill2017_1,Figueroa-OFarrill2017}.

Since inner product spaces of equal signature are isomorphic, the data required to specify a~flat model (super)algebra are the signature $(p,q)$, the integer $N$ or pair of integers $(N_+,N_-)$ and the Dirac current $\kappa$. The possible Dirac currents were classified by Alekseevsky and Cort\'es \cite{Alekseevsky1997} by reducing the problem to a classification of $\fso(V)$-invariant bilinears on the irreducible pinor module $\ssS$. More recently, this classification has been improved (at least in the superalgebra case) by describing a necessary and sufficient condition for two different Dirac currents to define isomorphic superalgebras \cite{Cortes2020}, and the possible $R$-symmetry groups -- a particular group of Lie superalgebra automorphisms which serves as an algebraic invariant -- for a wide class of Dirac currents have been described, the comparison of which provides a sufficient condition to distinguish non-isomorphic algebras which is simpler to check \cite{Gall2021}.

Our main focus in this work, or at least its latter half, will be on Lorentzian signature and ($N$-extended) Poincar\'e \emph{super}algebras, mainly due to important technical simplifications which arise in this case from the homogeneity theorem (see Theorem~\ref{thm:homogeneity}). Nonetheless, we will provide a significant amount of formalism in general signature. We will not consider conformal superalgebras or central extensions, nor any other generalisation, although flat model algebras extended by $R$-symmetry will be considered in upcoming work. There is a significant amount of literature on all of these generalisations: for the Lorentzian case, see the classic work of Strathdee (which also discusses the representation theory of the superalgebras) \cite{Strathdee1986} and the review of Van Proeyen \cite{VanProeyen1999}; for general signature, see \cite{Alekseevsky1997,Alekseevsky2005,DAuria2001_1} and also the more recent work \cite{Gall2021}.

\begin{Remark}
	The Dirac current can often be seen as a component of a polyvector-valued (or, dually, polyform-valued) \emph{squaring map} which has been used to study generalisations of Killing spinors algebro-geometrically \cite{Cortes2021}. The other components of this map are often called \emph{Dirac bilinears} or \emph{Dirac $p$-forms}; in the $N$-extended case, these forms typically ``carry internal indices'', i.e., they are valued not in $\Wedge^p V$ for some $p$ but in $\Wedge^p V\otimes\Delta$ where $\Delta$ is some auxiliary space.
	
	The extension of the (super)algebras described above to include these higher forms in the odd-odd components of the bracket has been considered, e.g., in~\cite{Alekseevsky2005,Figueroa-OFarrill2009}. The so-called FDAs (``free'' differential (graded) algebras) of the rheonomic formalism of supergravity \cite{DAuria1982,DAuria1980}, which have been put on a rigorous mathematical footing as $L_\infty$-algebras \cite{Schreiber2015,Schreiber2009,Stasheff2018}, also incorporate these forms, albeit in a different manner.
	All of these extensions can provide information about higher gauge theory data on supergravity backgrounds which is unfortunately invisible to the Lie (super)algebras discussed in this work. We hope to understand how the perspective proposed here can be extended to include such gauge-theoretic information in future work -- and we note that the natural setting for doing so would be higher Cartan geometry as in \cite{Schreiber2009} -- but such considerations are beyond our current scope.
	
	Notwithstanding the comments above, note that although higher Dirac forms are not used in our flat model (super)algebras, they nonetheless do appear in explicit descriptions of their filtered deformations (such as those found in \cite{Beckett2019,deMedeiros2016,deMedeiros2018,Figueroa-OFarrill2017} and the companion work \cite{Beckett2024_3}), as they provide a convenient way to parametrise bilinear maps on the spinor module (e.g., the maps denoted $\gamma$ appearing in \eqref{eq:22-cocycle-comps} and \eqref{eq:def-brackets-general}).
\end{Remark}

\subsection{Filtered deformations and Spencer cohomology}
\label{sec:spencer-thy}

Here we present some background material on the Spencer cohomology and filtered deformations of $\ZZ$-graded Lie superalgebras. Our main source for this section is the work of Cheng and Kac~\cite{Cheng1998}. The base field is assumed to be either $\RR$ or $\mathbb{C}$. We recommend that readers unfamiliar with gradings and filtrations of vector (super)spaces and (super)algebras read Appendix~\ref{sec:grading-filtration-superspace} before this section.

\subsubsection{Spencer cohomology}
\label{subsec:spencer-cohomology}

Let $\fg=\bigoplus_{k=-h}^\infty\fg_k$ be a $\ZZ$-graded Lie superalgebra of finite depth\footnote{The Spencer cohomology theory in its original form was only applicable to algebras of depth 1; for those of greater depth Cheng and Kac use the term \emph{generalised} Spencer cohomology. We do not make this distinction.}
$h$ such that $\dim\fg_k<\infty$ for all $k$. Ultimately, we will be interested in the case where $\dim\fg<\infty$, but our source \cite{Cheng1998} mainly treats examples where $\fg$ itself is infinite-dimensional, and the theory presented here applies in any case. We let $\fg_-=\bigoplus_{k<0}\fg_k$ and note that it is a finite-dimensional subalgebra.


We define the Spencer cochain complex $\qty(\ssC^\bullet(\fg_-;\fg),\partial^\bullet)$ as follows. For $p<0$, we set $\ssC^p(\fg_-;\fg)\allowbreak:=0$, while for $p\geq 0$, we have
$\ssC^p(\fg_-;\fg) := \qty(\Wedge^p \fg_-^*) \otimes \fg \cong \Hom\qty(\Wedge^p \fg_-,\fg)$,
where $\Wedge^p$ is meant in the super-sense (see equation \eqref{eq:super-wedge}). The Spencer differential $\partial^\bullet\colon \ssC^\bullet(\fg_-;\fg)\to \ssC^{\bullet+1}(\fg_-;\fg)$ is given by the following formula: for $\phi\in \ssC^p\qty(\fg_-;\fg)$, $X_1,\dots, X_r\in\fg_{\overline{0}}$ and $Y_1,\dots,Y_s\in\fg_{\overline{1}}$, where $r+s=p+1$,
\begin{gather}
	(\partial\phi)\qty(X_1,\dots,X_r,Y_1,\dots,Y_s)\nonumber\\
\qquad	= \sum_{i=1}^r (-1)^{i+1} [X_i,\phi\bigl(X_1,\dots,\hat{X}_i,\dots,X_r,Y_1,\dots,Y_s\bigr)]	\nonumber\\
\phantom{\qquad=}{}+ (-1)^{r} \sum_{s=1}^r [Y_j,\phi \bigl(X_1,\dots,X_r,Y_1,\dots,\hat{Y}_j,\dots,Y_s)]\nonumber	\\
\phantom{\qquad=}{}+ \sum_{1\leq i<j\leq r} (-1)^{i+j} \phi \bigl(\comm{X_i}{X_j},X_1,\dots,\hat{X}_i,\dots,\hat{X}_j,\dots,X_r,Y_1,\dots,Y_s\bigr)\nonumber	\\
\phantom{\qquad=}{}+ \sum_{i=1}^r\sum_{j=1}^s (-1)^{i} \phi \bigl(X_1,\dots,\hat{X}_i,\dots,X_r,\comm{X_i}{Y_j},Y_1,\dots,\hat{Y}_j,\dots,Y_s\bigr)	\nonumber\\
\phantom{\qquad=}{}+ \sum_{1\leq i<j\leq s} \phi \bigl(\comm{Y_i}{Y_j},X_1,\dots,X_r,Y_1,\dots,\hat{Y}_i,\dots,\hat{Y}_j,\dots,Y_s\bigr),\label{eq:CEdiffLiesuper}
\end{gather}
where the hat decorations denote omission of an entry. It can be checked that $\partial^2=0$, and so $\qty(\ssC^\bullet(\fg_-;\fg),\partial^\bullet)$ is indeed a cochain complex. This is nothing but the standard Chevalley--Eilenberg complex for $\fg_-$ with values in $\fg$, where the former acts on the latter via the restriction of the adjoint representation. However, since $\fg$ is graded, there is additional structure.

The grading on $\fg$ induces a grading on $\fg^*$ defined by $(\fg^*)_k:=(\fg_{-k})^*$. This in turn induces a grading on the space of $p$-cochains, in the usual manner for tensor products of graded spaces. Equivalently, we can consider a $p$-cochain to be homogeneous of degree $d$ if it has degree $d$ as a~map $\Wedge^p \fg_-\to\fg$. We then denote the space of $p$-cochains of degree $d$ as $\ssC^{d,p}(\fg_-,\fg)$,
\begin{equation}\label{eq:spencer-complex-grading}
	\ssC^{d,p}(\fg_-;\fg) = \qty(\Wedge^p \fg_-^*\otimes \fg)_d = \Hom\qty(\Wedge^p \fg_-;\fg)_d,
\end{equation}
so that
\[
	\ssC^p(\fg_-;\fg) = \bigoplus_{d\in\ZZ} \ssC^{d,p}(\fg_-;\fg) \qquad\text{and} \qquad\ssC^{d,\bullet}(\fg_-;\fg) := \bigoplus_{p\in\ZZ} \ssC^{d,p}(\fg_-;\fg).
\]
This grading is preserved by the Spencer differential, so by restricting the differential to the subspaces of cochains with degree $d$, for each $d\in\ZZ$ we have a subcomplex $\bigl(\ssC^{d,\bullet}(\fg_-;\fg),\partial^\bullet\bigr)$.

Finally, we note that the adjoint action of $\fg_0$ on $\fg$ preserves the grading, thus so does the natural action of $\fg_0$ on cochains induced by the adjoint representation (which also trivially preserves the homological grading); each $\ssC^p(\fg_-;\fg)$ is a graded $\fg_0$-module. This action also preserves the differential and thus all of the structure on the subcomplex $\bigl(\ssC^{d,\bullet}(\fg_-;\fg),\partial^\bullet\bigr)$.


We define the cocycles, coboundaries and cohomology of the $d$-subcomplex in the usual way
\begin{gather*}
	\ssZ^{d,p}(\fg_-;\fg) = \ker\qty(\partial^p\colon \ssC^p(\fg_-;\fg)\to \ssC^{p+1}(\fg_-;\fg)),	\\
	\ssB^{d,p}(\fg_-;\fg) = \Im\qty(\partial^p\colon \ssC^{p-1}(\fg_-;\fg)\to \ssC^{p}(\fg_-;\fg)),	\qquad
	\ssH^{d,p}(\fg_-;\fg) = \faktor{\ssZ^{d,p}(\fg_-;\fg)}{\ssB^{d,p}(\fg_-;\fg)},
\end{gather*}
and these are all naturally $\fg_0$-modules, since the action of $\fg_0$ preserves $\bigl(\ssC^{d,\bullet}(\fg_-;\fg),\partial^\bullet\bigr)$. Of course, the space of $p$-cocycles of the whole complex is simply the direct sum of the graded $(d,p)$-cocycles, and similarly for the coboundaries, whence we have a natural isomorphism of $\fg_0$-modules
\[
	\ssH^p(\fg_-;\fg) \cong \bigoplus_{d\in\ZZ}\ssH^{d,p}(\fg_-;\fg),
\]
where the former is the cohomology of the whole complex.

The definition of the differential \eqref{eq:CEdiffLiesuper} on its own is not very enlightening, so let us demonstrate it explicitly on cochains of low homological degree and give and interpretation of the spaces of cocycles and the cohomologies. Fix $d\in\ZZ$ and homogeneous elements $X,Y,Z\in\fg_-$; we denote their degree by $|X|$ etc. Then for $\phi\in\ssC^{d,0}(\fg_-;\fg)=\fg_d$, we have
$(\partial \phi)(X) = \comm{X}{\phi}$,
which shows that $(d,0)$-cocycles are $\fg_-$-invariants in $\fg_d$; equivalently, they are degree-$d$ elements of $\fg$ which centralise $\fg_-$, and $(d,1)$-coboundaries are (super)derivations $\fg_-\to\fg$ given by restricting inner derivations of $\fg$ induced by elements of degree $d$. For $\phi\in\ssC^{d,1}(\fg_-;\fg)=\Hom(\fg_-;\fg)_d$,
\[
	(\partial \phi)(X,Y) = \comm{X}{\phi(Y)} - \comm{Y}{\phi(X)} -\phi(\comm{X}{Y}),
\]
so $(d,1)$-cocycles are (super)derivations $\fg_-\to\fg$ of degree $d$. Thus we find that $\ssH^{d,0}(\fg_-;\fg) \cong (\fg_d)^{\fg_-}$ and $\ssH^{d.1}(\fg_-;\fg)=\operatorname{der}(\fg_-,\fg)_d/\fg_-$. For cocycles in higher degree, we do not have such simple descriptions, but for the calculations in this work we will need expressions only for homological degree $p\leq 2$. Thus we make the $p=2$ case our final example; for $\phi\in\ssC^{d,2}(\fg_-;\fg)=\Hom\bigl(\Wedge^2\fg_-;\fg\bigr)_d$ we have
\begin{align*}
	(\partial \phi)(X,Y,Z) &= \comm{X}{\phi(Y,Z)} - \phi(\comm{X}{Y},Z) + (-1)^{|X|(|Y|+|Z|)} \qty( \comm{Y}{\phi(Z,X)} - \phi(\comm{Y}{Z},X) )\\
		& + (-1)^{(|X|+|Y|)|Z|} \qty( \comm{Z}{\phi(X,Y)} - \phi(\comm{Z}{X},Y) ).
\end{align*}

\subsubsection{Filtered deformations}
\label{sec:filtered-def}

A \emph{filtered deformation} of a $\ZZ$-graded Lie superalgebra $\fg$ is a filtered Lie superalgebra $\fgtilde$ whose associated graded superalgebra $\Gr \fgtilde$ is isomorphic to $\fg$ as a~$\ZZ$-graded Lie superalgebra.\footnote{Filtrations of Lie superalgebras are defined in Appendix~\ref{sec:filtered-superspace-superalg}; see Appendix~\ref{sec:filtered-alg} for more on filtrations and associated graded structures.} We will justify this terminology by showing how the bracket of $\fgtilde$ can be viewed explicitly as a~deformation of that of $\fg$ on the same underlying vector space. In parts of this work, we will also say that a filtered Lie superalgebra $\fgtilde'$ is a \emph{filtered subdeformation} of $\fg$ if there exists an embedding of $\ZZ$-graded superalgebras $\Gr\fgtilde'\hookrightarrow \fg$; that is, if $\Gr\fgtilde'\cong\fg'$ as $\ZZ$-graded superalgebras for some graded subalgebra $\fg'$ of $\fg$, although we often do not specify the subalgebra $\fg'$.

Let us fix a filtered Lie superalgebra $\fgtilde=\fgtilde^{-h}\supseteq \fgtilde^{-h+1} \supseteq \cdots $ of finite depth $h$ and choose, \textit{non-canonically}, some complementary subspaces $W_k$ such that $\fgtilde^k=W_k\oplus\fgtilde^{k+1}$ and $W_k\subseteq\fgtilde_{\overline k}$ (which we can demand due to the superspace filtration compatibility conditions \eqref{eq:filter-super-1} and~\eqref{eq:filter-super-2}). As we have already seen, this induces a graded vector superspace isomorphism~${\Phi\colon\fgtilde=\bigoplus_{k=-h}^\infty W_k\to \Gr\fgtilde}$ where for $X\in W_k$, $\Phi(X)=\overline{X}\in\Gr_k\fgtilde$, where the overline denotes projection $\fgtilde^k\to\Gr_k\fgtilde=\fgtilde^k/\fgtilde^{k+1}$.

Now let $X\in W_k$ and $Y\in W_l$. Then since $\fgtilde^{k+l} = W_{k+l}\oplus\fgtilde^{k+l+1} = W_{k+l}\oplus W_{k+l+1}\oplus\cdots$ and the bracket is filtered, we have
\[
	\comm{X}{Y} = \pr_{k+l}\comm{X}{Y} + \pr_{k+l+1}\comm{X}{Y} + \cdots,
\]
where $\pr_j\colon\fgtilde^{k+l}\to W_j$ are the projection maps. We thus define a \emph{graded} bracket $\comm{-}{-}_0$ on $\fgtilde=\bigoplus_{k=-h}^\infty W_k$ as well as a sequence of maps $\bigl(\mu_d\colon\Wedge^2\fgtilde\to\fgtilde\bigr)_{d> 0}$, where $\mu_d$ has degree $d$, by declaring that, on homogeneous elements $X\in W_k$ and $Y\in W_l$,
\[
	\comm{X}{Y}_0 := \pr_{k+l}\comm{X}{Y},	\qquad \mu_d(X,Y) := \pr_{k+l+d}\comm{X}{Y}
\]
and extending linearly. It can be easily verified that $\comm{-}{-}_0$ is a Lie superalgebra bracket; alternatively, it follows from our next observation. Then on homogeneous elements $X\in W_k$ and~${Y\in W_l}$ again, we have
\[
	\comm{\Phi(X)}{\Phi(Y)} = \comm{\overline{X}}{\overline{Y}}_{\mathrm{Gr}} = \overline{\comm{X}{Y}} = \overline{\comm{X}{Y}_0} = \Phi(\comm{X}{Y}_0),
\]
where in each equality we have applied a definition, and in the last equality we have used the fact that $\comm{X}{Y}_0\in W_{k+l}$ is homogenous. Thus, since homogeneous elements span $\fgtilde$, $\Phi$ is a Lie superalgebra isomorphism\footnote{This shows in particular that although the graded bracket $\comm{-}{-}'_0$ defined by a different choice of complementary subspaces is not equal to $\comm{-}{-}_0$, the graded Lie superalgebra structures defined by the two brackets are isomorphic.}
$\bigl(\fgtilde,\comm{-}{-}_0\bigr)\cong \bigl(\Gr\fgtilde,\comm{-}{-}_{\mathrm{Gr}}\bigr)$.

Let us now change perspective slightly. What we have essentially shown here is that if $\fg=\bigoplus_{k=-h}^\infty\fg_k$ is a fixed $\ZZ$-graded superalgebra, a filtered deformation $\fgtilde$ of $\fg$ can be thought of as an alternative Lie algebra structure on the same underlying vector space as $\fg$; now denoting by $\comm{-}{-}$ the graded bracket on $\fg$ and by $\comm{-}{-}'$ the bracket of $\fgtilde$, for $X,Y\in\fg$ we have
\begin{equation}\label{eq:defining-sequence}
	\comm{X}{Y}' = \comm{X}{Y} + \sum_{d>0}\mu_d(X,Y)
\end{equation}
for some positive-degree maps $\bigl(\mu_d\colon\Wedge^2\fg\to\fg\bigr)_{d> 0}$; that is, the primed bracket can be thought of as a deformation of the unprimed bracket by these maps. We will call this sequence a~\emph{defining sequence} of the filtered deformation. This sequence is not unique; as we have already seen, describing the filtered deformation in this way depends on the choices of complementary subspaces. We also cannot simply define a deformed bracket by choosing some arbitrary sequence~$(\mu_d)_{d>0}$; the Jacobi identity for $\comm{-}{-}'$ imposes relations on the defining sequence. Thus, we need a particular algebraic tool to understand the possible filtered deformations of $\fg$, and that is Spencer cohomology.


\begin{Proposition}\label{prop:def-seq-spencer-cocycle}
	Let $\fg=\bigoplus_{k=-h}^\infty\fg_k$ be a graded Lie superalgebra and $\fgtilde$ a filtered deformation with defining sequence $\bigl(\mu_d\colon\Wedge^2\fg\to\fg\bigr)_{d> 0}$. Then
	\begin{enumerate}\itemsep=0pt
		\item[$(1)$] {\rm\cite[\emph{Proposition}~2.1]{Cheng1998}} The first non-zero term $\mu_k$ in $\qty(\mu_d)_{d> 0}$ is an even Chevalley--Eilenberg $2$-cocycle of $\fg$ with values in the adjoint representation
		$\mu_k\in\ssZ^{k,2}(\fg;\fg)_{\overline{0}}$,	
		where the adjoint Chevalley--Eilenberg complex of $\fg$ is $\ZZ$-graded in a similar fashion to the Spencer complex $($see \eqref{eq:spencer-complex-grading}$)$.
		\item[$(2)$] {\rm\cite[\emph{Proposition}~2.2]{Cheng1998}} Restricting $\mu_k$ to $\fg_-$ gives us a Spencer cocycle
		$\mu_k|_{\Wedge^2\fg_-}\in\ssZ^{k,2}(\fg_-;\fg)$.
		\item[$(3)$] {\rm\cite[\emph{Proposition}~2.2]{Cheng1998}} The cohomology class of this cocycle is $\fg_0$-invariant
			$\qty[\mu_k|_{\Wedge^2\fg_-}]\in \ssH^{k,2}(\fg_-;\fg)^{\fg_0}$.
	\end{enumerate}
\end{Proposition}
We will not prove this or any of the following statements here; we simply note that the above essentially follows by expanding the Jacobi identity for the deformed bracket and treating it degree-by-degree. We will do this for some explicit examples in Section~\ref{sec:filtered-def-poincare}. The results below are more involved but ultimately follow by the same kind of reasoning. We extract from \cite{Cheng1998} only the statements we will need in the rest of this work; see loc.\ cit.\ for full proofs and a more general treatment.

\begin{Definition}\label{def:fund-trans-prolong}
	Let $\fg$ be a graded Lie superalgebra of finite depth $h$. Then $\fg$ is said to be
	\begin{itemize}\itemsep=0pt
	\item \emph{fundamental} if $\fg_-$ is generated by $\fg_{-1}$;
	\item \emph{transitive} if for all $X\in\fg_k$ with $k\geq 0$, $\comm{X}{\fg_-}=0\implies X=0$;
	\item \emph{a full prolongation $\big($of $\bigoplus_{k=-h}^0\fg_k\big)$ of degree $k$} if $\ssH^{d,1}(\fg_-;\fg)=0$ for all $d\geq k$.
	\end{itemize}
\end{Definition}

We will not need fundamentality for the following result, but it is useful to record it here since we will often consider this property alongside the other two defined above.

\begin{Proposition}\label{prop:def-seq-redef}
	Let $\fgtilde$, $\fgtilde'$ be two filtered deformations of $\fg$ with defining sequences
\[
\bigl(\mu_d\colon\Wedge^2\fg\to\fg\bigr)_{d> 0}, \qquad v\bigl(\mu'_d\colon\Wedge^2\fg\to\fg\bigr)_{d> 0},
\]
 respectively. Then
	\begin{enumerate}\itemsep=0pt
		\item[$(1)$] {\rm\cite[\emph{Proposition}~2.3]{Cheng1998}} If for some $k\geq 1$, $\mu_k-\mu_k'|_{\Wedge^2\fg_-}\in\ssB^{k,2}(\fg_-;\fg)$, then $\fgtilde'$ has a defining sequence \smash{$\qty(\mu''_d)_{d> 0}$} such that $\mu''_i=\mu'_i$ for $i<k$ and \smash{$\mu''_k|_{\Wedge^2\fg_-}=\mu_k|_{\Wedge^2\fg_-}$}.
			\label{item:def-seq-redef-1}
		\item[$(2)$] {\rm\cite[\emph{Proposition}~2.5]{Cheng1998}} If $\fg$ is transitive, then the defining sequence $\qty(\mu_d)_{d> 0}$ is completely determined by its restriction to $\fg_-\otimes\fg$.
			\label{item:def-seq-redef-2}
		\item[$(3)$] {\rm\cite[\emph{Proposition}~2.6]{Cheng1998}} If $\fg$ is transitive and for some $k\geq 1$, $\fg$ is an almost full prolongation of degree $k$, $\mu_i|_{\fg_-\otimes\fg}=\mu'_i|_{\fg_-\otimes\fg}$ for $i<k$ and $\mu_k|_{\Wedge^2\fg_-}=\mu'_k|_{\Wedge^2\fg_-}$, then $\fgtilde'$ has a defining sequence \smash{$\qty(\mu''_d)_{d> 0}$} such that $\mu''_i=\mu'_i$ for $i\leq k$ and \smash{$\mu''_i|_{\Wedge^2\fg_-}=\mu_i'|_{\Wedge^2\fg_-}$} for all $i$.
			\label{item:def-seq-redef-3}
	\end{enumerate}
\end{Proposition}

This is a rather technical set of results, but it essentially allows us to redefine the defining sequence of a filtered deformation to partially align it with that of another; under certain conditions, one can completely align the two sequences in order to show that two deformations are isomorphic. We will make use of this in Section \ref{sec:filtered-def-poincare}.

\subsubsection{A note on ``polarised'' and ``depolarised'' expressions}
\label{sec:polarisation}

We finish this section by making note of a useful technique we will often employ in calculations involving Lie superalgebras and elsewhere maps on symmetric products of vector spaces appear.

Let $U$, $W$ be some (finite-dimensional) vector spaces and let $\Psi\in\Odot^n U\to W$ be some completely symmetric multilinear map. Then $\Psi$ is completely determined by ``polarised'' values of the form $\Psi(u,u,\dots,u)$ for elements $u\in U$. The most well-known application of this is the polarisation identity which is used to recover an inner product from its induced norm in the~${n=2}$ case, but it holds more generally. This greatly simplifies calculations involving such maps. When we seek to show that a symmetric tensor vanishes or that two such tensors agree, we will generally do this by showing that the desired property holds in its polarised form, where all inputs are taken to be the same. When we need to ``depolarise'' such an expression, we will indicate that this has been done but will not offer further explanation than what has been indicated here.

\section{Killing spinors and Killing (super)algebras}
\label{sec:killing-spinors-superalgebras}

Let $(M,g)$ be a connected pseudo-Riemannian spin manifold of signature $(p,q)$. We note in particular that $M$ is oriented and denote the special orthonormal frame bundle by $F_{\rm SO}\to M$ and the spin structure by $\varpi\colon P\to F_{\rm SO}$. We need not assume time-orientability,\footnote{Thus $(M,g)$ need not be \emph{strongly spin} in the sense of \cite{Cortes2021,Shahbazi2024_2,Shahbazi2024_1}; see Appendix~\ref{sec:spin-structures} for more details.}
but this assumption may be useful in practice (see Lemma~\ref{lemma:bundle-dirac-current-exist}).
We denote the Lie algebra of Killing vector fields\footnote{Note that this is \emph{not} necessarily the Lie algebra of the isometry group of $(M,g)$ which is in general a subalgebra of $\fiso(M,g)$. We will not need to discuss the isometry group itself so this should not cause confusion.}
by $\fiso(M,g)$ and the Levi-Civita connection by $\nabla$.

Let $S$ be a (possibly $N$-extended) spinor module of $\Spin(p,q)$ and define the \emph{spinor bundle}~${\Sbundle:=P\times_{\Spin(p,q)} S}$; whose sections are the \emph{spinor fields}. We denote the space of spinor fields by $\fS=\Gamma(\Sbundle)$.

We note that the adjoint bundle $\ad F_{\rm SO}\cong\ad P$ can be identified with the bundle of $g$-skew-symmetric endomorphisms of $TM$; we denote its space of sections by $\fso(M,g)$ and recall that it naturally embeds in the Clifford bundle $\Cl(M,g)$ (under natural identification of the latter with the exterior bundle $\Wedge^\bullet TM$) as $\Wedge^2 TM$. As a Lie algebra, the space of sections $\fso(M,g)$ acts on spinor fields via the Clifford action.

See Appendix~\ref{sec:lie-der} for more details of these constructions.

\subsection{Admissible connections and Killing (super)algebras}
\label{sec:admissible-ksa}

Let us now examine the additional structures on the spinor bundle, and the conditions upon them, which are required to define a Killing superalgebra.

\subsubsection{Connections on spinors}

Let $D$ be a connection on the spinor bundle $\Sbundle$. We denote by $\fS_D$ the real vector space of parallel spinors with respect to $D$:
$\fS_D := \qty{\epsilon\in\fS\mid D\epsilon = 0}$.
This is finite-dimensional; indeed, a parallel section can be reconstructed by parallel transport from its value at any point, so in particular we must have $\dim\fS_D\leq\rank\Sbundle=\dim S$. We say that $\fS_D$ is \emph{maximal} if it has the maximal dimension $\dim S$.

Since $\nabla$ and $D$ are both connections on $\Sbundle$, their difference is an $\End\Sbundle$-valued 1-form. We thus define $\beta\in\Omega^1(M;\End\Sbundle)$ by $\beta := \nabla - D$. We now recall that Killing vectors have a natural action on spinors via the \emph{spinorial Lie derivative} (see Appendix~\ref{sec:lie-der-spinors}) and consider whether this action preserves $\fS_D$. For all $X\in\fiso(M,g)$, $Y\in\fX(M)$ and $\epsilon\in\fS_D$, using properties of this Lie derivative (equation~\eqref{eq:lie-der-compat} and a Leibniz rule) for the second equality we have
\begin{align*}
	D_Y\eL_X \epsilon
	&= \nabla_Y\eL_X\epsilon -\beta(Y)\eL_X\epsilon\\
&= \eL_X\nabla_Y\epsilon -\nabla_{\eL_X Y}\epsilon -\eL_X\qty(\beta(Y)\epsilon) +\qty(\eL_X\beta)(Y)\epsilon +\beta(\eL_XY)\epsilon\\
	&= \eL_X D_Y\epsilon - D_{\eL_X Y}\epsilon +\qty(\eL_X\beta)(Y)\epsilon,
\end{align*}
so the spinorial Lie derivative along $X$ preserves $\fS_D$ if and only if $\eL_X\beta$ annihilates $\fS_D$. For the sake of simplicity, and because it will capture the case of supergravity which is our ultimate interest, we will choose to work with Killing vectors which preserve $\beta$; we denote by $\fV_D$ the real subspace of $\fiso(M,g)$ consisting of such vectors
\[
\fV_D := \qty{X\in\fiso(M,g)\mid \eL_X \beta = 0}.
\]
This is in fact an ideal of $\fiso(M,g)$ since it is the annihilator of $\beta$ in the representation of $\fiso(M,g)$ on~$\Omega^1(M;\End\Sbundle)$ defined by the Lie derivative. Clearly, the representation of $\fiso(M,g)$ on~$\fS$ (also via the Lie derivative) restricts to a representation of $\fV_D$ on $\fS_D$.

The curvature of the connection $D$ is the section $R^D \in \Omega^2(M;\End\Sbundle)$ defined by
\begin{gather}\label{eq:D-curvature-def}
	R^D(X,Y)\epsilon = D_XD_Y\epsilon - D_YD_X\epsilon - D_{\comm{X}{Y}}\epsilon,
\end{gather}
where $X,Y\in\fX(M)$ and $\epsilon\in\fS$. One can see immediately that $R^D\epsilon=0$ for all $\epsilon\in\fS_D$, and the former equation is very useful as an integrability condition for the existence of $D$-parallel spinors. In applications, we use the following result to compute $R^D$ and perform manipulations involving it.

\begin{Proposition}\label{prop:D-curvature}
	The curvature $R^D$ is given by the following formula, for $X,Y\in\fX(M)$
	\begin{equation}\label{eq:D-curvature-formula}
		R^D(X,Y)\epsilon = R(X,Y)\cdot\epsilon + \comm{\beta(X)}{\beta(Y)}\epsilon - (\nabla_X\beta)(Y)\epsilon - (\nabla_Y\beta)(X)\epsilon,
	\end{equation}
	where $R$ is the Riemann curvature $($understood as a $2$-form with values in $\ad F_{\rm SO})$ and $[\beta(X),\allowbreak\beta(Y)]$ is a commutator of endomorphisms.
\end{Proposition}

\begin{proof}
	For $\epsilon\in\fS$, we expand $D=\nabla-\beta$ in \eqref{eq:D-curvature-def} to find
	\begin{align*}
		 R^D(X,Y)\epsilon	
			={}& \nabla_X(\nabla_Y\epsilon-\beta(Y)\epsilon) - \beta(X)(\nabla_Y\epsilon-\beta(Y)\epsilon)	\\
				& - \nabla_Y(\nabla_X\epsilon-\beta(X)\epsilon) + \beta(Y)(\nabla_X\epsilon-\beta(X)\epsilon)
				- \nabla_{\comm{X}{Y}}\epsilon + \beta({\comm{X}{Y}})\epsilon	\\
			={}& \qty(\nabla_X\nabla_Y\epsilon - \nabla_Y\nabla_X\epsilon - \nabla_{\comm{X}{Y}}\epsilon)
				+ \qty(\beta(X)\beta(Y)-\beta(Y)\beta(X))\epsilon	\\
				& - \qty(\nabla_X(\beta(Y)\epsilon)-\beta(Y)(\nabla_X\epsilon)) + \qty(\nabla_Y(\beta(X)\epsilon)
				- \beta(X)(\nabla_Y\epsilon)) + \beta({\comm{X}{Y}})\epsilon \\
			={}& R(X,Y)\cdot\epsilon + \comm{\beta(X)}{\beta(Y)} - \nabla_X(\beta(Y))\epsilon + \nabla_Y(\beta(X))\epsilon + \beta({\comm{X}{Y}})\epsilon,
	\end{align*}
	where we have used the definition of the Riemann tensor and the Leibniz rule in the final line. A~further application of the Leibniz rule as well as a the torsion-free property of $\nabla$ gives \eqref{eq:D-curvature-formula}.
\end{proof}

\begin{Remark}\label{rem:riemann-spinor-action}
	It is worth further clarifying the notation $R(X,Y)\cdot\epsilon$ above. As alluded to in the statement of the proposition, the Riemann curvature can be considered as a 2-form with values in $\ad F_{\rm SO}$, meaning that for all $X,Y\in \fX(M)$, $R(X,Y)\in\fso(M,g)$, thus it acts on~$\fS$ via Clifford multiplication. In the proof, we have implicitly used the non-trivial fact that~${R(X,Y)\cdot\epsilon=R^\nabla(X,Y)\epsilon}$ where $R^\nabla\in\Omega^2(M;\End \Sbundle)$ is the curvature 2-form of the spin-lift of the Levi-Civita connection,
\[
R^\nabla(X,Y)\epsilon:=\nabla_X\nabla_Y\epsilon-\nabla_Y\nabla_X\epsilon-\nabla_{\comm{X}{Y}}\epsilon.
\]
Thanks to this relation to the Riemann tensor, we will never have to use this curvature explicitly.
\end{Remark}

\subsubsection{Dirac currents on spinor bundles}
\label{sec:dirac-current-bundle}

In order to define a Killing (super)algebra, we will need to be able to ``square'' (or pair) spinors to tangent vectors, just as we required a Dirac current to define the flat model (super)algebras $\fs$ in Section~\ref{sec:flat-model-algs}. For this purpose, we will define a bundle map $\Otimes^2\Sbundle\to TM$ using a Dirac current on $S$ by exploiting the associated bundle structure of $\Sbundle$ and $TM$ and the holonomy principle.

\begin{Lemma}\label{lemma:bundle-dirac-current-exist}
	Let $V=\RR^{p,q}$ and suppose there exists a Dirac current on the spinor module $S$; that is, an $\fso(p,q)$-equivariant map $\kappa\colon \Odot^2S\to V$ {\rm(}resp.\ $\kappa\colon \Wedge^2S\to V)$. Then there exists a~bundle Dirac current $\underline\kappa\colon \Odot^2\Sbundle\to TM$ {\rm(}resp.\ $\underline\kappa\colon \Wedge^2\Sbundle\to TM)$ which is covariantly constant and which satisfies $\eL_X\underline\kappa=0$ for all $X\in\fiso(M,g)$ if any of the following hold:
	\begin{itemize}\itemsep=0pt
	\item $(M,g)$ has definite signature,
	\item $\kappa$ is $\Spin(p,q)$-equivariant,
	\item $(M,g)$ has indefinite signature and is time-orientable.
	\end{itemize}
\end{Lemma}

\begin{proof}
	We prove the case with $\kappa\colon \Odot^2S\to V$; the other case is completely analogous. Let us consider $\kappa$ as an element $\kappa\in\Odot^2S^*\otimes V$. Fix $x\in M$. We will show that $\kappa$ defines a holonomy-invariant element in the fibre over $x$ of the bundle $\Odot^2\Sbundle^*\otimes TM$ which we note is canonically isomorphic to the associated bundle $E=P\times_{\Spin(p,q)}\bigl(\Odot^2S^*\otimes V\bigr)$. The holonomy principle then guarantees that a covariantly constant section $\underline\kappa$ of that bundle can be defined by parallel transport, which is the desired bundle map.
	
	If $(M,g)$ has definite signature, the spin group is connected, so $\kappa$ is invariant under the action of the spin group, thus the first of our special cases is subsumed by the second. If $(M,g)$ is time-orientable, then the spin structure $\varpi\colon P\to F_{\rm SO}$ is reducible to $\varpi_0\colon P_0\to F_{SO_0}$ (where~$P_0$ is to reduction of $P$ to a $\Spin_0(p,q)$-principal bundle) as discussed in the introduction to Section~\ref{sec:killing-spinors-superalgebras}. In order to treat all three cases in parallel, let us set $G=\Spin(p,q)$ in the first two cases, while in the third case $G=\Spin_0(p,q)$ and $P$ denotes the reduced bundle. In any case, $\kappa$ is $G$-equivariant.
	
	Let us now choose an arbitrary basepoint $u\in P_x$ in the fibre over $x$. This allows us to identify a module of $G$ with the fibre of the associated bundle over $x$; in particular, $\kappa\in\Odot^2S^*\otimes V$ defines an element $\underline\kappa_x:=[u,\kappa]\in E_x$. Moreover, it induces a Lie group isomorphism $\Aut_G(P_x)\cong G$ which allows us to view the holonomy group of the Levi-Civita connection on $P$ at $x$ as a~subgroup of $G$. Thus $G$-invariance of $\kappa$ implies holonomy invariance of $\underline\kappa_x$, and hence the existence of a~covariantly constant global section $\underline\kappa$, as required. We note that $G$-invariance of $\kappa$ also ensures that the definition of $\underline\kappa$ does not depend on the choice of basepoint $u\in P_x$. Finally, $\fg=\fso(p,q)$-invariance of $\kappa$ ensures that $\eL_X\underline\kappa=0$ for all Killing vectors $X$.
\end{proof}

To lighten the notation, we will omit the underline from the bundle map $\kappa$ and trust that this will not cause confusion. With the map $\kappa$ fixed, it will be useful for the proofs of our next few results (and to make a connection to the notation used in Section~\ref{sec:filtered-def-poincare}) to define a section $\gamma$ of the bundle $\Hom\bigl(\Otimes^2\Sbundle,\End(TM)\bigr)$ as follows:
\begin{equation}\label{eq:gamma-geom-def}
	\gamma(\epsilon,\zeta)X := - \kappa(\beta(X)\epsilon,\zeta) - \kappa(\epsilon,\beta(X)\zeta)
\end{equation}
for $\epsilon,\zeta\in\fS$ and $X\in\fX(M)$. Note that $\gamma$ has the same symmetry as $\kappa$. We then have the following.

\begin{Lemma}\label{lemma:dirac-current-endo}
	For all $\epsilon,\zeta\in\fS_D$, $\gamma(\epsilon,\zeta) = A_{\kappa(\epsilon,\zeta)}$.
\end{Lemma}

\begin{proof}
	Let $Z\in\fX(M)$. Then since $\nabla_Z\epsilon = \beta(Z)\epsilon$,
	\begin{gather*}
		A_{\kappa(\epsilon,\zeta)}Z
			= -\nabla_Z\kappa(\epsilon,\zeta)
			= - \kappa(\nabla_Z\epsilon,\zeta) - \kappa(\epsilon,\nabla_Z\zeta)= - \kappa(\beta(Z)\epsilon,\zeta) - \kappa(\epsilon,\beta(Z)\zeta)
			= \gamma(\epsilon,\zeta)
	\end{gather*}
	hence the claim.
\end{proof}

\subsubsection{Existence of the Killing (super)algebra}
\label{sec:exist-ksa}

We now give our first main result followed by a major definition inspired by it.

\begin{Theorem}[existence of (super)algebra associated to $D$]\label{thm:killing-algebra-exist}
	Let $(M,g)$ be a connected pseudo-Riemannian spin manifold, $\Sbundle\to M$ a bundle of spinors with the space of sections $\fS=\Gamma(\Sbundle)$ and let $\kappa\colon \Odot^2\Sbundle\to TM$ {\rm(}resp.\ $\kappa\colon \Wedge^2\Sbundle\to TM)$ be a bundle Dirac current as in Lemma~$\ref{lemma:bundle-dirac-current-exist}$. Let $D$ be a connection on $\Sbundle$ and define $\beta=\nabla-D$, where $\nabla$ is the Levi-Civita spin connection. Define the vector spaces
$\fS_D = \qty{\epsilon\in\fS| D\epsilon = 0}$, $
		\fV_D = \qty{X\in\fiso(M,g)| \eL_X \beta = 0}$.
	and the brackets
$\comm{X}{Y} = \eL_X Y$,
$\comm{X}{\epsilon} = \eL_X \epsilon$,
$\comm{\epsilon}{\zeta} = \kappa(\epsilon,\zeta)$,
	for $X,Y\in\fV_D$, $\epsilon,\zeta\in\fS_D$. Then~${\qty(\fK_D =\fV_D\oplus\fS_D,\comm{-}{-})}$ is a Lie superalgebra $($resp.\ algebra$)$ if and only if the following conditions are satisfied for all $\epsilon,\zeta,\eta\in\fS_D$
	\begin{gather}
		\gamma(\epsilon,\zeta)=A_{\kappa(\epsilon,\zeta)}\in\fso(M,g),\label{eq:killing-alg-exist-1}\\
		\eL_{\kappa(\epsilon,\zeta)}\beta = 0,\label{eq:killing-alg-exist-2}\\
		\beta(\kappa(\epsilon,\zeta))\eta + \beta(\kappa(\zeta,\eta))\epsilon + \beta(\kappa(\eta,\epsilon))\zeta
		+ \gamma(\epsilon,\zeta)\cdot\eta + \gamma(\zeta,\eta)\cdot\epsilon +\gamma(\eta,\epsilon)\cdot\zeta
		= 0 \label{eq:killing-alg-exist-3},
	\end{gather}
	where $\gamma$ is the map defined by equation~\eqref{eq:gamma-geom-def}.
\end{Theorem}

\begin{proof}
	We first check closure and then the Jacobi identities. We have $\comm{\fV_D}{\fV_D}\subseteq\fV_D$ and $\comm{\fV_D}{\fS_D}\subseteq \fS_D$ by construction, while $\comm{\fS_D}{\fS_D}\subseteq\fV_D$ if and only if
	\begin{equation}\label{eq:killing-closure-cond}
		\eL_{\kappa(\epsilon,\zeta)}g = 0
		\qquad\text{and}\qquad
		\eL_{\kappa(\epsilon,\zeta)}\beta = 0 \qquad \forall\epsilon,\zeta\in\fS_D.
	\end{equation}
	By equation~\eqref{eq:killing-condition}, $\eL_{\kappa(\epsilon,\zeta)}g = 0$ if and only if $A_{\kappa(\epsilon,\zeta)}\in\fso(M,g)$. By Lemma~\ref{lemma:dirac-current-endo}, this is equivalent to condition~\eqref{eq:killing-alg-exist-1}. The latter condition in \eqref{eq:killing-closure-cond} is just \eqref{eq:killing-alg-exist-2}.
	
	The Jacobi identity for three vectors is clearly satisfied. For the others, we denote by $[-,-,-]$ the Jacobiator and find
	\[
		[X,Y,\epsilon]= 0 \iff \eL_{\comm{X}{Y}} \epsilon = \eL_X \eL_Y \epsilon - \eL_Y \eL_X \epsilon,
	\]
	which simply says that $\fS_D$ is a representation of $\fV_D$, hence is true by construction;
	\[
		[X,\epsilon,\zeta]= 0 \iff \eL_X \kappa(\epsilon,\zeta) = \kappa(\eL_X\epsilon,\zeta) + \kappa(\epsilon,\eL_X\zeta),
	\]
	which is satisfied by the last part of Lemma~\ref{lemma:bundle-dirac-current-exist};
	\[
		[\epsilon,\zeta,\eta] = 0 \iff \eL_{\kappa(\epsilon,\zeta)} \eta + \eL_{\kappa(\eta,\epsilon)} \zeta + \eL_{\kappa(\zeta,\eta)} \epsilon = 0,
	\]
	which, using the equation $\nabla \epsilon=\beta\epsilon$ and Lemma~\ref{lemma:dirac-current-endo}, is condition~\eqref{eq:killing-alg-exist-3}.
\end{proof}

We interpret this result as giving us a set of conditions on the spinor connection $D$ for the existence of an associated (super)algebra. This inspires the following definition. Note that conditions (1) and (2) are stronger than \eqref{eq:killing-alg-exist-1} and \eqref{eq:killing-alg-exist-3} above; we demand that they hold for \emph{all} spinors, not just $D$-parallel spinors. This makes the definition simpler to check, and it will also be useful for some later results.

\begin{Definition}[admissible connections, Killing spinors and Killing (super)algebras]\label{def:killing-spinor}
	Using the notation of the theorem above, a connection $D=\nabla-\beta$ on $\Sbundle$ is \emph{admissible} if the following hold:
	\begin{enumerate}\itemsep=0pt
		\item The map $\gamma$ defined by equation~\eqref{eq:gamma-geom-def} takes values in $\ad F_{\rm SO}$ (that is, $\gamma(\epsilon,\zeta)\in\fso(M,g)$ for all $\epsilon,\zeta\in\fS$); 

		\item $
\beta(\kappa(\epsilon,\zeta))\eta + \beta(\kappa(\zeta,\eta))\epsilon + \beta(\kappa(\eta,\epsilon))\zeta
				+ \gamma(\epsilon,\zeta)\cdot\eta + \gamma(\zeta,\eta)\cdot\epsilon +\gamma(\eta,\epsilon)\cdot\zeta
				= 0$ for all $\epsilon,\zeta,\eta\in\fS$; 

		\item $\eL_{\kappa(\epsilon,\zeta)}\beta = 0$ for all $\epsilon,\zeta\in\fS_D$. 
	\end{enumerate}
	If $D$ is admissible, the differential equation $D\epsilon=0$ (equivalently $\nabla\epsilon =\beta\epsilon$) is called the \emph{Killing spinor equation} and $\fS_D$ the space of \emph{Killing spinors}. Furthermore, $\fV_D$ is the space of \emph{restricted Killing vectors}, and $\fK_D$ is the \emph{Killing $($super$)$algebra}.\footnote{Our use of the term ``Killing spinor'' here is a generalisation of the sense in which it is used in the supergravity literature and in particular those discussed in \cite{Beckett2021,deMedeiros2016,deMedeiros2018,Figueroa-OFarrill2017_1,Figueroa-OFarrill2017}. It does not subsume all geometric Killing spinors or generalisations thereof; it does however subsume the particular examples in \cite{Figueroa-OFarrill2008_except}.}
\end{Definition}

Let us make some comments to connect the above with treatments of Killing spinors and Killing superalgebras in the supergravity literature. In a supergravity theory, the Killing spinor equation $\nabla\epsilon=\beta\epsilon$ arises by demanding that the supersymmetry variation (with parameter $\epsilon$) of the gravitino vanishes, with the 1-form $\beta$ being parametrised by some bosonic fields $\Phi$. To prove the existence of a Killing superalgebra, one essentially checks that conditions (1) and (2) of the above definition hold and that the Killing spinor equation and bosonic equations of motion imply that $\eL_{\kappa_\epsilon}\Phi=0$, hence $\eL_{\kappa_\epsilon}\beta=0$, where $\kappa_\epsilon$ is the square of a Killing spinor $\epsilon$, which is condition~(3) above.

Often in the supergravity literature, the Killing superalgebra is considered to contain not all of the Killing vectors $X\in\fV_D$, but only those for which $\eL_X\Phi=0$ (hence $\eL_X\beta=0$), giving a subalgebra of what we call $\fK_D$ here in general. Some sources (e.g., \cite{Figueroa-OFarrill2009,Hustler2016}) call this the \emph{symmetry superalgebra}, giving an even more restrictive definition for the Killing superalgebra~-- namely that it contains only those Killing vectors in the span of those generated by pairs of spinors~-- which corresponds to our notion of the \emph{Killing ideal}. Our terminology here agrees with that of~\cite{Figueroa-OFarrill2017_1}. The Killing ideal is referred to as the \emph{transvection ideal} in \cite{Santi2022}.

\begin{Definition}[Killing ideal]
	Let $D$ be an admissible connection and $\fK_D$ its Killing (super)algebra. The \emph{Killing ideal} is the ideal $\overline\fK_D$ of $\fK_D$ generated by the Killing spinors; that is, $\overline\fK_D=\overline\fV_D\oplus\fS_D$ where $\overline\fV_D = \comm{\fS_D}{\fS_D}.$
\end{Definition}

\begin{Remark}
	The conditions \eqref{eq:killing-alg-exist-1}--\eqref{eq:killing-alg-exist-3} (and the related conditions in Definition~\ref{def:killing-spinor}) for existence of a Killing (super)algebra actually concern only the Killing ideal; that is, given a (not necessarily admissible) spinor connection $D$, $\fK_D$ is a Lie (super)algebra if and only if $\overline\fK_D$ is a~Lie (super)algebra. Moreover, we will see later that conditions~(1) and (2) of Definition~\ref{def:killing-spinor} are essentially Spencer cocycle conditions.
\end{Remark}

A similar general formulation of Killing superalgebras to the one given here was given in \cite{Hustler2016}, however our treatment here is more general because we allow for arbitrary metric signature and Dirac current $\kappa$ of arbitrary symmetry, and that work also makes the assumption mentioned above that $\beta$ is parametrised by some bosonic fields $\Phi$ and that the Killing vectors in the algebra satisfy $\eL_X\Phi=0$.

\subsection{Algebraic structure of Killing (super)algebras}
\label{sec:alg-structure-killing}

We now turn to describing the algebraic structure of the (super)algebras $\fK_D$ associated to the admissible connections $D$ we have just defined. Our treatment here is a direct generalisation of that of \cite{Figueroa-OFarrill2017_1}, which considered the case of Killing superalgebras of highly supersymmetric backgrounds of 11-dimensional supergravity. The methods used here were introduced (for Killing vectors) by Kostant \cite{Kostant1955} and independently by Geroch \cite{Geroch1969}, the latter of whom introduced the term \emph{Killing transport}.

\subsubsection{Killing transport}

Let $D$ be an admissible connection on the spinor bundle $\Sbundle$ and recall that the space of Killing spinors ($D$-parallel sections of $\Sbundle$) is denoted $\fS_D$. Let us define the vector bundle
$\eE := TM\oplus\Sbundle\oplus\ad F_{\rm SO}$,
where we recall that the space of sections of $\ad F_{\rm SO}$, the adjoint bundle to the special orthonormal frame bundle, can be identified with the Lie algebra $\fso(M,g)$ of skew-symmetric endomorphisms of the tangent bundle. We define a connection $\eD$ on $\eE$ as follows:
\[
	\eD_Y(X,\epsilon,A) := (\nabla_Y X + AY,D_Y\epsilon,\nabla_Y A + R(Y,X))
\]
for $X,Y\in\fX(M)$, $\epsilon\in\fS$, $A\in\fso(M,g)$.

\begin{Proposition}\label{prop:killing-transport}
	If $D$ is an admissible connection, the parallel sections of $\eD$ are triples $(X,\epsilon,A_X)$, where $X$ is a Killing vector, $\epsilon$ a Killing spinor and $A_X=-\nabla X$. Thus the $\RR$-vector space of $\eD$-parallel sections is isomorphic to $\fiso(M,g)\oplus\fS_D$.
\end{Proposition}
\begin{proof}
If the section $(X,\epsilon,A)$ of $\eE$ is parallel with respect to $\eD$, clearly $\epsilon$ is a Killing spinor and we have $AY = -\nabla_Y X$ for all $Y\in\fX(M)$, so $A = -\nabla X = A_X$. Since $A\in\fso(M,g)$, $X$ must be a Killing vector. Conversely, if $X$ is a Killing vector and $\epsilon$ a Killing spinor, we have $\nabla_Y A_X = -R(Y,X)$ (by equation~\eqref{eq:nabla-A_X}), and so $(X,\epsilon,A_X)$ is parallel with respect to $\eD$.
\end{proof}

Thanks to this result and since $M$ is connected, if one knows the values $(X,\epsilon,A_X)_p$ at any point~${p\in M}$ for any pair $(X,\epsilon)$ with $X\in\fiso(M,g)$ and $\epsilon\in\fS_D$ -- which we call a \emph{Killing pair} -- one can recover the section $(X,\epsilon,A_X)$ by parallel transport with respect to $\eD$. For this reason, $\eD$ is referred to as the \emph{Killing transport connection}, and $(X,\epsilon,A_X)_p$ the \emph{Killing transport data} of~$(X,\epsilon)$ at~$p$. We define an injective $\RR$-linear map $\Phi$ sending a Killing pair in~$\fK_D$ (so that~${X\in\fV_D}$) to its Killing transport data at $p$ as follows:
$\Phi_p\colon \fK_D\to\eE_p$,
	$\Phi_p(X,\epsilon):=(X_p,\epsilon_p,(A_X)_p)$.

\subsubsection{Localising the Killing algebra at a point}

If we restrict our attention to Killing pairs $(X,\epsilon)\in\fK_D$, it is possible to recover the Killing transport data of $\comm{(X,\epsilon)}{(Y,\zeta)}$ at a point from that of $(X,\epsilon)$ and $(Y,\zeta)$. We have
\[
	\comm{(X,\epsilon)}{(Y,\zeta)} = \qty(\comm{X}{Y} + \kappa(\epsilon,\zeta),\eL_X\zeta - \eL_Y\epsilon),
\]
and the corresponding parallel section of $\eE$ is
\begin{equation}\label{eq:killing-comm-triple}
	(\comm{X}{Y} + \kappa(\epsilon,\zeta),\eL_X\zeta - \eL_Y\epsilon,A_{\comm{X}{Y}} + A_{\kappa(\epsilon,\zeta)}).
\end{equation}
We now expand some of the terms in this expression. Equations~\eqref{eq:comm-formula} and \eqref{eq:A-comm} from Appendix~\ref{sec:lei-der-vector-tensor} give us
\[
	\comm{X}{Y} = \nabla_X Y - \nabla_Y X = A_X Y - A_Y X,\qquad
	A_{\comm{X}{Y}} = \comm{A_X}{A_Y} - R(X,Y).
\]
The Killing spinor equation in the form $\nabla \epsilon = \beta\epsilon$ gives
\[
	\eL_Y\epsilon = \nabla_Y\epsilon + A_Y\cdot\epsilon = \beta(Y)\epsilon + A_Y\cdot\epsilon,
\]
and by Lemma~\ref{lemma:dirac-current-endo},
$
	A_{\kappa(\epsilon,\zeta)} = - \kappa(\beta\epsilon,\zeta) - \kappa(\epsilon,\beta\zeta)$.
Thus the triple \eqref{eq:killing-comm-triple} is given by
\begin{gather*}
	\comm{X}{Y} + \kappa(\epsilon,\zeta) = A_X Y - A_Y X + \kappa(\epsilon,\zeta),\\
	\eL_X\zeta - \eL_Y\epsilon = \beta(X)\zeta - \beta(Y)\epsilon + A_X\cdot\zeta - A_Y\cdot\epsilon,\\
	A_{\comm{X}{Y}} + A_{\kappa(\epsilon,\zeta)} = \comm{A_X}{A_Y} - R(X,Y) - \kappa(\beta\epsilon,\zeta) - \kappa(\epsilon,\beta\zeta).
\end{gather*}
We now evaluate the expressions above at $p$ to obtain the Killing transport data of $\comm{(X,\epsilon)}{(Y,\zeta)}$. Let us define the inner product space $(V,\eta):=(T_pM,g_p)$ and identify $S:=\Sbundle_p$ as a ($N$-extended) spinor module of $\Spin(V)\subseteq\Cl(V)$ (where we omit $\eta$ as usual); we then have $\eE_p = V\oplus S\oplus\fso(V)$ as a $\Spin(V)$-module.\footnote{We note that there is some abuse of notation here; we previously denoted the spinor module of $\Spin(p,q)$ used to define $\Sbundle$ by~$S$; identification of these two $S$'s is not canonical and requires e.g., a choice of basepoint in the fibre of the spin structure over~$p$.}
For the sake of readability, we consider the data for the brackets~$\comm{X}{Y}$, $\comm{X}{\epsilon},$ and $\comm{\epsilon}{\zeta}$ separately; we find
\begin{align}
	\Phi_p(\comm{X}{Y})
		&= ((A_X)_p Y_p - (A_Y)_p X_p,0,\comm{(A_X)_p}{(A_Y)_p} - R_p(X_p,Y_p)),
		\label{eq:phi-p-XY}\\
	\Phi_p(\comm{X}{\epsilon})
		&= (0,\beta_p(X_p)\epsilon_p + (A_X)_p\cdot\epsilon_p,0),
		\label{eq:phi-p-Xepsilon}\\
	\Phi_p(\comm{\epsilon}{\zeta})
		&= (\kappa_p(\epsilon_p,\zeta_p),0,- \kappa_p(\beta_p\epsilon_p,\zeta_p) - \kappa_p(\epsilon_p,\beta_p\zeta_p)),
		\label{eq:phi-p-epsilonzeta}
\end{align}
where $\kappa_p\in\Otimes^2S^*\otimes V$, $\beta_p\in V^*\otimes\End S$, and $R_p\in\Wedge^2V^*\otimes \fso(V)$ are the values of $\kappa$, $\beta,$ and $R$ at $p\in M$. Thus we can indeed express the Killing transport data of any $\comm{(X,\epsilon)}{(Y,\zeta)}$ at $p\in M$ in terms of the Killing transport data of $(X,\epsilon)$ and $(Y,\zeta)$ at $p$ if we know the ``background'' data~$(\kappa_p,\beta_p,R_p)$.

\subsubsection{The structure theorem}

We will now prove our second main result. This is a generalisation of the first part of \cite[Theorem~12]{Figueroa-OFarrill2017_1}, which is particular to 11-dimensional supergravity. Recall from Definition~\ref{def:flat-model-alg} that for a fixed inner product space $(V,\eta)$, a flat model (super)algebra $\fs$ is defined by a spinor module $S$ and a Dirac current $\kappa$ on $S$. See Section~\ref{sec:filtered-def} for background on filtered deformations and Appendices~\ref{sec:filtered-alg} and~\ref{sec:filtered-superspace-superalg} for more on filtrations.

\begin{Theorem}[structure of Killing (super)algebras]\label{thm:killing-algebra-filtered}
	Let $D$ be an admissible connection $($see Definition~{\rm\ref{def:killing-spinor})} on a spinor bundle $\Sbundle$ with Dirac current $\kappa$ over a connected pseudo-Riemannian spin manifold $(M,g)$. Then the Killing $($super$)$algebra $\fK_D$ is a filtered deformation of a graded subalgebra $\fa$ of the flat model $($super$)$algebra $\fs$ with odd part $S$ and Dirac current $\kappa$.
\end{Theorem}

\begin{Remark}
Before proving this result, we note that there is a more na\"{\i}ve approach to it which unfortunately fails. We have already observed that $\eE_p = V\oplus S\oplus\fso(V)$, so one might be tempted to identify $\eE_p$ with the flat model algebra $\fs$ and proceed as follows: we have a~map~${\Phi_p\colon \fK_D\to \eE_p=\fs}$, and we can therefore compute brackets on $\Im\Phi_p$,
\begin{gather*}
	\comm{\Phi_p(X)}{\Phi_p(Y)}
		= \qty((A_X)_p Y_p - (A_Y)_p X_p,0,\comm{(A_X)_p}{(A_Y)_p}),\\
	\comm{\Phi_p(X)}{\Phi_p(\epsilon)}
		= \qty(0,(A_X)_p\cdot\epsilon_p,0),\qquad
	\comm{\Phi_p(\epsilon)}{\Phi_p(\zeta)}
		= \qty(\kappa_p(\epsilon_p,\zeta_p),0,0).
\end{gather*}
Comparing to equations~\eqref{eq:phi-p-XY}--\eqref{eq:phi-p-epsilonzeta}, we see that $\Phi_p$ is not a Lie algebra homomorphism in general, but perhaps it is the natural map carrying a filtered Lie algebra to its associated graded inside $\fs$.

The problem is that $\Im\Phi_p$ need not be closed under the bracket $\comm{-}{-}$ -- the expressions above need not be the transport data of any Killing pairs -- and it may not even be a graded subspace of $\eE_p$, since for example Killing vectors embed diagonally; $X\mapsto(X_p,0,(A_X)_p)$. A~more careful method is therefore required.
\end{Remark}

\begin{proof}[Proof of Theorem~\ref{thm:killing-algebra-filtered}]	
	We fix an arbitrary point $p\in M$ and set $V=T_pM$, $S=\Sbundle_p$ as before. Define the evaluation maps
\smash{$\evaluate^V_p\colon \fV_D \to V$}, $X\mapsto X_p$,
\smash{$\evaluate^S_p\colon \fS_D \to S$}, $\epsilon\mapsto \epsilon_p$.
	Let $V'$ be the subspace of $V$ consisting of values of restricted Killing vectors
$V' = \Im\evaluate^V_p = \qty{X_p\in V \mid X\in\fV_D }$,
	and let $\fh$ be the subspace of $\fso(V)$ consisting of values at $p$ of endomorphisms $A_X$ for Killing vectors which vanish at $X$
	\[
		\fh = \qty{(A_X)_p\in\fso(V) \mid X\in\fV_D\colon X_p=0}.
	\]
	Since a Killing vector $X$ is determined by its Killing transport data $(X_p,(A_X)_p)$, there is an isomorphism of vector spaces $\fh\simeq \ker\evaluate^V_p$. For $A\in\fh$, we denote by $X_A$ the corresponding Killing vector in $\ker\evaluate^V_p$. We also define
$S' = \Im\evaluate^S_p = \qty{\epsilon_p\in S \mid \epsilon\in\fS_D}
$
	and note that there exists a vector space isomorphism $S'\simeq \fS_D$ since a Killing spinor is defined by its value at $p$. For ease of notation, when there is no ambiguity we now identify these spaces; the symbol $\epsilon$ will denote both a Killing spinor and its value at $p$.
	
	Now let $\fa = V'\oplus S'\oplus\fh \subseteq \fs$. We claim that this graded subspace is in fact a graded subalgebra of $\fs$. If $A,B\in\fh$, by equation~\eqref{eq:phi-p-XY} we have
	\begin{gather*}
		\comm{X_A}{X_B}_p = 0,\qquad
		(A_{\comm{X_A}{X_B}})_p = \comm{A_{X_A}}{A_{X_B}}_p = \comm{A}{B},
	\end{gather*}
	so $\comm{X_A}{X_B}\in\ker\evaluate^V_p$, and the corresponding element in $\fh$ is $\comm{A}{B}$, showing that $\fh$ is closed under the bracket. Note that this also shows that $\fh\cong\ker\evaluate^V_p$ as Lie algebras. Now let $A\in\fh$ and $v\in V'$. Then there is some $Y\in\fV_D$ with $Y_p=v$. Again by equation~\eqref{eq:phi-p-XY}, we have
	\[
		\comm{X_A}{Y}_p = \qty(A_{X_A})_p Y_p = A v = \comm{A}{v},
	\]
	so $\comm{A}{v}$ is the value at $p$ of $\comm{X_A}{Y}$. Thus $\comm{\fh}{V'}\subseteq V'$. For $A\in\fh$ and $\epsilon\in S'\simeq\fS_D$, by equation~\eqref{eq:phi-p-Xepsilon} we have
	\[
		\comm{X_A}{\epsilon}_p = (A_{X_A})_p\cdot\epsilon_p = A\cdot\epsilon = \comm{A}{\epsilon},
	\]
	so $\comm{A}{\epsilon}\in S'$. Thus $\comm{\fh}{S'}\subseteq S'$. Finally, for $\epsilon,\zeta\in S$, we have
$\comm{\epsilon}{\zeta}_p
			= \kappa_p(\epsilon_p,\zeta_p)
			= \comm{\epsilon_p}{\zeta_p}$,
	so $\comm{S'}{S'}\subseteq V'$. Hence $\fa$ is indeed a graded Lie subalgebra of $\fs$.
	
	Using the notation and terminology of Appendices~\ref{sec:filtered-alg} and~\ref{sec:filtered-superspace-superalg}, we can equip $\fK_D$ with a~filtration $\fK_D^\bullet$ of depth 2 ($\fK_D^i = \fK_D$ for $i\leq-2$) as follows: we set $\fK_D^i = 0$ for $i>0$ and
	\[
		\fK_D^{-2}=\fK_D = \fV_D\oplus\fS_D \quad \supset \quad
		\fK_D^{-1} = \ker\evaluate^V_p\oplus \fS_D \quad \supset \quad
		\fK_D^0 = \ker\evaluate^V_p\cong \fh.
	\]
	The calculations above show that, as $\fh\cong\ker\evaluate^V_p$-modules,
	\begin{gather*}
		\Gr_{-2}\fK_D^\bullet = \faktor{\fK_D}{\ker\evaluate^V_p\oplus\fS_D} \cong V' \qquad \text{\big(using $\evaluate^V_p$\big)},	\\
		\Gr_{-1}\fK_D^\bullet = \faktor{\ker\evaluate^V_p\oplus\fS_D}{\ker\evaluate^V_p} \cong S'	\qquad \text{\big(using $\evaluate^{S}_p$\big)},	\\
		\Gr_0 \fK_D^\bullet = \ker\evaluate^V_p \cong \fh \qquad\text{(using $X\mapsto (A_X)_p$)},
	\end{gather*}
	thus $\Gr\fK_D^\bullet\cong\fa$ as $\fh$-modules. We will show that this isomorphism preserves Lie brackets, hence~$\fK_D$ is a filtered deformation of $\fa$. We denote the associated graded bracket on $\Gr\fK_D^\bullet$ by~$\comm{-}{-}_{\mathrm{Gr}}$ and use an overline to denote projections $\fK_D^i\to\Gr\fK_D^i$. For $A,B\in\fh$, $X,Y\in\fV_D$ and $\epsilon,\zeta\in\fS_D$, we have the following for our map $\Gr\fK_D^\bullet\to\fa$:
	\begin{gather*}
		\comm{X_A}{\overline{Y}}_{\mathrm{Gr}} = \overline{\comm{X_A}{Y}} = \overline{A_{X_A}Y}
			 \longmapsto 	 A Y_p = \comm{A}{Y_p}	\qquad (\text{in degree $-2$}),	\\
		\comm{\epsilon}{\zeta}_{\mathrm{Gr}} = \overline{\comm{\epsilon}{\zeta}} = \overline{\kappa(\epsilon,\zeta)}
			 \longmapsto 	 \kappa_p(\epsilon,\zeta) = \comm{\epsilon}{\zeta}	\qquad (\text{in degree $-2$}),	\\
		\comm{X_A}{\epsilon}_{\mathrm{Gr}} = \overline{\comm{X_A}{\epsilon}} = \overline{A_{X_A}\cdot\epsilon}
			 \longmapsto 	 A\cdot \epsilon = \comm{A}{\epsilon}	\qquad (\text{in degree $-1$}),	\\
		\comm{X_A}{X_B}_{\mathrm{Gr}} = \overline{\comm{X_A}{X_B}} = X_{\comm{A}{B}}
		 \longmapsto 	 \comm{A}{B}	\qquad (\text{in degree $0$}).
	\end{gather*}
	The only brackets we have not considered above are $\big[\overline{X},\epsilon\big]_{\mathrm{Gr}}$ and $\big[\overline{X},\overline{Y}\big]_{\mathrm{Gr}}$, but these lie in $\Gr_{-3}\fK_D^\bullet$ and $\Gr_{-4}\fK_D^\bullet$ respectively, both of which are both zero by definition.
	\end{proof}
	
We recall from the discussion in Section~\ref{sec:filtered-def} that filtered deformations can be explicitly described as a deformation of the bracket of the graded model algebra by a sequence of maps of positive degrees. It will be useful to give our Killing superalgebras such a presentation.

First note that the isomorphism $\fh\cong\ker\evaluate^V_p$ and the evaluation map $\evaluate^V_p$ itself give us the short exact sequence of vector spaces
\[
\begin{tikzcd}
	0 \ar[r] & \fh \ar[r] & \fV_D \ar[r] & V' \ar[r] & 0.
\end{tikzcd}
\]
Any short exact sequence of vector spaces splits, so there exists a (linear) splitting map $\Lambda\colon V'\to \fV_D$ with $\evaluate^V_p(\Lambda(v))=\Lambda(v)_p=v$. We can also parametrise the splitting in terms of a~map $\lambda\colon V'\to \fso(V)$ by setting $\lambda(v)=(A_{\Lambda(v)})_p$; the vector field $\Lambda(v)$ can be recovered from its Killing transport data $(v,\lambda(v))$. The splitting induces a linear isomorphism $V'\oplus\fh\to\fV_D$ given by $(v,A)\mapsto \Lambda(v)+X_A$.
This extends to a $\ZZ_2$-graded linear isomorphism $\Psi\colon \fa \to \fK_D$ given by~${\Psi(v,\epsilon,A) = (\Lambda(v)+X_A,\epsilon)}$,
where $v\in V'$, $A\in\fh$, and $\epsilon\in S'\simeq\fS_D$.
We use $\Psi$ to define a~new bracket $\comm{-}{-}'$ on $\fa$ by pulling back the bracket on $\fK_D$
\[
	\comm{(v,\epsilon,A)}{(w,\zeta,B)}'
		:= \Psi^{-1}(\comm{\Psi(v,\epsilon,A)}{\Psi(w,\zeta,B)})
\]
with respect to which $\Psi$ is a Lie (super)algebra isomorphism $(\fa,\comm{-}{-}')\cong\fK_D$ by construction. Explicitly, the bracket is given by
\begin{gather}
	\comm{A}{B}' = \comm{A}{B},\qquad
	 \comm{A}{v}' = Av + \comm{A}{\lambda(v)}-\lambda(Av),\nonumber	\\
	\comm{v}{w}' = \lambda(v)w-\lambda(w)v + \theta_\lambda(v,w),\qquad
	\comm{A}{\epsilon}' = A\cdot\epsilon,\nonumber	\\
	 \comm{v}{\epsilon}' = \beta_p(v)\epsilon + \lambda(v)\cdot \epsilon,	\qquad
	\comm{\epsilon}{\zeta}'
		 = \kappa_p(\epsilon,\zeta) + \gamma_p(\epsilon,\zeta)-\lambda\qty(\kappa_p(\epsilon,\zeta)),\label{eq:KSA-deformed-brackets}
\end{gather}
where $\gamma_p\colon\Otimes^2S\to\fso(V)$ is the evaluation at $p$ of the map $\gamma$ defined by equation~\eqref{eq:gamma-geom-def} and $\theta_\lambda:\Wedge^2V'\to \fh$ is the map
\[
	\theta_\lambda(v,w)
		= \comm{\lambda(v)}{\lambda(w)} - R_p(v,w) - \lambda(\lambda(v)w-\lambda(w)v).
\]
Using the notation of Lie algebra cohomology, we have a map $\partial\lambda\colon\Wedge^2\fs\to\fs$ of degree $+2$ defined as follows
\begin{gather*}
	\partial\lambda(A,v) = \comm{A}{\lambda(v)}-\lambda(Av),\qquad
	 \partial\lambda(v,w) = \lambda(v)w-\lambda(w)v,\qquad
	\partial\lambda(v,\epsilon)	= \lambda(v)\cdot\epsilon,\\
	 \partial\lambda(\epsilon,\zeta)	= -\lambda(\kappa_p(\epsilon,\zeta)),
\end{gather*}
whence we have
\[
	\comm{-}{-}' = \comm{-}{-} + \beta_p + \gamma_p + \partial\lambda + \theta_\lambda,
\]
where $\beta_p$, $\gamma_p$ and $\theta_\lambda$ have been trivially extended to maps $\Wedge^2\fs\to\fs$. The deforming maps are $\mu = \beta_p + \gamma_p+\partial\lambda$, of degree $+2$, and $\theta_\lambda$, of degree $+4$. That is, in the language of Section~\ref{sec:filtered-def} (see equation \eqref{eq:defining-sequence}), the deformation has defining sequence $(\mu,\theta_\lambda,0,\dots)$.

The splitting $\Lambda\colon V'\to \fV_D$ and the corresponding map $\lambda \colon V'\to\fh$ can be viewed as a ``correction'' to the failed na\"{\i}ve approach discussed in the remark before the proof of Theorem~\ref{thm:killing-algebra-filtered}. Note that these maps are not canonical; if $\Lambda'\colon V'\to\fV_D$ is another splitting and $\lambda'\colon V'\to\fso(V)$ is the corresponding map, we have
\[
	\Phi_p(\Lambda'(v)-\Lambda(v)) = (v-v,\lambda'(v)-\lambda(v)) = (0,\lambda'(v)-\lambda(v)),
\]
so $(\Lambda'-\Lambda)(v)\in\ker\evaluate^V_p$, and $(\lambda'-\lambda)(v)$ is the corresponding element of $\fh$. Thus $\lambda$, and therefore the splitting, is unique up to a choice of map $V'\to\fh$. We will see later that such a~map is a Spencer $(2,1)$-cocycle for $\fa$, so it is consistent that changing the defining sequence by the coboundary of such a map does not change the isomorphism class of the deformation, only its presentation. Also note that this choice of splitting is exactly a choice of complementary submodules as discussed near the beginning of Section~\ref{sec:filtered-def}.

\subsection{Constrained Killing spinors}
\label{sec:alg-kse}

So far in this section, we have only considered the differential Killing spinor equation ${D\epsilon=0}$ which in supergravity theories arises from demanding that the variation of the gravitino ${\partial_\epsilon\Psi=\!D\epsilon}$ vanishes. While in some pure minimal supergravity theories (in particular, those in dimensions~4,~5,~6 and~$11$, which have been studied using Spencer cohomology \cite{Beckett2021,deMedeiros2016,deMedeiros2018,Figueroa-OFarrill2016,Figueroa-OFarrill2017}) this is the only local condition imposed by supersymmetry, in most supergravity theories (including any extended, gauged or matter-coupled theories and even some pure minimal theories), there are additional fermions present. A background of such a theory is called supersymmetric if the variations of \emph{all} fermions vanish for some supersymmetry parameter. The variations of the additional fermions are always algebraic, rather than differential, in the spinor. The vanishing of these variations are referred to as \emph{algebraic} or \emph{linear constraints} and we will call the spinors satisfying them as well as $D\epsilon=0$ \emph{constrained Killing spinors}.\footnote{Our terminology follows \cite{Cortes2021,Lazaroiu2016} (omitting the word ``generalised''). The more recent work \cite{Shahbazi2024_1} uses the term ``constrained differential spinor''.}
Similar considerations for Killing superalgebras are found in \cite{Hustler2016}.

We will make some comments here about how these might be brought into our present framework; we leave a full treatment for future work.

\subsubsection{Linear constraints}

Schematically, the linear constraints are of the form $\eP\epsilon=0$ where $\eP$ is some fibrewise linear operator on the spinor bundle $\Sbundle$. More concretely, all fermions may take values in the same representation as the gravitini, in which case $\eP\in\End\Sbundle$. In a more general situation, the additional fermions take values in some vector bundle $E\to M$ (for example, an associated bundle for some gauge symmetry). Then the constraint equation is $\eP\epsilon=0$, where
\[
\eP\in\Hom(\Sbundle,\Sbundle\otimes E)\cong \End\Sbundle\otimes E.
\]
This bundle map is parametrised by the bosonic fields of the background, just as~${\beta=D-\nabla}$~is.\footnote{Let us note that the fermion variations $\eP\epsilon$ and the expression $\beta\epsilon$ are fibrewise linear in $\epsilon$ but not necessarily in the bosonic background fields. For example, gauge field strengths (hence derivatives of gauge connections) and derivatives of scalars typically appear as coefficients in both, and may also appear in nonlinear combinations.}

Let $D=\nabla-\beta$ be a spinor connection which is not necessarily admissible in the sense of Definition~\ref{def:killing-spinor} and let us denote the space of $D$-parallel sections of $\Sbundle$ which are annihilated by~$\eP$~by
\[
	\fS_{D,\eP} := \qty{\epsilon\in\fS\mid D\epsilon = 0 \text{ and } \eP\epsilon = 0}.
\]
We previously discussed the integrability condition $R^D\epsilon=0$ for the existence of $D$-parallel spinors (see equation~\eqref{eq:D-curvature-def} and the subsequent discussion). In the case of constrained Killing spinors $\epsilon\in\fS_{D,\eP}$, this can be supplemented by further integrability conditions, the first-order such conditions being schematically (see, e.g., \cite{Gran2019})
$\comm{\eP}{\eP}\epsilon = 0$, \smash{$\widehat{D}\eP \epsilon := \comm{D}{\eP}\epsilon = 0$},
where~${\comm{-}{-}}$ in both expressions is a commutator. These conditions require some interpretation which depends on the structure on the bundle $E$ in which the additional fermions take their values. If $\eP\in\End\Sbundle$, so $E\to M$ is the trivial line bundle, then the meaning of the expressions is clear. The same holds if $E$ is simply the sum of such bundles. In the more general situation, we may interpret $\eP$ as acting on $\Sbundle\otimes E$ as $\eP\otimes\Id_E$ and extend $D$ to a connection on $\Sbundle\otimes E$ which acts trivially on the second factor. In practice though, there may be some more natural interpretation in which these operations are ``twisted'' by some structure on $E$.

\subsubsection{Constrained Killing superalgebra}

As we did for $\fS_D$ in Section~\ref{sec:admissible-ksa}, we now seek to find a subalgebra of the algebra of Killing vectors which preserves this space under the action of the Lie derivative; a natural choice is
\[
	\fV_{D,\eP} := \big\{X\in\fiso(M,g) \mid \eL_X \beta = 0 \text{ and } \eLhat_X \eP = 0\big\},
\]
where again we must interpret the expression $\eLhat_X \eP$. In particular, we assume that there is some structure on $E$ such which gives rise to some generalisation $\eLhat$ of the Lie derivative which acts on $\Hom\Sbundle\otimes E$ -- if $E$ happens to be a trivial bundle, the spinorial Lie derivative naturally generalises, while if it is an associated bundle to the principal bundle for some gauge symmetry, a covariant derivative constructed from a gauge connection will suffice. We note that in the latter case, $\fV_{D,\eP}$ may not close under the Lie bracket, so we may need to refine its definition.

If $\fV_{D,\eP}$ (or some more carefully chosen subspace) is indeed a subalgebra of $\fiso(M,g)$, its action via the spinorial Lie derivative preserves $\fS_{D,\eP}$.
We now seek to give $\fK_{D,\eP} = \fV_{D,\eP}\oplus\fS_{D,\eP}$ the structure of a superalgebra using the Lie derivative and Dirac current as in Theorem~\ref{thm:killing-algebra-exist}. But as in the proof of that theorem, all that needs to be checked is that the $\comm{\fS_{D,\eP}}{\fS_{D,\eP}}$ bracket closes on $\fV_{D,\eP}$, and that the $\qty[\fS_{D,\eP},\fS_{D,\eP},\fS_{D,\eP}]$ Jacobi identity is satisfied. But we quickly see that this is the case if equations \eqref{eq:killing-alg-exist-1}--\eqref{eq:killing-alg-exist-3} as well as
\begin{equation}\label{eq:alg-constraint-cons}
	\eLhat_{\kappa(\epsilon,\zeta)}\eP = 0
\end{equation}
are satisfied for all $\epsilon,\zeta\in\fS_{D,\eP}$. We note that in a supergravity theory, $\beta$ and $\eP$ are defined in terms of background fields which we expect to be preserved by the Lie derivatives along Killing vectors generated by spinors in $\fS_{D,\eP}$, hence the condition above is also satisfied.

Assuming that the expressions above can be made precise as suggested, we should then generalise Definition~\ref{def:killing-spinor} to the constrained case, considering not just admissible connections $D$ but admissible \emph{pairs} $(D,\eP)$, incorporating the condition \eqref{eq:alg-constraint-cons} and possibly demanding that the admissibility conditions of Definition~\ref{def:killing-spinor} should only be imposed on $\ker\eP$, rather than the whole spinor bundle.

\subsubsection{Killing transport and filtered deformation}

We now consider whether an analogue of Theorem~\ref{thm:killing-algebra-filtered} holds, that is, whether it is the case that the Killing superalgebra is a filtered deformation of the Poincar\'e superalgebra in the constrained case. Recall that an important part of this proof was the localisation of the Killing vectors and spinors in $\fK_D$ at a point via Killing transport. Since constrained Killing spinors are still $D$-parallel with respect to a connection $D$, the same arguments hold, the only modification being that the Killing transport data of a constrained Killing spinor at $p\in M$ must lie in $\ker\eP_p$. Note that this does not mean that any spinor at $p$ annihilated by $\eP$ can be Killing-transported to define a constrained Killing spinor; as in the unconstrained case, a $D$-parallel global section with this transport data may not exist, and in the constrained case, even if one does exist, it need not be $\eP$-parallel.

One should thus expect an analogue to Theorem~\ref{thm:killing-algebra-filtered} for $\fK_{D,\eP}$ with essentially the same proof, the main difference being that the filtered subdeformation of $\fs$ one thus obtains must have $S'\simeq \evaluate^S_p\fS_{D,\eP}\subseteq\ker\eP_p\subseteq\Sbundle_p$.

We will not consider constrained Killing spinors in the rest of this work, primarily because, in contrast to the differential Killing spinor equation, the linear constraints apparently do not have an interpretation in terms of Lie algebra cohomology. We thus leave the details of the formalism suggested above to be fully treated in future publications.

\section{Filtered subdeformations of the Poincar\'e superalgebra}
\label{sec:filtered-def-poincare}

Theorem~\ref{thm:killing-algebra-filtered} establishes that any Killing (super)algebra is a filtered deformation of a graded subalgebra of a flat model (super)algebra $\fs$ (see Definition~\ref{def:flat-model-alg}). Let us now consider such deformations from a homological point of view. Throughout, let $(V,\eta)$ be an inner product space and let $S$ be a (possibly $N$-extended) spinor module of $\Spin(V)$. For ease of exposition, we now specialise to the case of a symmetric Dirac current $\kappa\colon \Odot^2 S \to V$, whence we work with \emph{super}algebras; the skew-symmetric case is entirely analogous. Eventually, we will further specialise to Lorentzian signature, hence Poincar\'e superalgebras, and make some further assumptions on the Dirac current and the graded subalgebras we work with in order to take advantage of some technical simplifications offered by the homogeneity theorem (see Theorem~\ref{thm:homogeneity}). We comment on how dependent our results are on these assumptions in Remark~\ref{rem:relaxing-homogeneity-assumption}. We will make extensive use of the background material on Spencer cohomology and filtered deformations presented in Section~\ref{sec:spencer-thy}.

\subsection{The Spencer (2,2)-cohomology}
\label{sec:poincare-spencer-22-cohomology}

We saw in Section~\ref{sec:filtered-def} that filtered deformations of graded Lie superalgebras are governed by their Spencer cohomology. We therefore need to understand the Spencer cohomology of the flat model superalgebra and its graded subalgebras before discussing their deformations.

\subsubsection{The complex}

A graded subalgebra $\fa$ of the flat model superalgebra $\fs$ defined by the data $(V,S,\kappa)$ takes the form $\fa=\fa_{-2}\oplus\fa_{-1}\oplus\fa_0$, where $\fa_{-1}=S'$ is a vector subspace of $S$, $\fa_{-2}=V'$ is a subspace of $V$ containing $\Im \kappa|_{\Odot^2S'}$ and $\fa_0=\fh$ is a subalgebra of $\fso(V)$ preserving $S'$ and $V'$. By Proposition~\ref{prop:def-seq-spencer-cocycle}, to study filtered deformations of $\fa$ we must first understand the degree-2 Spencer complex $\bigl(\ssC^{2,\bullet}(\fa_-;\fa),\partial\bigr)$. Recalling the definition of the Spencer complex from Section~\ref{subsec:spencer-cohomology}, the space of $(2,p)$-cochains consists of degree-2 maps $\Wedge^p\fa\to\fa$, whence we deduce that the degree-2 complex can be written as
\begin{align*}
	0 &\longrightarrow \Hom(V',\fh)
\longrightarrow \Hom\bigl(\Wedge^2V',V'\bigr)\oplus\Hom(V'\otimes S',S')\oplus\Hom\bigl(\Odot^2S',\fh\bigr)\\
	&\longrightarrow \Hom\bigl(V'\otimes\Odot^2S',V'\bigr)\oplus\Hom\bigl(\Odot^3S',S'\bigr)\longrightarrow 0.
\end{align*}
We denote the projections to the components of $\ssC^{2,2}(\fa_-;\fa)$ as follows
\begin{gather*}
	\pi_1\colon\ \ssC^{2,2}(\fa_-;\fa) \longrightarrow \Hom\bigl(\Wedge^2V',V'\bigr),\qquad
	\pi_2\colon\ \ssC^{2,2}(\fa_-;\fa) \longrightarrow \Hom(V\otimes S',S'),\\
	\pi_3\colon\ \ssC^{2,2}(\fa_-;\fa) \longrightarrow \Hom\bigl(\Odot^2S',\fh\bigr).
\end{gather*}

We now record a pair of results which will be useful for studying the cohomology of a class of graded subalgebras which we will be particularly interested in later on.

\begin{Lemma}\label{lemma:V-h-maps}
	Let $\fh$ be a subalgebra of $\fso(V)$ and $V'$ an $\fh$-submodule of $V$ such that
	\begin{enumerate}\itemsep=0pt
		\item[$(1)$] the restriction of $\eta$ to $V'$ is non-degenerate,
		\item[$(2)$] the action of $\fh$ on $V'$ is faithful,
	\end{enumerate}
	$($which hold in particular for $V'=V)$ and suppose that $\lambda\colon V'\to\fh$ is a linear map such that
	\[
		\lambda(v)w - \lambda(w)v = 0, \qquad \forall v,w\in V'.
	\]
	Then $\lambda=0$.
\end{Lemma}

\begin{proof}
	Using that $\lambda$ takes values in $\fh\subseteq \fso(V)$, for all $u,v,w\in V'$ we have
	\begin{align*}
		\eta(u,\lambda(v)w)
		&= \eta(u,\lambda(w)v)
		= -\eta(\lambda(w)u,v)
		= -\eta(\lambda(u)w,v) \\
		&
		= \eta(w,\lambda(u)v)
		= \eta(w,\lambda(v)u)
		= -\eta(\lambda(v)w,u).
	\end{align*}
	Thus by non-degeneracy we have $\lambda(v)w=0$ for all $v,w\in V'$, so since $\fh$ acts faithfully on $V'$, $\lambda=0$.
\end{proof}

\begin{Corollary}\label{coro:21-22-injective}
	Let $\fa=V\oplus S'\oplus \fh$ be a graded subalgebra of $\fs$ with $\fa_{-2}=V$. Then the composition
$\pi_1\circ\partial\colon \Hom(V,\fh) \to \Hom\bigl(\Wedge^2 V,V\bigr)
$
	is injective, whence ${\ssH^{2,1}(\fa_-;\fa)=\ssZ^{2,1}(\fa_-;\fa)=0}$. If~${\fh=\fso(V)}$, $\pi_1\circ\partial$ is an isomorphism.
\end{Corollary}
\begin{proof}
Let $\lambda\in\ker(\pi_1\circ\partial)$. Then for all $v,w\in V$,
$(\pi_1\circ\partial)(\lambda)(v,w) = \partial\lambda(v,w) = \lambda(v)w - \lambda(w)v = 0$.
The first claim then follows from Lemma~\ref{lemma:V-h-maps}. The second claim follows by a~dimension count.
\end{proof}

\begin{Remark}
	The second claim in the result above can be restated as follows: for any linear map~${\alpha\colon\Wedge^2 V\to V}$, there exists a unique linear map $\lambda\colon V\to\fso(V)$ satisfying $\alpha(v,w)=\partial\lambda(v,w) = \lambda(v)w-\lambda(w)v$. We will use this fact repeatedly.
\end{Remark}

\subsubsection{The cohomology}

The discussion above gives us a characterisation of the $(2,2)$-cohomology group of the whole flat model superalgebra $\fs$ which we will later be able to bootstrap in order to study a particular class of filtered subdeformations.

\begin{Lemma} \label{lemma:H22-full-flat-model}
	Let $\alpha+\beta+\gamma$ be a cocycle in $\ssC^{2,2}(\fs_-;\fs)$ where
	\begin{align}\label{eq:22-cocycle-comps}
		\alpha\in\Hom\bigl(\Wedge^2V,V\bigr),
		\qquad\beta\in\Hom(V\otimes S,S),
		\qquad\gamma\in\Hom\bigl(\Odot^2S,\fso(V)\bigr).
	\end{align}
	Then the cohomology class $[\alpha+\beta+\gamma]\in \ssH^{2,2}(\fs_-;\fs)$ has a unique representative $\beta'+\gamma'$ with no $\Hom\bigl(\Wedge^2 V,V\bigr)$ component. We denote the space of such \emph{normalised cocycles} by $\cH^{2,2}$
	\[
		\cH^{2,2}
		= \qty{\beta+\gamma \in \ssZ^{2,2}(\fs_-;\fs)
			 \mid
				\beta\in\Hom(V\otimes S,S), \gamma\in\Hom\big(\Odot^2S,\fso(V)\big)
			},
	\]
	and then we have
$\ssH^{2,2}(\fs_-;\fs) \cong \cH^{2,2}
$
	as $\fso(V)$-modules.
\end{Lemma}

\begin{proof}
	Using the isomorphism $\pi_1\circ\partial$, there exists a unique map $\lambda\in\Hom(V,\fso(V))$ such that~${\partial\lambda(v,w)=\lambda(v)w - \lambda(w)v = \alpha(v,w)}$.
	Thus, if we define $\beta'$ and $\gamma'$ by
	\begin{gather*}
		\beta'(v,s) = \beta(v,s) - \partial\lambda(v,s)
			= \beta(v,s) - \lambda(v)\cdot s,\\
		\gamma'(s,s') = \gamma(s,s') - \partial\lambda(s,s')
			= \gamma(s,s) + \lambda(\kappa(s,s')),
	\end{gather*}
	we have
$[\alpha+\beta+\gamma]
		= [\alpha+\beta+\gamma-\partial\lambda]
		= [\beta' + \gamma']$,
	which proves the first claim. It immediately follows that, as $\fso(V)$-modules, we have $\ssZ^{2,2}(\fs_-;\fs)=\cH^{2,2}\oplus\ssB^{2,2}(\fs_-;\fs)$ and $\ssH^{2,2}(\fs_-;\fs) \cong \cH^{2,2}$.
\end{proof}

It is worth pausing to explicitly examine the ``normalised'' Spencer cocycle conditions for $\beta+\gamma\in\cH^{2,2}$. Written in polarised form (explained in Section~\ref{sec:polarisation}), these are
\begin{gather}
	2\kappa(s,\beta(v,s)) + \gamma(s,s)v = 0, \label{eq:cocycle-1}	\\
	\beta(\kappa_s,s) + \gamma(s,s)\cdot s = 0	\label{eq:cocycle-2}
\end{gather}
for $v\in V$ and $s\in S$, where we introduce the notation $\kappa_s:=\kappa(s,s)$. We will refer to these as the \emph{first} and \emph{second} \emph{$($normalised$)$ Spencer cocycle conditions} for $\fs$ respectively. We note that the first equation uniquely defines $\gamma\colon\Odot^2S\to\fso(V)$ in terms of $\beta\colon V\otimes S\to S$ and constrains the endomorphism $v\mapsto \kappa(s,\beta(v,s))$ of $V$ to lie in $\fso(V)$ for all $s\in S$. In explicit cohomology calculations (see \cite{Beckett2021,deMedeiros2016,deMedeiros2018,Figueroa-OFarrill2017}), we first determine the space of maps $\beta$ satisfying this constraint and then substitute the resulting expressions for $\beta$ and $\gamma$ into the second equation and solve it.

We also note here that $\beta+\gamma\in\cH^{2,2}$ (or the corresponding cohomology class) is invariant under the action of a subalgebra $\fh\in\fso(V)$ if and only if $\beta$ is $\fh$-invariant (that is, $\beta$ is $\fh$-equivariant when considered as a map $V\otimes S\to S$); indeed, if $A\in\fh$, $v\in V$ and $s\in S$ then we have~${(A\cdot(\beta+\gamma))(v,s) = (A\cdot\beta)(v,s) = A\cdot (\beta(v,s)) - \beta(Av,s) - \beta(A,v\cdot s)}
$
so clearly~$\beta$ is invariant if $\beta+\gamma$ is, while on the other hand if $\beta$ is invariant then using the first cocycle condition (once in the form above and once in depolarised form) as well as $\fso(V)$-invariance of~$\kappa$, we have\looseness=-1
\begin{align*}
(A\cdot\gamma)(s,s)v &= \comm{A}{\gamma(s,s)}v - 2\gamma(A\cdot s,s)v \\
	&= 2A(\kappa(s,\beta(v,s))) - 2\kappa(s,\beta(Av,s)) - 2\kappa(A\cdot s,\beta(v,s)) - 2\kappa(s,\beta(v,A\cdot s)) \\
	&= 2\kappa(s,A\cdot\beta(v,s)) - 2\kappa(s,\beta(Av,s)) - 2\kappa(s,\beta(v,A\cdot s))= 0,
\end{align*}
thus $\gamma$ is $\fh$-invariant, whence $\beta+\gamma$ is too.

\begin{Remark}
	Let us pause briefly to relate these cocycle conditions back to our theory of Killing superalgebras and to explain the coincidences in notation. First note that there is a canonical homomorphism $\Hom(V\otimes S,S)\cong V^*\otimes\End S$, so we can view the cocycle component $\beta$ as a~spinor endomorphism-valued 1-form on $V$. Under this interpretation, we write $\beta(v)s:=\beta(v,s)$.
	
	The definition of admissible connections, Definition~\ref{def:killing-spinor}, involves three conditions on the 1-form $\beta=D-\nabla\in\Omega^1(M;\End\Sbundle)$ and the section $\gamma\in\Gamma\bigl(\Odot^2\Sbundle^*\otimes TM\bigr)$ which is defined in terms of $\beta$ by equation~\eqref{eq:gamma-geom-def}. Conditions~(1) and~(2) are algebraic in $\beta$, $\gamma$ and condition~(3) is differential. For a symmetric Dirac current, we can polarise the definition of $\gamma$ and the algebraic conditions. Then, taking them pointwise, the definition of $\gamma$ and condition~(1) are equivalent to cocycle condition \eqref{eq:cocycle-1} above, while condition~(2) is precisely \eqref{eq:cocycle-2} (where the fibre of $\Sbundle$ is~$S$ and the same Dirac current is chosen). Thus, if we wish to find admissible connections on a~spinor bundle $\Sbundle$ equipped with a Dirac current $\kappa$ over a spin manifold $(M,g)$, we need only compute~$\cH^{2,2}$ to determine a local ansatz for the allowed 1-forms $\beta$ and then check whether or under what circumstances the remaining condition (namely, $\eL_{\kappa_\epsilon}\beta=0$ for all $D$-parallel spinors~$\epsilon$) holds, which we typically do using integrability conditions.
\end{Remark}

\subsection{General filtered deformations and cohomology}
\label{sec:filtered-def-cohom}

Recall from Section~\ref{sec:filtered-def} that we can describe filtered deformations of a $\ZZ$-graded superalgebra as deformed brackets on the same underlying vector space as the original algebra. In particular, a filtered deformation $\fatilde$ of a graded subalgebra $\fa=V'\oplus S'\oplus \fh$ of $\fs$ has the brackets
\begin{gather}
	\comm{A}{B} = AB-BA,\qquad
	\comm{s}{s} = \kappa_s + \gamma(s,s),\qquad
	\comm{A}{v} = Av + \delta(A,v),
	\nonumber\\
	\comm{v}{s} = \beta(v,s),\qquad\comm{A}{s} = A\cdot s
	, \qquad\comm{v}{w} = \alpha(v,w) + \theta(v,w),\label{eq:def-brackets-general}
\end{gather}
where $A,B\in\fh$, $v,w\in V'$, $s\in S'$, and
$\alpha\colon \Wedge^2V' \to V'$, $
	\beta\colon V'\otimes S' \to S'$, $
	\gamma\colon \Odot^2 S' \to \fh$, $
	\delta\colon \fh\otimes V' \to \fh
$
are the components of the degree-2 deformation map, while
$\theta\colon \Wedge^2 V' \to \fh
$
is the degree-4 deformation map. We denote the full degree-2 deformation by
${\mu=\alpha+\beta+\gamma+\delta\colon \fa\otimes\fa \to \fa}
$
so the defining sequence (as in equation~\eqref{eq:defining-sequence}) for this deformation is $(\mu,\theta,0,\dots)$. From Proposition~\ref{prop:def-seq-spencer-cocycle}, we have the following homological information:
\begin{gather}
	\mu \in \ssZ^{2,2}(\fa;\fa),	\label{eq:mu-CE-cocycle}	\\
	\mu|_{\fa_-\otimes\fa_-} \in \ssZ^{2,2}(\fa_-;\fa),	\label{eq:mu-spencer-cocycle}	\\
	\qty[\mu|_{\fa_-\otimes\fa_-}] \in \ssH^{2,2}(\fa_-;\fa)^\fh, \label{eq:mu-invt-class}
\end{gather}
where in the first expression we have graded the full Chevalley--Eilenberg complex in the same way that we graded the Spencer complex in equation~\eqref{eq:spencer-complex-grading}.

\subsubsection{Unpacking the cocycle conditions}

Let us consider these conditions in more detail. Equation~\eqref{eq:mu-CE-cocycle} is equivalent to the following system of linear equations:
\begin{gather}
	\alpha(\kappa_s,v) + 2\kappa(s,\beta(v,s)) + \gamma(s,s)v = 0,\label{eq:22-spencer-cocycle-vss} \\
	\beta(\kappa_s,s) + \gamma(s,s)\cdot s = 0,\label{eq:22-spencer-cocycle-sss}	\\
	A\alpha(v,w) - \alpha(Av,w) - \alpha(v,Aw) + \delta(A,w)v - \delta(A,v)w = 0,\label{eq:alpha-delta} \\
	A\cdot(\beta(v,s)) - \beta(Av,s) - \beta(v,A\cdot s) - \delta(A,v)\cdot s = 0,\label{eq:beta-delta} \\
	\comm{A}{\gamma(s,s)} - 2\gamma(A\cdot s,s) + \delta(A,\kappa_s) = 0,\label{eq:gamma-delta} \\
	\delta(\comm{A}{B},v) - \comm{A}{\delta(B,v)} + \comm{B}{\delta(A,v)} - \delta(A,Bv) + \delta(B,Av) = 0 \label{eq:delta-coycle-cond}
\end{gather}
for all $A,B\in\fh$, $v,w\in V'$ and $s\in S'$, where again $\kappa_s=\kappa(s,s)$. Equation~\eqref{eq:mu-spencer-cocycle} is equivalent to the first pair of equations in this system -- we will refer to these as the \emph{Spencer cocycle conditions} for $\fa$. Equation~\eqref{eq:mu-invt-class} is equivalent to the following: for each $A\in\fh$, there exists a~cochain $\chi_A\in \ssC^{2,1}(\fa_-;\fa)=\Hom(V',\fh)$ such that
\begin{equation}\label{eq:chiA-def}
	A\cdot(\alpha+\beta+\gamma) = \partial\chi_A.
\end{equation}
This is implied by equation~\eqref{eq:mu-CE-cocycle}: simply taking $\chi_A=\imath_A\delta$, the above equation is equivalent to equations~\eqref{eq:alpha-delta}--\eqref{eq:gamma-delta}. In the case that $V'=V$, we have a converse: it follows from Corollary~\ref{coro:21-22-injective} that \eqref{eq:chiA-def} defines $\chi_A$ uniquely and that the assignment $A\mapsto\chi_A$ is linear. Thus we can define a unique linear map $\delta\colon \fh\otimes V\to\fh$ by $\delta(A,v):=\chi_A(v)$. Then
\begin{gather*}
	\delta(\comm{A}{B},v) - \comm{A}{\delta(B,v)} + \comm{B}{\delta(A,v)} - \delta(A,Bv) + \delta(B,Av)\\
	\qquad= \chi_{\comm{A}{B}}(v) - \comm{A}{\chi_B(v)} + \comm{B}{\chi_A(v)}	- \chi_A(Bv) + \chi_B(Av)	\\
	\qquad= \qty(\chi_{\comm{A}{B}} - A\cdot\chi_B + B\cdot\chi_A)(v).
\end{gather*}
On the other hand, by \eqref{eq:chiA-def} we have
$\partial(\chi_{\comm{A}{B}} - A\cdot\chi_B + B\cdot\chi_A) = 0$,
so $\chi_{\comm{A}{B}} - A\cdot\chi_B + B\cdot\chi_A\in\ssZ^{2,1}(\fa_-;\fa)$. But this space is zero by Corollary~\ref{coro:21-22-injective}. Hence $\delta$ satisfies \eqref{eq:delta-coycle-cond}. Thus we have shown the following.

\begin{Lemma}\label{lemma:invt-spencer-cocycle}
	Let $\mu\in\ssZ^{2,2}(\fa;\fa)$. Then $\mu|_{\fa_-\otimes\fa_-}\in\ssZ^{2,2}(\fa_-;\fa)$ with $\fh$-invariant cohomology class. If $V'=V$, we conversely have that if $\mu_-\in\ssZ^{2,2}(\fa_-;\fa)$ is a Spencer cocycle such that with $\fh$-invariant cohomology class then there exists a unique cocycle $\delta\in\ssZ^1(\fh;V^*\otimes\fh)\subseteq\Hom(\fh\otimes V,\fh)$ such that $\mu=\mu_-+\delta\in\ssZ^{2,2}(\fa;\fa)$.
\end{Lemma}

We can be more explicit about the form of $\delta$ in the $V=V'$ case: by Corollary~\ref{coro:21-22-injective}, there exists a unique map $\lambda\colon V\to\fso(V)$ such that
$\alpha(v,w) = \lambda(v)w - \lambda(w)v$.
One can then easily check that
$\delta(A,v) = \comm{A}{\lambda(v)} - \lambda(Av)
$
solves equations \eqref{eq:alpha-delta}--\eqref{eq:gamma-delta}.

\subsubsection{Jacobi identities}

We now consider the conditions imposed on the deforming maps by the Jacobi identities. These split into three types: those which are independent of the deforming maps, which are identically satisfied; those which involve only the deforming maps of degree~2, which are linear in those maps and which are equivalent to the homological conditions~\eqref{eq:22-spencer-cocycle-vss}--\eqref{eq:delta-coycle-cond}; and those which involve the degree-4 map $\theta$, all of which are quadratic in the deforming maps. We will not explicitly write the first two types, only the third type, since the equations are new. We denote by $[\fatilde_i,\fatilde_j,\fatilde_k]$ the Jacobi identity on the $(i, j,k)$-th homogeneous subspace of $\fa$. The Jacobi identities break down as follows:
\begin{itemize}\itemsep=0pt
	\item $[\fatilde_0,\fatilde_0,\fatilde_0]$: Identically satisfied by the Jacobi identity on $\fh\subseteq \fso(V)$.
	\item $[\fatilde_0,\fatilde_0,\fatilde_{-1}]$: Identically satisfied since $S'$ is a representation of $\fh$.
	\item $[\fatilde_0,\fatilde_0,\fatilde_{-2}]$: Equivalent to the cocycle condition \eqref{eq:delta-coycle-cond}, using that $V'$ is a representation of $\fh$.
	\item $[\fatilde_0,\fatilde_{-1},\fatilde_{-1}]$: Equivalent to cocycle condition \eqref{eq:gamma-delta}, using the $\fso(V)$-equivariance of $\kappa$.
	\item $[\fatilde_0,\fatilde_{-1},\fatilde_{-2}]$: Equivalent to cocycle condition~\eqref{eq:beta-delta}.
	\item $[\fatilde_0,\fatilde_{-2},\fatilde_{-2}]$: This identity has a $V'$-component, which is equivalent to the cocycle condition \eqref{eq:alpha-delta}, and an $\fh$-component which is quadratic in the deformation maps
		\begin{equation}\label{eq:jacobi-022}
			(A\cdot\theta)(v,w) = \delta(\delta(A,v),w) - \delta(\delta(A,w),v) - \delta(A,\alpha(v,w))
		\end{equation}
	for $A\in\fh$ and $v,w\in V'$.
	\item $[\fatilde_{-1},\fatilde_{-1},\fatilde_{-1}]$: Equivalent to the second Spencer cocycle condition \eqref{eq:22-spencer-cocycle-sss}.
	\item $[\fatilde_{-1},\fatilde_{-1},\fatilde_{-2}]$: Like the $[\fatilde_0,\fatilde_{-2},\fatilde_{-2}]$ identity, there is a $V'$-component and a $\fh$-component to this Jacobi identity. The former is the first Spencer cocycle condition \eqref{eq:22-spencer-cocycle-vss}, while the latter is the quadratic condition
		\begin{equation}\label{eq:jacobi-112}
			\theta(\kappa_s,v) + \delta(\gamma(s,s),v) + 2\gamma(s,\beta(v,s)) = 0
		\end{equation}
	for $s\in S'$ and $v\in V'$.
	\item $[\fatilde_{-1},\fatilde_{-2},\fatilde_{-2}]$: This identity gives the quadratic condition
		\begin{equation}\label{eq:jacobi-122}
			\theta(v,w)\cdot s + \beta(\alpha(v,w),s) - \beta(v,\beta(w,s)) + \beta(w,\beta(v,s)) =0,
		\end{equation}
	for $s\in S'$, $v,w\in V'$.
	\item $[\fatilde_{-2},\fatilde_{-2},\fatilde_{-2}]$: This has components in $V'$ and $\fh$, both of which are quadratic
		\begin{gather}
			\alpha(\alpha(u,v),w)
				+ \theta(u,v)w
				+ \text{cyclic perms.} = 0,\label{eq:jacobi-222a}\\
			\theta(\alpha(u,v),w)
				+ \delta(\theta(u,v),w)
				+ \text{cyclic perms.} =0\label{eq:jacobi-222b}
		\end{gather}
	for $u,v,w\in V'$.
\end{itemize}

In particular, recalling Lemma~\ref{lemma:invt-spencer-cocycle}, we have the following.

\begin{Proposition}\label{prop:gen-filtered-defs-cocycle}
	Let $\fa=V'\oplus S'\oplus\fh\subseteq\fs$ be a graded subalgebra, let $\alpha+\beta+\gamma+\delta\in\ssZ^{2,2}(\fa;\fa)$ and suppose we have a map $\theta\colon\Wedge^2V'\to\fh$. Then the brackets \eqref{eq:def-brackets-general} define a filtered deformation~$\fatilde$ of $\fa$ if and only if $\theta$ satisfies the equations~\eqref{eq:jacobi-022}--\eqref{eq:jacobi-222b}.
\end{Proposition}

Thus the task of determining the possible filtered deformations of a graded subalgebra $\fa\subseteq\fs$ consists of solving the linear equations \eqref{eq:22-spencer-cocycle-vss}--\eqref{eq:delta-coycle-cond} for $\mu=\alpha+\beta+\gamma+\delta$ and then the quadratic equations \eqref{eq:jacobi-022}--\eqref{eq:jacobi-222b} for $\theta$. We can interpret a $(2,2)$-cocycle $\mu$ as an \emph{infinitesimal filtered deformation} of $\fa$, and this deformation \emph{integrates} to an actual filtered deformation $\fatilde$ of $\fa$ if and only if there exists a solution $\theta$ to the equations \eqref{eq:jacobi-022}--\eqref{eq:jacobi-222b}. These equations are highly non-trivial in general, and even if an infinitesimal deformation integrates, the full deformation need not be unique. It is also difficult to give a homological characterisation of when two pairs~$(\mu,\theta)$ and $(\mu',\theta')$ give rise to isomorphic deformations.

\subsection{Homogeneity and high supersymmetry}
\label{sec:homogeneity}

We can get much better homological control over the problem of determining filtered deformations by passing to a particular class of graded subalgebras $\fa$. We have already seen in Lemmas~\ref{lemma:V-h-maps} and \ref{lemma:invt-spencer-cocycle} that taking $V=V'$ gives us significant simplifications. We would now like to guarantee that $S'$ is sufficiently large so that $V=\comm{S'}{S'}$, which will allow us to fully characterise filtered deformations of $\fa$ in terms of Spencer cohomology. To do this, we will make additional assumptions which allow us to use the homogeneity theorem, but we note that the essential ingredient for the results of this section is the fact that $V=\comm{S'}{S'}$.

\subsubsection{The homogeneity theorem}

We say that a symmetric Dirac current $\kappa$ is \emph{causal} if $\kappa_s=\kappa(s,s)$ is either timelike or null for all~${s\in S}$. In Lorentzian signature, a symmetric, causal Dirac current on the real pinor module $\ssS$ always exists in any dimension and can be used to construct Dirac currents on extended spinor modules. The existence of symmetric currents can be deduced from the Alekseevsky--Cort\'es classification \cite{Alekseevsky1997}; that they can be chosen to be causal must be shown case-by-case \cite{Hustler2016}. These can be used to build symmetric Dirac currents on ($N$-extended) spinor modules $S$ by restriction and tensoring with bilinears on an auxiliary module. However such a map is constructed, we have the following.

\begin{Theorem}[homogeneity theorem \cite{Figueroa-OFarrill2012,Hustler2016}]\label{thm:homogeneity}
	Let $(V,\eta)$ be a Lorentzian vector space of dimension $\dim V>2$ and $S$ a $($possibly $N$-extended$)$ spinor module of $\fso(V)$ with a symmetric, causal Dirac current $\kappa\colon \Odot^2 S \to V$. If $S'$ is a vector subspace of $S$ with $\dim S'>\frac{1}{2}\dim S$, then~${\kappa|_{\Odot^2S'}}$ is surjective onto $V$.
\end{Theorem}

For all $A\in\fso(V)$ and $s\in S$, $A\cdot s=0$ implies that $A \kappa_s=0$ by $\fso(V)$-equivariance of the Dirac current. Since $\fso(V)$ acts faithfully on $V$, the homogeneity theorem has the following corollary.

\begin{Corollary}\label{coro:homog-faithful}
	If the Dirac current $\kappa$ is causal and $\dim S'>\frac{1}{2}\dim S$, then the annihilator of~$S'$ in $\fso(V)$ is trivial. In particular, any subalgebra $\fh$ of $\fso(V)$ which preserves $S'$ acts faithfully on~$S'$.
\end{Corollary}

Guided by the homogeneity theorem (and by physics), let us now fix $\fs$ to be the Poincar\'e superalgebra associated to a Lorentzian inner product space $(V,\eta)$ and spinor module $S$ with symmetric causal Dirac current $\kappa$. Note that for such algebras we have $\comm{S'}{S'}=V$ for $\dim S'>\frac{1}{2}\dim S$, which prompts the following definition.

\begin{Definition}
	A graded subalgebra $\fa=V\oplus S'\oplus\fh$ of $\fs$ with $\dim S'> \frac{1}{2}\dim S$ is said to be \emph{highly supersymmetric}, and it is \emph{maximally supersymmetric} if $S'=S$.
	
	A filtered subdeformation $\fatilde$ of $\fs$ is said to be \emph{highly supersymmetric} if $\dim \fatilde_{\overline 1}> \frac{1}{2}\dim S$, and it is \emph{maximally supersymmetric} if $\dim \fatilde_{\overline 1} = \dim S$.
\end{Definition}

Clearly a filtered subdeformation $\fatilde$ is highly (maximally) supersymmetric if and only if its associated graded superalgebra $\fa$ is highly (maximally) supersymmetric as a graded subalgebra of $\fs$.

\begin{Remark}\label{rem:relaxing-homogeneity-assumption}
	The significance of the homogeneity theorem and the assumption of high supersymmetry is that they ensure that $\comm{S'}{S'}=V$, or $\comm{\fa_{-1}}{\fa_{-1}}=\fs_{-2}$. If one wishes to relax the assumptions which were made on the signature and Dirac current in order to invoke the theorem, one can simply impose this restriction on the graded subalgebra $\fa$ instead of the assumption of high supersymmetry and recover the results of the remainder of Sections~\ref{sec:filtered-def-poincare} and~\ref{sec:high-susy-spin-mfld}. One might term this \emph{sufficiently supersymmetry}.
\end{Remark}

\subsubsection{Homological properties}

Recall Definition~\ref{def:fund-trans-prolong} of some homological properties of graded Lie superalgebras. We now show that highly supersymmetric graded subalgebras of $\fs$ satisfy some of these.

\begin{Lemma}[\cite{Figueroa-OFarrill2017_1}]\label{lemma:high-susy-subalg-coho}
	Let $\fa$ be a highly supersymmetric graded subalgebra of $\fs$. Then
	\begin{itemize}\itemsep=0pt
	\item $\fa$ is fundamental and transitive,
	\item $\fa$ is a full prolongation of degree $2$,
	\item $\ssH^{d,2}(\fa_-;\fa)=0$ for all even $d>2$.
	\end{itemize}
\end{Lemma}
\begin{proof}
	Theorem~\ref{thm:homogeneity} precisely says that $\fa$ is fundamental, and transitivity follows from faithfulness of the action of $\fh=\fa_0\subseteq\fso(V)$ on $V=\fa_{-2}$ or $S'=\fa_{-1}$. Corollary~\ref{coro:21-22-injective} then says that~${\ssH^{2,1}(\fa_-;\fa)=0}$. For $d>2$, $\ssH^{d,1}(\fa_-;\fa)=0$ trivially since $\ssC^{d,1}(\fa_-;\fa)=0$. The beginning of the degree-4 Spencer complex is
	\begin{align*}
		0	&\longrightarrow \Hom\bigl(\Wedge^2V,\fh\bigr)	
			\longrightarrow \Hom\bigl(\Wedge^3 V,V\bigr)\oplus \Hom\bigl(\Wedge^2V\otimes S',S'\bigr) \oplus \Hom\bigl(V\otimes \Odot^2 S', \fh\bigr) \\
&\longrightarrow \cdots,
	\end{align*}
	and for $\theta\in \ssC^{4,2}(\fa_-;\fa)=\Hom\bigl(\Wedge^2V,\fh\bigr)$, we have
$\partial\theta (v,s,s) = \theta(v,\kappa_s)$
	for all $v\in V$, $s\in S'$, so again by homogeneity, $\partial\theta=0$ if and only if $\theta=0$. This shows that $\ssH^{d,4}(\fa_-;\fa)=0$, and we have~${\ssH^{d,2}(\fa_-;\fa)=0}$ for all $d>4$ since $\ssC^{d,2}(\fa_-;\fa)=0$.
\end{proof}

The following pivotal observation is a straightforward generalisation to all spacetime dimensions of \cite[Proposition 6]{Figueroa-OFarrill2017_1} and also of the analogous results for \emph{maximally} supersymmetric graded subalgebras (those for which $\fa_{-1}=S$) in \cite{deMedeiros2016,Figueroa-OFarrill2017}. For the proof, we recall Propositions~\ref{prop:def-seq-spencer-cocycle} and \ref{prop:def-seq-redef} which characterise the defining sequences of filtered deformations in terms of Spencer cohomology and tell us how much freedom we have to redefine them.

\begin{Proposition}[\cite{Figueroa-OFarrill2017_1}]\label{prop:highly-susy-filtered-def}
	Let $\fa=V'\oplus S'\oplus\fh$ be a graded subalgebra of $\fs$ and $\fatilde$ a filtered deformation of $\fa$ with defining sequence $(\mu,\theta,0,\dots)$ as above. Then
	\begin{enumerate}\itemsep=0pt
		\item[$(1)$] $\mu|_{\fa_-\otimes\fa_-}$ is a cocycle in $\ssC^{2,2}(\fa_-;\fa)$, and
$[\mu|_{\fa_-\otimes\fa_-}] \in \ssH^{2,2}(\fa_-;\fa)^\fh$.
		Furthermore, $\mu$ is a~cocycle in $\ssC^2(\fa;\fa)$.	
		\item[$(2)$] If $\fatilde$ is highly supersymmetric and $\fatilde'$ is another filtered deformation of $\fa$ with degree-$2$ deformation map $\mu'$ such that $[\mu'|_{\fa_-\otimes\fa_-}]=[\mu|_{\fa_-\otimes\fa_-}]$ then $\fatilde\cong\fatilde'$ as filtered Lie superalgebras.	
	\end{enumerate}
\end{Proposition}

\begin{proof}
	The first statement follows by Proposition~\ref{prop:def-seq-spencer-cocycle}. For the second, if $[\mu'|_{\fa_-\otimes\fa_-}]=[\mu|_{\fa_-\otimes\fa_-}]$ then $(\mu'-\mu)|_{\fa_-\otimes\fa_-}$ is a Spencer coboundary, so by point \ref{item:def-seq-redef-1} of Proposition~\ref{prop:def-seq-redef} there is a new defining sequence $\{\mu'',\theta'',0,\dots\}$ of $\fatilde'$ such that $\mu''|_{\fa_-\otimes\fa_-}=\mu|_{\fa_-\otimes\fa_-}$. Then, since $\fa$ is a full prolongation of degree 2 by Lemma~\ref{lemma:high-susy-subalg-coho}, by point \ref{item:def-seq-redef-3} of Proposition~\ref{prop:def-seq-redef} we can again find a new defining sequence $\{\mu''',\theta''',0,\dots\}$ of $\fatilde'$ such that $\mu'''=\mu$ and $\theta'''=\theta'''|_{\fa_-\otimes\fa_-}=\theta|_{\fa_-\otimes\fa_-}=\theta$.
\end{proof}

Note that we already saw the first statement of the proposition above in Section~\ref{sec:filtered-def-cohom}. The second statement allows us to strengthen Proposition~\ref{prop:gen-filtered-defs-cocycle}; by Lemma~\ref{lemma:invt-spencer-cocycle}, a choice of $\fh$-invari\-ant~$\mu_{\fa_-\otimes\fa_-}$ determines $\mu$ for $V=V'$; we now see that the cohomology class $[\mu_{\fa_-\otimes\fa_-}]$ determines the entire deformation up to isomorphism in the highly supersymmetric case.

What this means in practice is that in order to find the infinitesimal deformations of a~highly supersymmetric graded subalgebra $\fa$, we compute $\ssH^{2,2}(\fa_-;\fa)$ and its $\fh$-invariants. The other deforming maps $\delta$ and $\theta$ will then be determined by the Jacobi identities. Note, though, that the Jacobi identities which are quadratic in the deformation maps may still obstruct the ``integration'' of the infinitesimal deformation $[\mu_{\fa_-\otimes\fa_-}]$ to a filtered deformation $\fatilde$ of $\fa$.

\subsubsection{Reparametrising the deformation}
\label{sec:reparam-defs}

It will be useful for us to relate $\ssH^{2,2}(\fa_-;\fa)$ to some other cohomology groups. The inclusion $i\colon\fa\hookrightarrow\fs$ of a graded subalgebra into $\fs$ induces the following maps of cochains
\begin{align*}
	i_*\colon\ \ssC^{\bullet,\bullet}(\fa_-;\fa) \to \ssC^{\bullet,\bullet}(\fa_-;\fs),
		\qquad \phi\mapsto i_*\phi = i\circ\phi,\\
	i^*\colon\ \ssC^{\bullet,\bullet}(\fs_-;\fs) \to \ssC^{\bullet,\bullet}(\fa_-;\fs),
		\qquad \phi\mapsto i^*\phi = \phi\circ i,
\end{align*}
where $\fs$ is viewed as representation of $\fa_-$ via the restriction of the adjoint representation, and the cochain complex $\ssC^{\bullet,\bullet}(\fa_-;\fs)$ has been given a grading in the same manner as the Spencer complex (equation~\eqref{eq:spencer-complex-grading}); the differential respects this grading and the maps $i_*$, $i^*$ respect both this grading and the homological grading -- that is, they commute with the differential $\partial$. Thus these maps in turn induce maps in cohomology also denoted $i_*$ and $i^*$.

Let us consider these maps in degree $(2,2)$ for the highly supersymmetric case. Recall that, by Lemma~\ref{lemma:H22-full-flat-model}, $\ssH^{2,2}(\fs_-;\fs)\cong\cH^{2,2}$ where the latter is the space of normalised Spencer $(2,2)$-cocycles of $\fs$.

\begin{Lemma}\label{lemma:maps-in-cohomology}
	Let $\fa$ be a highly supersymmetric graded subalgebra of $\fs$. Then we have a diagram with exact rows and columns
	\[
	\begin{tikzcd}
		& & & 0 \arrow[d]\\
		& & & \ssH^{2,2}(\fa_-;\fa) \arrow[d,"i_*"] \\
		0 \arrow[r]
		& \cK^{2,2}(\fa_-) \arrow[r]
		& \ssH^{2,2}(\fs_-;\fs) \cong \cH^{2,2} \arrow[r, "i^*"]
		& \ssH^{2,2}(\fa_-;\fs),
	\end{tikzcd}
	\]
	where
$\cK^{2,2}(\fa_-) := \big\{ \beta+\gamma\in \cH^{2,2} \mid \beta|_{V\otimes S'}=0 \big\} \cong \ker i^*$.
\end{Lemma}

\begin{proof}
	For $\beta\colon V\otimes S\to S$, $\gamma\colon \Odot^2S\to\fso(V)$ such that $\beta+\gamma\in\cH^{2,2}$, we have $i^*[\beta+\gamma]=0$ if and only if $\beta|_{V\otimes S'}+\gamma|_{\Odot^2S'}=\partial\lambda$, where $\lambda\colon V\to\fso(V)$. This implies that $\partial\lambda|_{\Wedge^2 V}=0$, so $\lambda=0$ by Lemma~\ref{lemma:V-h-maps}, thus $\beta|_{V\otimes S'}=0$ and $\gamma|_{\Odot^2S'}=0$. Note that $\beta|_{V\otimes S'}=0$ implies $\gamma|_{\Odot^2S'}=0$ by the cocycle condition for $\beta+\gamma$, so the isomorphism $\cH^{2,2}\cong\ssH^{2,2}(\fs_-;\fs)$ restricts to give $\cK^{2,2}(\fa_-)\cong \ker i^*$.
	
	For $\alpha\colon \Wedge^2V\to V$, $\beta\colon V\otimes S'\to S'$, $\gamma\colon \Odot^2S'\to\fh$ with $\alpha+\beta+\gamma\in\ssZ^{2,2}(\fa_-;\fa)$, $i_*[\alpha+\beta+\gamma]=0$ if and only if $i_*(\alpha+\beta+\gamma)=\partial\lambda$ where $\lambda:V\to\fso(V)$. Then $\gamma(s,s)=-\lambda(\kappa(s,s))$ for all $s\in S'$; in particular, $\lambda(V)=\lambda(\comm{S'}{S'})\subseteq\fh$, or equivalently $\lambda=i_*\lambda'$ for some $\lambda'\in\ssC^{2,1}(\fa_-;\fa)=\Hom(V,\fh)$. It follows that $\alpha+\beta+\gamma=\partial\lambda'$, so $[\alpha+\beta+\gamma]=0$, showing that $i_*\colon \ssH^{2,2}(\fa_-;\fa)\to \ssH^{2,2}(\fa_-;\fs)$ is injective.
\end{proof}

By the same argument as in the proof of Lemma~\ref{lemma:H22-full-flat-model}, in the highly supersymmetric case any cohomology class in $\ssH^{2,2}(\fa_-;\fs)$ has a unique representative in the space of normalised cocycles
\[
	\cH^{2,2}(\fa_-)
	= \bigl\{\beta+\gamma \in \ssZ^{2,2}(\fa_-;\fs)
		\mid
			\beta\in\Hom(V\otimes S',S),\,
			\gamma\in\Hom\bigl(\Odot^2S',\fso(V)\bigr)\bigr\},
\]
and we have a splitting of $\fh$-modules $\ssZ^{2,2}(\fa_-;\fs)=\cH^{2,2}(\fa_-)\oplus\ssB^{2,2}(\fa_-;\fs)$, whence as $\fh$-modules $\ssH^{2,2}(\fa_-;\fs) \cong \cH^{2,2}(\fa_-)$. We thus have an exact sequence at the level of cocycles corresponding to the bottom row of the diagram in the lemma above
\[
\begin{tikzcd}
	0 \arrow[r]
	& \cK^{2,2}(\fa_-) \arrow[r]
	& \cH^{2,2} \arrow[r, "i^*"]
	& \cH^{2,2}(\fa_-).
\end{tikzcd}
\]
On the other hand, any cohomology class in $\ssH^{2,2}(\fa_-;\fa)$ defines a unique normalised cocycle in~$\cH^{2,2}(\fa_-)$ by pushing forward to $\ssH^{2,2}(\fa_-;\fs)$ and then identifying the normalised cocycle. Since~$i_*$ is injective in cohomology, if a cocycle in $\cH^{2,2}(\fa_-)$ is obtained from a class in $\ssH^{2,2}(\fa_-;\fa)$, that class is unique. We can be much more explicit about this: if $\alpha+\beta+\gamma\in\ssZ^{2,2}(\fa_-;\fa)$ then there exists a unique element $\lambda\in\ssC^{2,2}(\fa_-;\fs)=\Hom(V,\fso(V))$ such that
\begin{equation}\label{eq:alpha-lambda}
	\alpha(v,w) = \lambda(v)w - \lambda(w)v;
\end{equation}
the components $\betatilde\colon V\otimes S' \to S$ and $\gammatilde\colon \Odot^2 S' \to \fso(V)$ of the unique normalised cocycle $\betatilde+\gammatilde\in\cH^{2,2}(\fa_-)$ such that $i_*[\alpha+\beta+\gamma]=[\betatilde+\gammatilde]$ are then given by
\begin{align}
	\betatilde(v,s) = \beta(v,s) - \lambda(v)\cdot s,\qquad
	\gammatilde(s,s) = \gamma(s,s) + \lambda(\kappa_s), \label{eq:gamma-tilde}
\end{align}
for $v\in V$, $s\in S'$. The maps $\lambda$, $\betatilde$, $\gammatilde$ are unique for a given cocycle $\alpha+\beta+\gamma$, while a different choice of representative for the same cohomology class corresponds to shifting $\lambda$ by a map $V\to\fh$. We note that
$i_*(\alpha+\beta+\gamma) = \betatilde+\gammatilde+\partial\lambda$,
and that we have the constraints
$\betatilde(v,s) + \lambda(v)\cdot s \in S'$, $
	\gammatilde(s,s) - \lambda(\kappa_s) \in \fh$
for all $v\in V$, $s\in S'$. Finally, observe that $[\alpha+\beta+\gamma]$ is $\fh$-invariant if and only if $\betatilde+\gammatilde\in\cH^{2,2}(\fa_-)$ is (which is the case if and only if $\betatilde$ is $\fh$-invariant). If this is the case, the discussion following Lemma~\ref{lemma:invt-spencer-cocycle} tells us that the map $\delta\colon \fh\otimes V\to\fh$ defined by
\begin{equation}\label{eq:delta-lambda}
	\delta(A,v) = \comm{A}{\lambda(v)} - \lambda(Av)
\end{equation}
is the unique such map so that $\alpha+\beta+\gamma+\delta\in\ssZ^{2,2}(\fa;\fa)$.

The discussion above allows us to reparametrise the deformed brackets \eqref{eq:def-brackets-general} in the highly supersymmetric case.

\begin{Proposition}\label{prop:filtered-def-high-susy}
	If $\fatilde$ is a highly supersymmetric filtered deformation of a graded subalgebra~${\fa=V\oplus S'\oplus \fh}$ of $\fs$, then $\fatilde$ has a presentation of the form
	\begin{gather}
		\comm{A}{B} = \underbrace{AB-BA}_\fh,\qquad
		\comm{s}{s} = \underbrace{\kappa_s}_V + \underbrace{\gammatilde(s,s) - \lambda(\kappa_s)}_\fh,\nonumber\\
		\comm{A}{v} = \underbrace{Av}_V + \underbrace{\comm{A}{\lambda(v)} - \lambda(Av)}_\fh,\qquad
		\comm{v}{s} = \underbrace{\betatilde(v,s) + \lambda(v)\cdot s}_{S'},\qquad
		\comm{A}{s} = \underbrace{A\cdot s}_{S'},\nonumber
	\\
		\comm{v}{w} = \underbrace{\lambda(v)w-\lambda(w)v}_V + \underbrace{\thetatilde(v,w) - \lambda(\lambda(v)w-\lambda(w)v) + \comm{\lambda(v)}{\lambda(w)}}_\fh\label{eq:def-brackets-V}
	\end{gather}
	for $A,B\in\fh$, $v,w\in V$, $s\in S'$, where the deforming maps are
	\[
\betatilde\colon\ V\otimes S' \to S,\qquad
		\gammatilde\colon\ \Odot^2 S' \to \fso(V),\qquad
		\lambda\colon\ V\to\fso(V),\qquad
		\thetatilde\colon\ \Wedge^2 V \to \fso(V),
	\]
	where $\betatilde+\gammatilde\in\cH^{2,2}(\fa_-)^\fh$, $\thetatilde$ is $\fh$-invariant and the remaining Jacobi identities are equivalent to
	\begin{gather}
		\thetatilde(\kappa_s,v) + 2\gammatilde\bigl(s,\betatilde(v,s)+\lambda(v)\cdot s\bigr) - \comm{\lambda(v)}{\gammatilde(s,s)} = 0, \label{eq:tilde-theta-kappa}\\
		\thetatilde(v,w)\cdot s
		= \betatilde\bigl(v,\betatilde(w,s)+\lambda(w)\cdot s\bigr)
			- \lambda(w)\cdot \betatilde(v,s) + \betatilde(\lambda(w)v,s)\nonumber\\
			\phantom{\thetatilde(v,w)\cdot s
		=}{} - \betatilde\bigl(w,\betatilde(v,s)+\lambda(v)\cdot s\bigr)
			+ \lambda(v)\cdot \betatilde(w,s) - \betatilde(\lambda(v)w,s),\label{eq:tilde-theta-lambda}\\
		\thetatilde(u,v)w + \thetatilde(v,w)u + \thetatilde(w,u)v =0,\label{eq:tilde-theta-bianchi}\\
		\bigl(\lambda(u)\cdot\thetatilde\bigr)(v,w) + \bigl(\lambda(v)\cdot\thetatilde\bigr)(w,u) +			\bigl(\lambda(w)\cdot\thetatilde\bigr)(u,v) = 0 \label{eq:tilde-theta-lambda-bianchi}
	\end{gather}
	for all $u,v,w\in V$ and $s\in S'$.
	Two deformations $\fatilde$, $\fatilde'$ described in this way are isomorphic as filtered Lie superalgebras if and only if they determine the same element $\betatilde+\gammatilde\in\cH^{2,2}(\fa_-)^\fh$.
\end{Proposition}

\begin{proof}
	We begin with a presentation of the form \eqref{eq:def-brackets-general}, so that $\alpha+\beta+\gamma$ is a $(2,2)$-Spencer cocycle for $\fa$ with $\fh$-invariant cohomology class and then define $\lambda,\betatilde,\gammatilde$ by equations \eqref{eq:alpha-lambda}--\eqref{eq:gamma-tilde}; $\delta$ must then be given by \eqref{eq:delta-lambda}. The brackets then take the claimed form if we define
$\thetatilde(v,w) = \theta(v,w) + \lambda(\lambda(v)w-\lambda(w)v) - \comm{\lambda(v)}{\lambda(w)}$.
	It remains to check the Jacobi identities~\eqref{eq:jacobi-022}--\eqref{eq:jacobi-222b} in our new variables. Substituting the definition of $\delta$ into the first identity \eqref{eq:jacobi-022} gives us $A\cdot \thetatilde=0$ for all $A\in\fh$. Next, we find that~\eqref{eq:jacobi-112} implies \eqref{eq:tilde-theta-kappa}
	\begin{gather*}
		\theta(\kappa_s,v) + \delta(\gamma(s,s),v) + 2\gamma(s,\beta(v,s)) \\
		\qquad= \theta(\kappa_s,v) + \comm{\gammatilde(s,s)}{\lambda(v)} - \lambda(\gammatilde(s,s)v) - \comm{\lambda(\kappa_s)}{\lambda(v)} + \lambda(\lambda(\kappa_s)v) \\
		\phantom{\qquad=}{}+ 2\gammatilde(s,\betatilde(v,s)+\lambda(v)\cdot s) - 2\lambda(\kappa(s,\betatilde(v,s)+\lambda(v)\cdot s)) \\
	\qquad	= \theta(\kappa_s,v) + \lambda\qty(\lambda(\kappa_s)v - 2\kappa(s,\lambda(v)\cdot s)) - \comm{\lambda(\kappa_s)}{\lambda(v)} + 2\gammatilde(s,\betatilde(v,s)+\lambda(v)\cdot s)\\
			\phantom{\qquad=}{}- \comm{\lambda(v)}{\gammatilde(s,s)}
			- \lambda(\kappa(s,s)v + 2\kappa(s,\betatilde(v,s)))\\
		\qquad= \thetatilde(\kappa_s,v) + 2\gammatilde(s,\betatilde(v,s)+\lambda(v)\cdot s) - \comm{\lambda(v)}{\gammatilde(s,s)},
	\end{gather*}
	where in the last line we have used the definition of $\thetatilde$, a cocycle condition for $\betatilde+\gammatilde$ and $\fso(V)$-invariance of $\kappa$. Substitution of the change of variables followed by some trivial rearrangement shows that \eqref{eq:jacobi-122} and \eqref{eq:jacobi-222a} are equivalent to \eqref{eq:tilde-theta-lambda} and \eqref{eq:tilde-theta-bianchi} respectively. Finally, \eqref{eq:jacobi-222b} gives~\smash{$\bigl(\lambda(u)\cdot\thetatilde\bigr)(v,w)
			+ \lambda\bigl(\thetatilde(u,v)w\bigr)
	+ \text{cyclic perms.}= 0$},
	but then substitution of \eqref{eq:tilde-theta-bianchi} yields~\eqref{eq:tilde-theta-lambda-bianchi}. The final claim follows by Proposition~\ref{prop:highly-susy-filtered-def} since the cocycle $\betatilde+\gammatilde\in\cH^{2,2}(\fa_-)^\fh$ is uniquely determined by the cohomology class $[\alpha+\beta+\gamma]\in\ssH^{2,2}(\fa_-;\fa)$.
\end{proof}

\subsection{Realisable subdeformations}
\label{sec:realisable-subdef}

We have just seen that if $\alpha+\beta+\gamma\in\ssZ^{2,2}(\fa_-;\fa)$ is the Spencer cocycle associated to some presentation of a filtered deformation of a highly supersymmetric graded subalgebra $\fa\subseteq\fs$, then there is an $\fh$-invariant cocycle $\betatilde+\gammatilde\in\ssZ^{2,2}(\fa_-;\fs)$ such that $i_*[\alpha+\beta+\gamma]=[\betatilde+\gammatilde]$. We will now restrict our attention to cocycles obeying a stronger version of this and the corresponding filtered subdeformations of $\fs$ (where they exist). Essentially, we demand that $\betatilde+\gammatilde$ be the pull-back of some $\fh$-invariant element $\betahat+\gammahat\in\cH^{2,2}$.

\subsubsection{Admissibility}

\begin{Definition}\label{def:admissible-class}
	A cohomology class $[\alpha+\beta+\gamma]\in \ssH^{2,2}(\fa_-;\fa)$ for a highly supersymmetric graded subalgebra $\fa$ of $\fs$ is \emph{admissible} if we have
$i_*[\alpha+\beta+\gamma]=i^*\big[\betahat+\gammahat\big]
$
	in $\ssH^{2,2}(\fa_-;\fs)$ for some $\fh$-invariant cocycle \smash{$\betahat+\gammahat\in \cH^{2,2}$}.
	
	A cocycle in $\ssZ^{2,2}(\fa_-;\fa)$ is \emph{admissible} if its cohomology class is admissible.
\end{Definition}

Note that an admissible cohomology class is in particular $\fh$-invariant. Unpacking the definition a little, for admissible cocycles we have
\begin{equation}\label{eq:admiss-cocycle-param}
	i_*(\alpha + \beta + \gamma) = i^*(\betahat + \gammahat) + \partial\lambda
\end{equation}
for a unique $\lambda\in\ssZ^{2,2}(\fa_-;\fs)=\Hom(V,\fso(V))$ and $\betahat+\gammahat\in\bigl(\cH^{2,2}\bigr)^\fh$ determined uniquely up to elements in $\cK^{2,2}(\fa_-)$. If we fix an admissible cohomology class but not a particular representative, $\lambda$ is only determined up to elements in $\ssZ^{2,2}(\fa_-;\fa)=\Hom(V,\fh)$.

Recalling that a highly supersymmetric subdeformation of $\fs$ uniquely determines an element of $\ssH^{2,2}(\fa_-;\fa)^\fh$ by Proposition~\ref{prop:highly-susy-filtered-def}, we make the following definition.

\begin{Definition}[geometric realisability]\label{def:geom-real}
	A highly supersymmetric filtered subdeformation $\fatilde$ of $\fs$ is \emph{$($geometrically$)$ realisable} if the associated cohomology class in $\ssH^{2,2}(\fa_-;\fa)^\fh$ is admissible.
\end{Definition}

The definitions above are inspired by, although \emph{not} equivalent to \cite[Definition~9]{Figueroa-OFarrill2017_1} in the context of 11-dimensional supergravity. There, one finds that $\cH^{2,2}\cong \Wedge^4V$, and it is the 4-form~$\varphi$ associated to $\betahat+\gammahat$ which is called admissible if the conditions above are met \emph{and} an additional algebraic condition related to the fact that $\varphi$ should considered to be the value at a point of a~closed 4-form field strength on some supergravity background is imposed. Our definition here is weaker because it is not clear \textit{a priori} how this condition should be generalised.

Let us make some particular remarks about these definitions in the maximally supersymmetric case, where $\fa_-=\fs_-$. Here, the map $i^*\colon\ssC^{2,2}(\fs_-;\fs)\to\ssC^{2,2}(\fa_-;\fs)$ whose induced map in cohomology appears in the definition of admissibility, is the identity map. An admissible cocycle for $\fa$ is then simply a cocycle in $\ssZ^{2,2}(\fa_-;\fa)$ whose cohomology class is $\fa_0$-invariant. It follows that \emph{all} maximally supersymmetric filtered deformations of $\fs$ are realisable.

We will justify the term ``geometrically realisable'' in Section~\ref{sec:high-susy-spin-mfld-homogeneous} (see Theorem~\ref{thm:homog-spin-mfld}). Our goal in this section will be to determine the conditions under which an admissible cohomology class of a highly supersymmetric graded subalgebra $\fa\subseteq\fs$ is \emph{integrable}, in the sense that there exists a corresponding geometrically realisable filtered deformation $\fatilde$. We first record the following. Note that the homogeneity theorem (see Theorem~\ref{thm:homogeneity}) is required for this result as well as for many in Section~\ref{sec:integrability}.

\begin{Lemma}\label{lemma:admiss-props}
	Let $\alpha+\beta+\gamma$ be an admissible cocycle parametrised as in \eqref{eq:admiss-cocycle-param}. Then
	\begin{itemize}\itemsep=0pt
		\item for all $v\in V$, $s\in S'$, $\betahat(v,s)+\lambda(v)\cdot s \in S'$;
		\item for all $s,s'\in S'$, $\gammahat(s,s')-\lambda(\kappa(s,s'))\in\fh$;
		\item for all $A\in\fh$, $v\in V$, $\comm{A}{\lambda(v)}-\lambda(Av)\in\fh$.
	\end{itemize}
\end{Lemma}

\begin{proof}
	The first two claims follow directly from the definition. For the third, by Theorem~\ref{thm:homogeneity} we can take $v=\kappa_s$ for some $s\in S'$. Recalling that the action of $\fh$ preserves $S'$ and using~$\fso(V)$-equivariance of $\kappa$ and $\fh$-invariance of $\gammahat$, we have
	\begin{align*}
		\comm{A}{\lambda(\kappa_s)} - \lambda(A\kappa_s)
		 &= \comm{A}{\lambda(\kappa_s)} - 2\lambda(\kappa(As,s))\\
		 &= \comm{A}{\lambda(\kappa_s)-\gamma(s,s)} + 2(\gammahat(As,s)-\lambda(\kappa(As,s))),
	\end{align*}
	and both of the terms in the last expression lie in $\fh$ by the second claim of the lemma.
\end{proof}

\subsubsection{Integrability}
\label{sec:integrability}

Suppose we have an admissible cocycle parametrised as in \eqref{eq:admiss-cocycle-param} and take $\betatilde=i^*\betahat$, $\gammatilde=i^*\gammahat$. We wish to know whether it is possible to define a (geometrically realisable) filtered subdeformation as in Proposition~\ref{prop:filtered-def-high-susy} using these maps in the definition of the bracket \eqref{eq:def-brackets-V}. By the above lemma, it remains only to check whether there exists an $\fh$-invariant map $\thetatilde\colon \Wedge^2V\to \fso(V)$ such that
\begin{equation}\label{eq:theta-lambda-in-h}
	\thetatilde(v,w) - \lambda(\lambda(v)w-\lambda(w)v) + \comm{\lambda(v)}{\lambda(w)} \in \fh
\end{equation}
satisfying the equations \eqref{eq:tilde-theta-kappa} to \eqref{eq:tilde-theta-lambda-bianchi}. We first observe that equations~\eqref{eq:tilde-theta-kappa} and \eqref{eq:tilde-theta-lambda} simplify slightly
\begin{gather}
	\thetatilde(v,\kappa_s) = 2\gammahat\bigl(s,\betahat(v,s)\bigr) - \bigl(\lambda(v)\cdot\gammahat\bigr)(s,s),
	\label{eq:tilde-theta-def-1}\\
		\thetatilde(v,w)\cdot s
		= \betahat\bigl(v,\betahat(w,s)\bigr) - \betahat\bigl(w,\betahat(v,s)\bigr) + \bigl(\lambda(v)\cdot\betahat\bigr)(w,s) - \bigl(\lambda(w)\cdot\betahat\bigr)(v,s),
		\label{eq:tilde-theta-act-s}
\end{gather}
and by Theorem~\ref{thm:homogeneity}, either equation fully determines $\thetatilde$ (if it exists). In fact, $\thetatilde$ is uniquely determined by the admissible cohomology class: it depends only on $i^*\betahat$, $i^*\gammahat$ and $\lambda$, so a shift in~${\betahat+\gammahat}$ by an element of $\cK^{2,2}(\fa_-)$ does not change it, while by $\fh$-invariance of $\betahat$ and $\gammahat$, changing~$\lambda$ by addition of a map $V\to\fh$ also does not change $\thetatilde$. There may, however, be obstructions to the existence of $\thetatilde\colon\Wedge^2 V\to \fso(V)$ even for an admissible cocycle. For now, we will postpone considering this issue and show that a map $\thetatilde$ determined by the equations above automatically satisfies all of the other required properties.

\begin{Lemma}\label{lemma:theta-lambda-in-h}
	Any map $\thetatilde\colon\Wedge^2V\to \fso(V)$ satisfying equation~\eqref{eq:tilde-theta-def-1} also satisfies~\eqref{eq:theta-lambda-in-h}.
\end{Lemma}

\begin{proof}
	By Theorem~\ref{thm:homogeneity}, it is sufficient to show the result for $w=\kappa_s$ for some $s\in S'$; we have
	\begin{gather*}
		\thetatilde(v,\kappa_s) - \lambda(\lambda(v)\kappa_s-\lambda(\kappa_s)v) + \comm{\lambda(v)}{\lambda(\kappa_s)}\\
			\qquad= 2\big\{\qty(\gammahat-\lambda\circ\kappa)(s,\betahat(v,s)+\lambda(v)\cdot s)\big\}\\
				\phantom{\qquad= }{} + \qty{\comm{(\gammahat-\lambda\circ\kappa)(s,s)}{\lambda(v)} - \lambda((\gammahat-\lambda\circ\kappa)(s,s)v)},
	\end{gather*}
	and since $(\gammahat-\lambda\circ\kappa)|_{\Odot^2S'}$ takes values in $\fh$, Lemma~\ref{lemma:admiss-props} implies that both of the expressions within the braces lie in $\fh$.
\end{proof}

\begin{Lemma}\label{lemma:tilde-theta-def-2}
	Any map $\thetatilde\colon\Wedge^2 V\to \fso(V)$ satisfying equation~\eqref{eq:tilde-theta-act-s} also satisfies
	\begin{gather*}
		\thetatilde(v,w)\kappa_s
			= 2\gammahat\bigl(\betahat(v,s),s\bigr)w
				- 2\gammahat\bigl(\betahat(w,s),s\bigr)v- \bigl(\lambda(v)\cdot\gammahat\bigr)(s,s)w
				+ \qty(\lambda(w)\cdot\gammahat)(s,s)v\label{eq:tilde-theta-def-2}
	\end{gather*}
	for all $v,w\in V$, $s\in S'$, and this equation determines $\thetatilde$ uniquely.
\end{Lemma}

\begin{proof}
	First, it will be useful to depolarise one of the cocycle conditions for $\betahat+\gammahat$
	\begin{equation}\label{eq:22-spencer-cocycle-vss-depol}
		2\kappa\bigl(\betahat(v,s),s'\bigr) + 2\kappa\bigl(\betahat(v,s'),s\bigr) + 2\gammahat(s,s')v = 0.
	\end{equation}
	Using $\fso(V)$-invariance of $\kappa$, equation~\eqref{eq:tilde-theta-act-s} gives
	\begin{align*}
		\thetatilde(v,w)\kappa_s
		= 2\kappa\bigl(\thetatilde(v,w)\cdot s,s\bigr)
		={}& 2\kappa\bigl(\betahat(v,\betahat(w,s)\bigr),s)
			+ 2\kappa\bigl(\lambda(v)\cdot\betahat(w,s),s\bigr)\\
			& - 2\kappa\bigl(\betahat(\lambda(v)w,s),s\bigr)
			- 2\kappa\bigl(\betahat(w,\lambda(v)\cdot s),s\bigr)
			- (v \leftrightarrow w).
	\end{align*}
	We now use equation~\eqref{eq:22-spencer-cocycle-vss-depol} with $s'=s,\lambda(v)\cdot s$ and $\betahat(w,s) $ to replace some of the terms above and then use $\fso(V)$-invariance of $\kappa$ again to conclude that
	\begin{align*}
		\thetatilde(v,w)\kappa_s
		={}& - 2\gammahat\bigl(\betahat(w,s),s\bigr)v
			- 2\kappa\bigl(\betahat(v,s),\betahat(w,s)\bigr) + 2\kappa\bigl(\lambda(v)\cdot\betahat(w,s),s\bigr)
\\
			&+ \gammahat(\lambda(v)\cdot s,s)v + \gammahat(s,s)(\lambda(v)w)
			+ 2\kappa\bigl(\betahat(w,s),\lambda(v)\cdot s\bigr)
			- (v \leftrightarrow w)\\
		={}&2\gammahat\bigl(\betahat(v,s),s\bigr)w
			+ \lambda(v)\kappa\bigl(\betahat(w,s),s\bigr)\\
			& + \gamma(s,s)(\lambda(v)w)
			+ \gammahat(\lambda(v)\cdot s,s)w - (v \leftrightarrow w),
	\end{align*}
	which is the desired expression. In the final equality, we have used the skew-symmetry of the expression in $v$ and $w$. By Theorem~\ref{thm:homogeneity}, this defines $\thetatilde$ uniquely.
\end{proof}

\begin{Lemma}\label{lemma:quadr-jacobi-dependency}
	For any map $\thetatilde\colon\Wedge^2 V\to\fso(V)$, we have the following:
	\begin{itemize}\itemsep=0pt
		\item Either equation~\eqref{eq:tilde-theta-def-1} or \eqref{eq:tilde-theta-act-s} implies that $\thetatilde$ is $\fh$-invariant.
		\item Assuming equation~\eqref{eq:tilde-theta-def-1}, equation~\eqref{eq:tilde-theta-def-2} is equivalent to the algebraic Bianchi identity~\eqref{eq:tilde-theta-bianchi}. In particular, equations \eqref{eq:tilde-theta-def-1} and \eqref{eq:tilde-theta-act-s} together imply the Bianchi identity.
		\item Equations~\eqref{eq:tilde-theta-def-1} and \eqref{eq:tilde-theta-act-s} together imply \eqref{eq:tilde-theta-lambda-bianchi}.
	\end{itemize}
\end{Lemma}

\begin{proof}
	The first claim follows by some fairly straightforward but long manipulations which use $\fso(V)$-invariance of $\kappa$, $\fh$-invariance of $\betahat$ and $\gammahat$ and the fact that $\lambda(Av)-\comm{A}{\lambda(v)}\in\fh$ for all~${A\in\fh}$. For the second claim, let us assume equation \eqref{eq:tilde-theta-def-1} holds. Then equation~\eqref{eq:tilde-theta-def-2} gives~${\thetatilde(v,w)\kappa_s = \thetatilde(v,\kappa_s)w - \thetatilde(w,\kappa_s)v}
$
	for all $v,w\in V$, $s\in S'$, which by the homogeneity theorem is the Bianchi identity. Conversely, if we assume the Bianchi identity, the above equation along with \eqref{eq:tilde-theta-def-1} gives us \eqref{eq:tilde-theta-def-2}.
	
	For the final claim let us assume that $\thetatilde$ is given by \eqref{eq:tilde-theta-def-1} and note that by Theorem~\ref{thm:homogeneity} we can assume without loss of generality that $w=\kappa_s$ for arbitrary $s\in S'$. Then we can use $\fso(V)$-invariance of $\kappa$ to compute
	\begin{align*}
		\bigl(\lambda(u)\cdot\thetatilde\bigr)(v,\kappa_s) + \bigl(\lambda(v)\cdot\thetatilde\bigr)(\kappa_s,u)
			={}& \bigl((\lambda(\lambda(u)v-\lambda(v)u)
				- \comm{\lambda(u)}{\lambda(v)})\cdot\gammahat\bigr)(s,s) \\
				& + 2\gammahat\bigl(s,(\lambda(u)\cdot\betahat)(v,s)
				- (\lambda(v)\cdot\betahat)(u,s)\bigr)\\
				& + 2(\lambda(u)\cdot\gammahat)\bigl(s,\betahat(v,s)\bigr)
				- 2(\lambda(v)\cdot\gammahat)\bigl(s,\betahat(u,s)\bigr).
	\end{align*}
	On the other hand, we recall that $\gammahat(s,s)-\lambda(\kappa_s)\in\fh$, and that $\thetatilde$ is $\fh$-invariant by the first part of the lemma, whence $\lambda(\kappa_s)\cdot\thetatilde=\gammahat(s,s)\cdot\thetatilde$. Using \eqref{eq:tilde-theta-def-1} and a cocycle condition, we find\footnote{Actually, it's necessarily to depolarise \eqref{eq:tilde-theta-def-1} in order to evaluate expressions like $\thetatilde(v,\kappa(s,s'))$ for $v\in V$ $s,s'\in S'$.}
	\begin{align*}
		\bigl(\lambda(\kappa_s)\cdot\thetatilde\bigr)(u,v) ={}& [\gammahat(s,s),\thetatilde(u,w)]
			+ 2\gammahat\bigl(s,\betahat(u,\betahat(v,s))-\betahat(v,\betahat(u,s))\bigr)\\
			& - 2(\lambda(u)\cdot\gammahat)\bigl(s,\betahat(v,s)\bigr)
			+ 2(\lambda(v)\cdot\gammahat)\bigl(s,\betahat(u,s)\bigr).
	\end{align*}
	Putting the two calculations together and then using \eqref{eq:tilde-theta-act-s}, we find
	\begin{gather*}
		\bigl(\lambda(u)\cdot\thetatilde\bigr)(v,\kappa_s) + \bigl(\lambda(v)\cdot\thetatilde\bigr)(\kappa_s,u) + \bigl(\lambda(\kappa_s)\cdot\thetatilde\bigr)(u,v)\\
		\qquad = \bigl((\lambda(\lambda(u)v-\lambda(v)u) - \comm{\lambda(u)}{\lambda(v)})\cdot\gammahat\bigr)(s,s)
		-\big[\thetatilde(u,w),\gammahat(s,s)\big] \\
		\phantom{\qquad =}{} + 2\gammahat\bigl(s,\betahat\bigl(u,\betahat(v,s)\bigr)-\betahat\bigl(v,\betahat(u,s)\bigr)\bigr)+ 2\gammahat\bigl(s,\bigl(\lambda(u)\cdot\betahat\bigr)(v,s) - \bigl(\lambda(v)\cdot\betahat\bigr)(u,s)\bigr)\\
		\qquad= -\bigl(\bigl(\thetatilde(u,w) + \comm{\lambda(u)}{\lambda(v)} - \lambda(\lambda(u)v-\lambda(v)u)\bigr)\cdot\gammahat\bigr)(s,s).
	\end{gather*}
	The final line vanishes by Lemma~\ref{lemma:theta-lambda-in-h} and $\fh$-invariance of $\gammahat$.
\end{proof}

Thus the problem of integrating an admissible cocycle reduces to the question of whether a~map satisfying the equations \eqref{eq:tilde-theta-def-1} and \eqref{eq:tilde-theta-act-s} exists. Let us first consider equation \eqref{eq:tilde-theta-def-1}. We would like to interpret this as the definition of a map $\thetatilde\colon \Wedge^2V\to \fso(V)$, but \textit{a priori} we only have that
\begin{equation}\label{eq:Theta-def}
	\Theta(v,s,s) = 2\gammahat\bigl(s,\betahat(v,s)\bigr) - \bigl(\lambda(v)\cdot\gammahat\bigr)(s,s)
\end{equation}
defines a map $\Theta\colon V\otimes \Odot^2 S' \to \fso(V)$. So we require  that $\Theta$ factorises as follows
\[
\begin{tikzcd}
	V\otimes \Odot^2 S' \ar[rr,"\Theta"]\ar[rd,"\Id\otimes \kappa"'] & & \fso(V)\\
		& V \otimes V \ar[ru,"\thetatilde"'] &
\end{tikzcd}
\]
and that $\thetatilde$ is skew-symmetric. We define the \emph{Dirac kernel} of $S'$ by $\fD=\ker\kappa|_{\Odot^2S'}$ and note that $\Theta$ factorises as above if and only if it annihilates the Dirac kernel, i.e., $\Theta(V\otimes\fD)=0$. We then find that $\thetatilde$ is automatically skew-symmetric.

\begin{Lemma}\label{lemma:theta-well-def}
	Suppose that $\Theta$ annihilates the Dirac kernel and thus factors as described above. Then equation~\eqref{eq:tilde-theta-def-1} defines a map $\thetatilde\colon\Wedge^2V\to \fso(V)$.
\end{Lemma}

\begin{proof}
	The factorisation of $\Theta$ gives us a map $\thetatilde\colon V\otimes V \to \fso(V)$ satisfying \eqref{eq:tilde-theta-def-1}. To show that $\thetatilde$ is alternating, it is sufficient to show that $\thetatilde(\kappa_s,\kappa_s)=0$ for all $s\in S'$. Using the fact that $\gammahat(s,s)-\lambda(\kappa_s)\in\fh$ for $s\in S'$ and $\fh$-invariance of $\gammahat$, we have $\lambda(\kappa_s)\cdot\gammahat = \gammahat(s,s)\cdot\gammahat$, so
	\begin{align*}
		\thetatilde(\kappa_s,\kappa_s)
		&= 2\gammahat\bigl(s,\betahat(\kappa_s,s)\bigr) - \bigl(\gammahat(s,s)\cdot\gammahat\bigr)(s,s)\\
		& = 2\gammahat\bigl(s,\betahat(\kappa_s,s) + \gammahat(s,s)\cdot s\bigr) - \comm{\gammahat(s,s)}{\gammahat(s,s)}
		= 0,
	\end{align*}
	where we have used a cocycle condition in the last equality.
\end{proof}

\begin{Definition}\label{def:integrable-cocycle}
	Let $[\alpha+\beta+\gamma]$ be an admissible cohomology class and, recalling that it does not depend on the choice of cocycle representative, let $\Theta\colon V\otimes\Odot^2S'\to\fso(V) $ be the map defined by equation~\eqref{eq:Theta-def}. We say that $[\alpha+\beta+\gamma]$ is \emph{integrable} if
	\begin{enumerate}\itemsep=0pt
	\item $\Theta$ annihilates the Dirac kernel $\fD$ of $\Odot^2S'$,
	\item the map $\thetatilde\colon\Wedge^2 V\to \fso(V)$ defined by $\thetatilde\circ\kappa=\Theta$ satisfies	\eqref{eq:tilde-theta-act-s}.
	\end{enumerate}
	An admissible cocycle is \emph{integrable} if its cohomology class is.
\end{Definition}

By construction, we have the following refinement of Proposition~\ref{prop:filtered-def-high-susy}.

\begin{Theorem}[integration of admissible cocycles]\label{thm:admiss-coycle-integr}
	Let $\alpha+\beta+\gamma$ be an admissible, integrable Spencer $(2,2)$-cocycle for a highly supersymmetric graded subalgebra $\fa$ of $\fs$ with $i_*(\alpha+\beta+\gamma)=i^*\bigl(\betahat+\gammahat\bigr)+\partial\lambda$ and let $\thetatilde\colon\Wedge^2V\to\fso(V)$ be the map defined by equation~\eqref{eq:tilde-theta-def-1}. Then the brackets
	\begin{gather*}
			\comm{A}{B} = \underbrace{AB-BA}_\fh,\qquad
			\comm{s}{s} = \underbrace{\kappa_s}_V + \underbrace{\gammahat(s,s) - \lambda(\kappa_s)}_\fh,\\
			\comm{A}{v}= \underbrace{Av}_V + \underbrace{\comm{A}{\lambda(v)} - \lambda(Av)}_\fh,\qquad
			\comm{v}{s} = \underbrace{\betahat(v,s) + \lambda(v)\cdot s}_{S'},\qquad
			\comm{A}{s} = \underbrace{A\cdot s}_{S'},
	\\
		\comm{v}{w} = \underbrace{\lambda(v)w-\lambda(w)v}_V + \underbrace{\thetatilde(v,w) - \lambda(\lambda(v)w-\lambda(w)v) + \comm{\lambda(v)}{\lambda(w)}}_\fh \label{eq:admiss-cocycle-brackets}
	\end{gather*}
	define a geometrically realisable filtered deformation of $\fa$. Two such deformations $\fatilde$, $\fatilde'$ are isomorphic if and only if $\bigl(\betahat'+\gammahat'\bigr)-\bigl(\betahat+\gammahat\bigr)\in\cK^{2,2}(\fa_-)$.
\end{Theorem}

Where $\fa$ is a \emph{maximally} supersymmetric graded subalgebra of the (minimal) Poincar\'e superalgebra $\fs$ in 4, 5, 6 or 11 spacetime dimensions, there is no obstruction to an admissible cocycle being integrable \cite{Beckett2021,deMedeiros2016,deMedeiros2018,Figueroa-OFarrill2017}. We already observed that admissible cohomology classes are simply elements of $\ssH^{2,2}(\fs_-;\fa)^{\fa_0}$ if $\fa$ is highly supersymmetric, so \emph{any} such cohomology class defines a realisable filtered deformation in these cases. The question of whether obstructions exist in the highly supersymmetric (but sub-maximal) case is still open even in these cases.

\subsection[The classification problem for highly supersymmetric odd-generated filtered\\ subdeformations]{The classification problem for highly supersymmetric\\ odd-generated filtered subdeformations}
\label{sec:classification-odd-gen}

We will now outline a scheme for the classification of highly supersymmetric odd-generated realisable filtered subdeformations $\fatilde$ of $\fs$. Odd-generated here means that $\fatilde$ is generated by its odd subspace; $\fatilde_{\overline{0}}=\comm{\fatilde_{\overline{1}}}{\fatilde_{\overline{1}}}$. We note that this is very close to a classification of \emph{all} highly supersymmetric realisable filtered subdeformations, since if $\fatilde$ is a general such subdeformation, the ideal subalgebra $\fatilde'=\comm{\fatilde_{\overline{1}}}{\fatilde_{\overline{1}}}\oplus\fatilde_{\overline{1}}\subseteq\fatilde$ (which we will call the \emph{odd-generated ideal}) is odd-generated realisable, and for the associated graded subalgebras of $\fs$ we have $\fa_{-2}'=\fa_{-2}$, ${\fa_{-1}'=\fa_{-1}}$ and $\fa_0'\subseteq\fa_0$. Thus only the zero-graded part is not determined by the odd-generated subalgebra, and by inspection of \eqref{eq:admiss-cocycle-brackets}, it is clear that the deformation maps for $\fatilde$ are also fully determined by those of $\fatilde'$. Thus the only data which distinguishes $\fatilde$ from its odd-generated ideal is the choice of some subalgebra $\fa_0\subseteq\fso(V)$ which contains $\fa_0'$ and preserves both $S'$ and~${\beta+\gamma\in\cH^{2,2}}$.

We generalise the results of \cite[Section~4]{Figueroa-OFarrill2017_1} and \cite[Section~5]{Santi2022} which treat the 11-dimensional case.

\subsubsection{Isomorphisms and embeddings}

We first generalise Definitions 5 and 11 from \cite{Figueroa-OFarrill2017_1}. For this, consider $\fs=V\oplus S\oplus\fso(V)$ under the action of its graded automorphism group
\[
	\GrAut\fs = \qty{ g\in\GL(\fs)\mid g\cdot\fs_i = \fs_i \text{ and } \forall X,Y\in\fs, g\cdot\comm{X}{Y}=\comm{g\cdot X}{g\cdot Y}}.
\]
We note that there is a natural action of $\Spin_0(V)$ as ``inner'' graded automorphisms, in the sense that this group integrates the Lie algebra of inner graded derivations induced by the action of $\fs_0=\fso(V)$, and this gives us an embedding $\Spin_0(V)\hookrightarrow \GrAut\fs$.
If the Dirac current (which by definition is $\fso(V)$-equivariant) is equivariant under the whole spin group, then this inner automorphism group can be enlarged to $\Spin(V)$. The group of graded automorphisms also contains a one-dimensional subgroup corresponding to dilatations and may additionally contain an ``$R$-symmetry'' group (for a recent review and classification see \cite{Gall2021}.

The role played by $\GrAut\fs$ in our formalism is played by $\Spin(V)$ in \cite{Figueroa-OFarrill2017_1}, which is a somewhat more natural choice in that more concrete context. We note that with the choices made in loc.\ cit., $\GrAut\fs$ contains the whole spin group, while the $R$-symmetry group is isomorphic to~$\ZZ_2$ and embeds into $\GrAut\fs$ as the centre of the isomorphic image of $\Spin(V)$. In this case, $\Gr\Aut\fs\cong\Spin(V)\times\RR_{>0}$, where the latter factor is a multiplicative group acting as dilatations.

Note that if $\fa$ is a graded subalgebra of $\fs$ then $g\cdot\fa$ is also a graded subalgebra for all $g\in G$, with left multiplication by $g$ defining a graded isomorphism $\fa\cong g\cdot\fa$.

\begin{Definition}\label{def:subdef-isom-embedding}
	Two filtered subdeformations $\fatilde$, $\fatilde'$ of the Poincar\'e superalgebra $\fs$ are \emph{isomorphic} if there exists a (strict) filtered isomorphism of Lie superalgebras $\Phi\colon \fatilde\to\fatilde'$ such that the associated graded morphism $\Gr\Phi\colon \fa\to\fa'$ is given for all $x\in \fa$, by $\Gr\Phi(x)=g\cdot x$ for some~${g\in\GrAut\fs}$.
	
	An \emph{embedding} of subdeformations $\fatilde\hookrightarrow\fatilde'$ of $\fs$ is an injective (strict) filtered morphism of Lie superalgebras which is an isomorphism of subdeformations onto its image.
\end{Definition}

Note that two filtered subdeformations of $\fs$ being isomorphic \emph{as subdeformations} is stronger than them simply being isomorphic as filtered Lie superalgebras.

With this in mind, we can strengthen point (2) of Proposition~\ref{prop:highly-susy-filtered-def} by noting that two highly supersymmetric deformations of the same subalgebra $\fa\subseteq\fs$ inducing the same cohomology class in $\ssH^{2,2}(\fa_-;\fa)$ are isomorphic \emph{as subdeformations} and not just as filtered Lie superalgebras. Indeed, this follows immediately from the proof of the proposition, which also shows in particular that the associated graded automorphism of $\fa$ is the identity map.

\subsubsection{Kernels and sections of the Dirac current}

Let $S'\subseteq S$ with $\dim S'>\frac{1}{2}\dim S$ and recall that by Theorem~\ref{thm:homogeneity}, $\kappa|_{\Odot^2S'}$ is surjective, so we have a short exact sequence of vector spaces
\[
\begin{tikzcd}
	0 \ar[r] & \fD \ar[r] & \Odot^2S' \ar[r,"\kappa|_{\Odot^2S'}"] & V \ar[r] & 0,
\end{tikzcd}
\]
where $\fD=\ker\kappa|_{\Odot^2S'}$ is the Dirac kernel of $S'$ as in Section~\ref{sec:integrability}. If $\fh\subseteq\fso(V)$ is a subalgebra preserving $S'$, then this is a sequence of $\fh$-modules. We define a \emph{section of the Dirac current $\kappa$} (understood to mean $\kappa$ restricted to $S'$) as a linear map $\Sigma\colon V\to\Odot^2S'$ which splits the above as a sequence of vector spaces, i.e., $\kappa\circ\Sigma=\Id_V$, giving a non-canonical decomposition of vector spaces (not of $\fh$-modules) $\Odot^2S'\simeq V\oplus\fD$. Two such sections differ by a map $V\to \fD$.

In the following, we recall from Section~\ref{sec:reparam-defs} the maps of cochains $i_*\colon\ssC^{\bullet,\bullet}(\fa_-;\fa) \to \!\ssC^{\bullet,\bullet}(\fa_-;\fs)$, $\phi\mapsto i\circ\phi$ and $i^*\colon\ssC^{\bullet,\bullet}(\fs_-;\fs) \to \ssC^{\bullet,\bullet}(\fa_-;\fs)$, $\phi\mapsto \phi\circ i$ (and the induced maps of cochains and cohomology classes) induced by the inclusion $i\colon\fa\hookrightarrow\fs$.

\begin{Lemma}\label{lemma:admissible-cocycles-sections}
	Let $\fa=V\oplus S'\oplus\fh$ be a highly supersymmetric graded subalgebra of $\fs$ and let $\betahat+\gammahat\in\cH^{2,2}$. Then there exists a cohomology class $[\alpha+\beta+\gamma]\in\ssH^{2,2}(\fa_-;\fa)$ such that $i_*[\alpha+\beta+\gamma]=i^*[\betahat+\gammahat]$ if and only if
	\begin{enumerate}\itemsep=0pt
	\item[$(1)$] $\gammahat(\fD)\subseteq\fh$,
	\item[$(2)$] for all $v\in V$, $\bigl(\imath_v\betahat+\gammahat(\Sigma(v))\bigr)(S')\subseteq S'$ for some $($hence any$)$ section $\Sigma$ of $\kappa$.
	\end{enumerate}	
	Moreover, if such a cohomology class exists, for a fixed section $\Sigma$ there exists a unique representative $\alpha+\beta+\gamma$ such that $i_*(\alpha+\beta+\gamma)=i^*\bigl(\betahat+\gammahat\bigr)+\partial\bigl(\gammahat\circ\Sigma\bigr)$.
\end{Lemma}

\begin{proof}
	Suppose $i_*[\alpha+\beta+\gamma]=i^*\big[\betahat+\gammahat\big]$, so that $i_*(\alpha+\beta+\gamma)=i^*\bigl(\betahat+\gammahat\bigr)+\partial\lambda$ for some unique $\lambda\in\ssC^{2,1}(\fa_-;\fs)=\Hom(V,\fso(V))$. On one hand, restricting to $\fD\subseteq\Odot^2S'$, we find that~${\gammahat|_{\fD}=\gamma|_{\fD}}$, which takes values in $\fh$. Now let $\Sigma\colon V\to \Odot^2 S'$ be a section of $\kappa$ and note that we have
$\gammahat\circ\Sigma - \lambda
			= \qty(\gammahat - \lambda\circ\kappa)\circ\Sigma
			= \gamma\circ\Sigma
$
	which also takes values in $\fh$. On the other hand, restricting to $V\otimes S'$, for all $v\in V$ and $s\in S'$ we have $\betahat(v,s)+\lambda(v)\cdot s = \beta(v,s)\in S'$, i.e., $\imath_v\betahat+\lambda(v)$ preserves $S'$, but then
	\[
		\imath_v\betahat + \gammahat(\Sigma(v)) = (\imath_v\betahat + \lambda(v)) + (\gammahat\circ\Sigma-\lambda)(v)
	\]
	also preserves $S'$ since the action of $\fh$ preserves $S'$.
	
	Conversely, suppose that $\gammahat(\fD)\subseteq\fh$ and $\imath_v\betahat+\gammahat(\Sigma(v))$ preserves $S'$ for all $v\in V$ for some (hence any\footnote{Since any two sections differ by a map $V\to\fD$, the conditions together imply that $\imath_v\betahat+\gammahat(\Sigma'(v))$ preserves~$S'$ for any other section $\Sigma'$ as well.})
	section $\Sigma$. Then define $\lambda\colon V\to\fso(V)$ by $\lambda=\gammahat\circ\Sigma$, giving us $(\gammahat-\lambda\circ\kappa)(\fD)=\gammahat(\fD)\in\fh$ and $(\gammahat-\lambda\circ\kappa)\circ\Sigma=\gammahat\circ\Sigma-\lambda=0$ which together imply that $(\gammahat-\lambda\circ\kappa)\bigl(\Odot^2S'\bigr)\subseteq \fh$. Thus we can define a cocycle $\alpha+\beta+\gamma\in\ssZ^{2,2}(\fa_-;\fa)$ with \smash{$i_*[\alpha+\beta+\gamma]=i^*\big[\betahat+\gammahat\big]$} by
	\begin{gather*}
		\alpha(v,w)	= \gammahat(\Sigma(v))w - \gammahat(\Sigma(w))(v),\qquad
		\beta(v,s)	= \betahat(v,s) + \gammahat(\Sigma(v))\cdot s,\\
		\gamma(s,s) = \gammahat(s,s) - \gammahat(\Sigma(\kappa_s)),
	\end{gather*}
	for $v,w\in V$, $s\in S'$. This also shows the final claim.
\end{proof}

Note that in the lemma above, a change of section $\Sigma$ does not change the cohomology class of a cocycle $i_*(\alpha+\beta+\gamma)=i^*\bigl(\betahat+\gammahat\bigr)+\partial\bigl(\gammahat\circ\Sigma\bigr)$. Indeed, if $\Sigma'$ is another section, then
\[
	i^*\bigl(\betahat+\gammahat\bigr)+\partial(\gammahat\circ\Sigma') - \bigl(i^*\bigl(\betahat+\gammahat\bigr)+\partial\bigl(\gammahat\circ\Sigma\bigr)\bigr)
	= \partial(\gammahat\circ(\Sigma'-\Sigma))
\]
and since $(\Sigma'-\Sigma)$ takes values in $\fD$, $\gammahat\circ(\Sigma'-\Sigma)$ takes values in $\fh$.

\subsubsection{Lie pairs, envelopes and odd-generated realisable subeformations}

We now introduce some technology for building subdeformations from normalised cocycles. Note that since all admissible cocycles will now be discussed only in terms of a fixed element of $\cH^{2,2}$, we no longer decorate the components of such elements with hats, trusting that this will not cause confusion.

\begin{Definition}[envelopes and Lie pairs]\label{def:envelope}
	Fix a subspace $S'\subseteq S$ with $\dim S'>\frac{1}{2}\dim S$ and~${\beta+\gamma\in\cH^{2,2}}$. We define the \emph{envelope of $(S',\beta+\gamma)$} as the subspace
$\fh_{(S',\beta+\gamma)}=\gamma(\fD) \subseteq \fso(V)$.
	We say that $(S',\beta+\gamma)$ is a \emph{Lie pair} if for all $A\in\fh_{(S',\beta+\gamma)}$
	\begin{enumerate}\itemsep=0pt
	\item $A\cdot\beta=0$,
	\item $A\cdot S' \subseteq S'$.
	\end{enumerate}
\end{Definition}

Our use of the term ``Lie pair'' follows our sources \cite{Figueroa-OFarrill2017_1,Santi2022} and is motivated by the first part of the following result, reconciling it with the Lie pairs defined in the appendices (see Definition~\ref{def:lie-pair-algs}).

\begin{Proposition}\label{prop:lie-pair-admiss-cocycle}
	If $(S',\beta+\gamma)$ is a Lie pair, then $\fh_{(S',\beta+\gamma)}$ is a Lie subalgebra of $\fso(V)$. If additionally $\imath_v\beta+\gamma(\Sigma(v))$ preserves $S'$ for some section $\Sigma$ of $\kappa$, then $i^*(\beta+\gamma)+\partial(\gamma\circ\Sigma)$ defines an admissible cocycle for the graded subalgebra $\fa_{(S',\beta+\gamma)}=V\oplus S'\oplus\fh_{(S',\beta+\gamma)}$ of $\fs$. Furthermore, if this cocycle is integrable, the associated realisable deformation $\fatilde_{(S',\beta+\gamma)}$ is odd-generated.
	
	Conversely, for any odd-generated realisable filtered subdeformation $\fatilde$ of $\fs$ there exists a Lie pair $(S',\beta+\gamma)$ such that $\Gr\fatilde\cong\fa_{(S',\beta+\gamma)}$.
\end{Proposition}

\begin{proof}
	Let $A,B\in\fh_{(S',\beta+\gamma)}$. Then there exists $\omega\in\fD$ such that $B=\gamma(\omega)$; note that $A\cdot\beta=0\implies A\cdot\gamma=0$, so $\comm{A}{B}=\comm{A}{\gamma(\omega)}=\gamma(A\cdot\omega)$. On the other hand, the action of $A$ preserves~$S'$, whence $A\cdot\omega\in\Odot^2S'$, and $\kappa(A\omega)=A\kappa(\omega)=0$, so $A\cdot\omega\in\fD$. Thus $\comm{A}{B}\in\fh_{(S',\beta+\gamma)}$, so the latter is indeed a subalgebra of $\fso(V)$. Since the action of $\fh_{(S',\beta+\gamma)}$ preserves $S'$, $\fa_{(S',\beta+\gamma)}$ is a~graded subalgebra of $\fs$.
	
	If $\imath_v\beta+\gamma(\Sigma(v))$ preserves $S'$ for some (hence any) section $\Sigma$, Lemma~\ref{lemma:admissible-cocycles-sections} tells us that $i^*(\beta+\gamma)+\partial(\gamma\circ\Sigma)$ is the image under $i_*$ of a $(2,2)$-cocycle of $\fa_{(S',\beta+\gamma)}$; since $\beta$ is $\fh_{(S',\beta+\gamma)}$-invariant, this cocycle is admissible. If it also integrable, we must show that $[S,S]=V\oplus\fh_{(S',\beta+\gamma)}$ where the bracket is the deformed one defined by $i^*(\beta+\gamma)+\partial(\gamma\circ\Sigma)$. The odd-odd bracket is nothing but the map $\kappa+\gamma-\gamma\circ\Sigma\circ\kappa\colon \Odot^2S'\to V\oplus\fh_{(S',\beta+\gamma)}$. Recall that the section gives~${\Odot^2S'=\fD\oplus\Sigma(V)}$ as a vector space, so since we have
	\[
		(\kappa+\gamma-\gamma\circ\Sigma\circ\kappa)(\fD) = \gamma(\fD) = \fh_{(S',\beta+\gamma)}
	\]
	by definition, and for $v\in V$,
	\[
		(\kappa+\gamma-\gamma\circ\Sigma\circ\kappa)(\Sigma(v))
			= v + \gamma(\Sigma(v)) - \gamma(\Sigma(v))
			= v,
	\]
	we have $(\kappa+\gamma-\gamma\circ\Sigma\circ\kappa)\bigl(\Odot^2S'\bigr)=V\oplus\fh_{(S',\beta+\gamma)}$ as required.
	
	Now assume that $\fatilde$ is an odd-generated realisable filtered subdeformation with associated graded $\fa=V\oplus S'\oplus\fh$; by Lemma~\ref{lemma:admissible-cocycles-sections}, the associated admissible cohomology class has a~re\-presentative $\alpha+\beta+\gamma$ satisfying
\[
i_*(\alpha+\beta+\gamma)=i^*\bigl(\betahat+\gammahat\bigr)+\partial\bigl(\gammahat\circ\Sigma\bigr)
\]
 for some \smash{$\betahat+\gammahat\in\cH^{2,2}$} (which~we take to be $\fh$-invariant by admissibility) and some section $\Sigma$. We will show that~$\smash{\fh=\fh_{(S',\betahat+\gammahat)}}$, and we note that this implies that $\bigl(S',\betahat+\gammahat\bigr)$ is a~Lie pair and~that \smash{$\fa=\fa_{(S',\betahat+\gammahat)}$}. The odd-odd bracket $\kappa+\gammahat-\gammahat\circ\Sigma\circ\kappa$ surjects onto $V\oplus\fh$. Therefore, if $A\in\fh$, there exists $\omega_A\in\Odot^2S'$ such that
	\[
		A = \kappa(\omega_A) + \gammahat(\omega_A) - \gammahat(\Sigma(\kappa(\omega_A))) \in V\oplus\fh,
	\]
	but since $A\in\fh$, we must have $\kappa(\omega_A)=0$, whence $\omega_A\in\fD$. Conversely, for any $\omega\in\fD$ we must have $\gammahat(\omega)\in\fh$, whence $\fh=\gammahat(\fD)=\fh_{(S',\betahat+\gammahat)}$.
\end{proof}

With the result above in mind, we say that a Lie pair $(S',\beta+\gamma)$ is \emph{admissible} if $\imath_v\beta+\gamma(\Sigma(v))$ preserves $S'$ and \emph{integrable} if it is admissible and the corresponding admissible cocycle is integrable.

\subsubsection{Equivalence of Lie pairs and correspondence}

Now fix $S'$ and a section $\Sigma\colon V\to\Odot^2S'$, and suppose that $\beta+\gamma\in\cK^{2,2}(\fa_-)\subseteq\cH^{2,2}$; that is, $\beta|_{V\otimes S'}=0$, which also implies that $\gamma|_{\Odot^2S'}=0$. Then $\fh_{(S',\beta+\gamma)}=0$, and $(S',\beta+\gamma)$ is trivially an integrable Lie pair; we will call it a \emph{trivial Lie pair}. For such a pair, the associated deformation obtained from Proposition~\ref{prop:lie-pair-admiss-cocycle} is trivial; $\fatilde_{(S',\beta+\gamma)}\cong\fa_{(S',\beta+\gamma)}$ which has underlying vector space~${V\oplus S'}$.

Note that if we have two admissible Lie pairs $(S',\beta+\gamma)$ and $(S',\beta'+\gamma')$ such that $(\beta'-\beta)+(\gamma'-\gamma)\in\cK^{2,2}(V\oplus S')$ (that is, $\beta'|_{V\otimes S'}=\beta|_{V\otimes S'})$ then $\fh_{(S',\beta'+\gamma')}=\fh_{(S',\beta+\gamma)}$, and $(\beta'-\beta)+(\gamma'-\gamma)$ is invariant under the action of this algebra. Both Lie pairs define the same admissible cocycle in Proposition~\ref{prop:lie-pair-admiss-cocycle}; it follows that if one of these pairs is integrable, so is the other, and $\fatilde_{(S',\beta'+\gamma')}=\fatilde_{(S',\beta+\gamma)}$. Conversely, two Lie pairs which have the same envelope $\fh$ and define the same admissible cocycle for the resulting graded subalgebra of $\fs$ define an $\fh$-invariant element in $\cK^{2,2}(V\oplus S')$.

Now recall the natural action of $G=\GrAut\fs$ on $\fs$. This action preserves the graded components of $\fs$ as $\fso(V)$-modules and the Dirac current $\kappa$, whence it follows that $G$ acts on the Spencer cochain spaces of $\fs$ (via induced actions on duals and tensor products of $G$-modules) and that this action preserves the Spencer differential $\partial$, so $G$ acts as automorphisms of the Spencer complex.

If $(S',\beta+\gamma)$ is a (resp.\ integrable, admissible) Lie pair and $g\in G$, then $(g\cdot S',g\cdot(\beta+\gamma))$ is also a (resp.\ integrable, admissible) Lie pair with
$\fh_{(g\cdot S',g\cdot(\beta+\gamma))}=g\cdot\fh_{(S',\beta+\gamma)}$ and $ \fa_{(g\cdot S',g\cdot(\beta+\gamma))}=g\cdot\fa_{(S',\beta+\gamma)}$.
Note that the map $g\cdot\Sigma:V\to\Odot^2(g\cdot S')$ defined by $(g\cdot\Sigma)(v)=g\cdot\Sigma\bigl(g^{-1}\cdot v\bigr)$ is a section of the Dirac current for $g\cdot S'$. Using $\Sigma$ and $g\cdot\Sigma$ to construct the deformations of Proposition~\ref{prop:lie-pair-admiss-cocycle}, in the integrable case we can construct an isomorphism of subdeformations of~$\fs$
\[
	\Phi\colon\ \fatilde_{(S',\beta+\gamma)} \to \fatilde_{(g\cdot S',g\cdot(\beta+\gamma))}
\]
as follows. Since we treat a filtered deformation $\fatilde$ as an alternative bracket structure on the same underlying space as its associated graded $\fa$, we have $\Phi=\Gr\Phi$ as linear maps. Thus we may define $\Phi\colon\fa_{(S',\beta+\gamma)}\to g\cdot\fa_{(S',\beta+\gamma)}$ on the underlying graded vector spaces by $\Phi(X)=g\cdot X$ for $X\in\fa_{(S',\beta+\gamma)}$. This preserves the grading, hence also the natural filtration (in the strong sense), and one can check that we have
\[
	\Phi(\comm{X}{Y}_{\beta+\gamma+\partial(\gamma\circ\Sigma)}) = \comm{\Phi(X)}{\Phi(Y)}_{g\cdot\beta+g\cdot\gamma+\partial((g\cdot\gamma)\circ(g\cdot\Sigma))},
\]
where the brackets are the deformed ones using the cocycle in the subscript. The inverse of $\Phi$ is of course given by the action of $g^{-1}$, which also preserves all structure, hence $\Phi$ is indeed an isomorphism of subdeformations.

The observations above suggest putting an equivalence relation on Lie pairs; in fact, we will define two different such notions.

\begin{Definition}\label{def:lie-pair-equiv}
	Two Lie pairs $(S',\beta+\gamma)$ and $(S'',\beta'+\gamma')$ are \emph{weakly equivalent} if $(S'',\beta'+\gamma')=(g\cdot S',g\cdot(\beta+\gamma)+\phi)$ for some $g\in G$ and $\phi\in\cK^{2,2}(V\oplus S'')^{\fh_{(S'',\beta'+\gamma')}}$, and they are \emph{strongly equivalent} if they are weakly equivalent and $\phi=0$.
\end{Definition}

These notions of equivalence are a generalisation of the equivalence relation on Lie pairs for the 11-dimensional Poincar\'e superalgebra found in \cite{Figueroa-OFarrill2017_1,Santi2022}. As previously noted, in those works~${G=\Spin(V)}$ rather than $\GrAut$ and the $\phi$ term does not appear in the equivalence relation since $\cK^{2,2}(V\oplus S')=0$ for all $S'$, so there is no distinction between strong and weak equivalence. We will see that weak equivalence is a more natural concept than strong equivalence, but strong equivalence will give us some additional information.

Fixing a Lie pair $(S',\beta+\gamma)$, we define
\[
	\fh^{\mathrm{max}}_{(S',\beta+\gamma)}:=\{A\in\fso(V)\mid A\cdot S'\subseteq S' \text{ and } A\cdot\beta=0\} \subseteq \fso(V);
\]
note that by definition of a Lie pair, this algebra contains the envelope $\fh_{(S',\beta+\gamma)}$, and
\[
	\fa^{\mathrm{max}}_{(S',\beta+\gamma)} = V\oplus S'\oplus\fh^{\mathrm{max}}_{(S',\beta+\gamma)}
\]
is a graded subalgebra of $\fs$ containing $\fa_{(S',\beta+\gamma)}$; in fact, it is the maximal such subalgebra whose odd part is $S'$ and whose zero-graded part preserves $\beta$.

If $(S',\beta+\gamma)$ is integrable, then the integrable cohomology class of $\fa_{(S',\beta+\gamma)}$ corresponding to the Lie pair pushes forward to an integrable cohomology class of \smash{$\fa^{\mathrm{max}}_{(S',\beta+\gamma)}$} which defines a~geometrically realisable filtered deformation \smash{$\fatilde^{\mathrm{max}}_{(S',\beta+\gamma)}$} which we call \emph{maximal} (with respect to~${(S',\beta+\gamma)}$), since it is maximal among realisable subdeformations into which $\fatilde_{(S',\beta+\gamma)}$ embeds as the odd-generated ideal.

\begin{Proposition}\label{prop:lie-pair-subdef-corresp}
	There is a one-to-one correspondence between weak equivalence classes of integrable Lie pairs $(S',\beta+\gamma)$ and isomorphism classes of highly supersymmetric odd-generated realisable subdeformations of $\fs$.
	
	If $(S',\beta+\gamma),(S'',\beta'+\gamma')$ are \emph{strongly} equivalent Lie pairs, then their maximal realisable subdeformations $\fatilde^{\mathrm{max}}_{(S'',\beta'+\gamma')},\fatilde^{\mathrm{max}}_{(S',\beta+\gamma)}$ are isomorphic as subdeformations.
\end{Proposition}

\begin{proof}
	The discussion above Definition~\ref{def:lie-pair-equiv} shows that weakly equivalent integrable Lie pairs define isomorphic odd-generated realisable subdeformations. By the converse part of Proposition~\ref{prop:lie-pair-admiss-cocycle}, every odd-generated realisable subdeformation is induced by a Lie pair. It remains to show that if the Lie pairs $(S',\beta+\gamma)$ and $(S'',\beta'+\gamma')$ are such that $\fatilde_{(S',\beta+\gamma)}\cong\fatilde_{(S'',\beta'+\gamma')}$ as subdeformations then they are weakly equivalent. If $\Phi\colon\fatilde_{(S',\beta+\gamma)}\to\fatilde_{(S'',\beta'+\gamma')}$ is such an isomorphism, $\Gr\Phi\colon\fa_{(S',\beta+\gamma)}\to\fa_{(S'',\beta'+\gamma')}$ is given by $x\mapsto g\cdot x$ for $g\in G$, so
	\[
	 	\fa_{(S'',\beta'+\gamma')} = g\cdot\fa_{(S',\beta+\gamma)}=\fa_{(g\cdot S',g\cdot(\beta+\gamma))}.
	\]
	But we have already seen that $\fatilde_{(S',\beta+\gamma)}\cong \fatilde_{(g\cdot S',g\cdot(\beta+\gamma))}$ for all $g\in G$, with the associated graded isomorphism of graded superalgebras given by the action of $g$, so we must also have $\fatilde_{(S'',\beta'+\gamma')}\cong \fatilde_{(g\cdot S',g\cdot(\beta+\gamma))}$ with the associated graded isomorphism being the identity map on $\fa':=g\cdot\fa_{(S',\beta+\gamma)}$. It follows from Proposition~\ref{prop:highly-susy-filtered-def} that the admissible cohomology classes in $\ssH^{2,2}(\fa'_-,\fa')$ associated to the two deformations are equal, whence their images under $i_*\colon\ssH^{2,2}(\fa'_-,\fa')\to\ssH^{2,2}(\fa'_-,\fs')$ are also equal, i.e., we have
	\[
		i^*[\beta'+\gamma'] = i^*[g\cdot\beta+g\cdot\gamma] \in \ssH^{2,2}(\fa'_-,\fs').
	\]
	Thus we find that $\phi:=(\beta'-g\cdot\beta)+(\gamma'-g\cdot\gamma)\in\cK^{2,2}(\fa'_-=V\oplus S'')$, and $\phi$ is invariant under the action of the envelope $\fh_{(S'',\beta'+\gamma')}$ since $\beta'+\gamma'$ and $g\cdot\beta+g\cdot\gamma$ both are.
	
	For the strong case, we note that $g\cdot\fh^{\mathrm{max}}_{(S',\beta+\gamma)}=\fh^{\mathrm{max}}_{(g\cdot S',g\cdot(\beta+\gamma))}$, whence it follows that strongly equivalent Lie pairs define isomorphic maximal subdeformations by an entirely analogous argument to the weak case.
\end{proof}

We note that there does not seem to be a converse statement to the second part of the result above; for this to hold, we would have to show that $\phi$ as in the proof vanishes if it is~$g\cdot\fh^{\mathrm{max}}_{(S',\beta+\gamma)}$-invariant, which there does not seem to be a reason to expect. Nonetheless, we have reduced the problem of classifying odd-generated realisable subdeformations of $\fs$ up to isomorphism to classifying integrable Lie pairs up to weak equivalence, and we can also obtain all maximal realisables with respect to a weak equivalence class from the strong equivalence classes it contains.

\subsubsection{Maximally supersymmetric case}

The case $S'=S$ brings a number of simplifications to the above discussion. First, for any $\beta+\gamma\in\cH^{2,2}$, $(S,\beta+\gamma)$ is a Lie pair if and only if the envelope $\fh_{(S,\beta+\gamma)}$ preserves $\beta$, and any Lie pair is admissible (recall our previous observation that if $\fa$ is maximally supersymmetric, all $\fa_0$-invariant $(2,2)$-cohomology classes are admissible), although one still must check that it is integrable. Since $\cK^{2,2}(V\otimes S)=0$, there is no distinction between strong and weak equivalence of Lie pairs; to any equivalence class of integrable Lie pairs there corresponds a unique (up to isomorphism) odd-generated realisable subdeformation and a unique maximal realisable subdeformation containing~it.

\section{Highly supersymmetric Lorentzian spin manifolds}
\label{sec:high-susy-spin-mfld}

In this section, we invert the perspective of Section~\ref{sec:killing-spinors-superalgebras}: rather than constructing a filtered subdeformation of a flat model $\fs$ from sections on a spin manifold (namely the Killing superalgebra), we seek to begin with such a filtered subdeformation and construct a spin manifold with admissible connection on which the elements of the algebra can be realised as sections, in doing so providing a partial converse to Theorem~\ref{thm:killing-algebra-filtered}. We restrict ourselves to the highly supersymmetric Lorentzian case, in which the homogeneity theorem (see Theorem~\ref{thm:homogeneity}) provides significant simplifications as it did in Section~\ref{sec:filtered-def-poincare} (although see Remark~\ref{rem:relaxing-homogeneity-assumption}), and our results will ultimately justify the approach of that section by demonstrating that our notion of geometric realisability (see Definition~\ref{def:geom-real}) is a good one. We will essentially generalise the results of \cite[Section~3.3]{Figueroa-OFarrill2017_1}. We make extensive use of the conventions and results related to homogeneous spaces discussed in Appendix~\ref{sec:homog-spaces}. We begin by establishing how we will deal with spinors on homogeneous manifolds.

\subsection{Homogeneous spin structures}
\label{sec:homogeneous-spin-mfld}

Spin structures on homogeneous manifolds were first treated systematically by B\"ar \cite{Bar1992} and by Cahen et al.\ \cite{Cahen1993}; our present treatment follows \cite{Figueroa-OFarrill2017_1}.

Let $(G,H,\eta)$ be a metric Klein pair (of any signature), let $(M=G/H,g)$ be the associated homogeneous space, which we assume to be oriented, and let $V=T_oM$ where $o=H$. Suppose that there exists a lift of the isotropy representation $H\to\SO(V)$ to $\Spin(V)$; that is, a Lie group morphism $H\to\Spin(V)$ making the following diagram commute
\[
\begin{tikzcd}
	&	& \Spin(V) \ar[d]	\\
	& H	\ar[r] \ar[ur,dashed]& \SO(V).
\end{tikzcd}
\]
Then the canonical map $\Spin(V)\to\SO(V)$ induces a $\Spin(V)$-equivariant bundle map
\[
	G\times_H\Spin(V)\to G\times_H\SO(V) \cong F_{\rm SO},
\]
which is a spin structure on $M$, where we recall that $F_{\rm SO}$ is the special orthonormal frame bundle. In order to define the isomorphism we have implicitly chosen an orthonormal frame over $o$ (as discussed in Appendix~\ref{sec:homog-spaces-principal-bundles}). We call this the \emph{homogeneous spin structure associated to the lift $H\to\Spin(V)$}. Note that $G$ (left) acts naturally on both of these bundles compatibly with the right action by the spin group, and that the bundle isomorphism is $G$-equivariant.

Using the isotropy representation and related representations of $H$, many natural bundles on $M$, particularly those associated to $F_{\rm SO}$ via representations of $\SO(V)$, can be viewed as associated bundles to the $H$-principal bundle $G\to M$. Homogeneous spin structures allow us to do the same with representations of the spin group; for example, if $S$ is a spinor module of~$\Spin(V)$, there is a $G$-equivariant isomorphism of vector bundles
\[
	\Sbundle := (G\times_H\Spin(V))\times_{\Spin(V)} S \cong G\times_H S.
\]

If $G$ is simply connected, the construction above classifies the spin structures on $(M,g)$ in definite signature.

\begin{Lemma}[\cite{Bar1992,Cahen1993}]\label{lemma:homogeneous-spin}
	If $(M,g)$ is a homogeneous Riemannian $G$-space where $G$ is simply connected, the equivalence classes of spin structures on $M$ are in one-to-one correspondence with spin lifts of the isotropy representation.
\end{Lemma}

The proof of this result relies on the ability to lift the action of $G$ on $F_{\rm SO}$ induced by the action on $M$ to an action on an arbitrary spin structure $P\to M$ which can then be localised at $o$ to find a lift of the isotropy representation. In definite signature, this follows from the fact that the spin group is connected and simply connected (for $\dim V \geq 3$), but it fails in indefinite signature due to the non-trivial topology of the spin group; in signature $(p\geq 1,q\geq 1)$, the spin group is disconnected, and its connected component is not simply connected for $(p\geq 2,q\geq 2)$. However, in the case that $M$ is a \emph{reductive} $G$-space, the above still holds in arbitrary signature~\cite{Alekseevsky2019,Cahen1993}.

To avoid these issues, one can simply hypothesise that the action of $G$ lifts to the spin structure; we then say that we have a \emph{homogeneous spin structure}. Such structures are in one-to-one correspondence with lifts of the isotropy representation.

\subsection{Reconstruction of highly supersymmetric backgrounds}
\label{sec:high-susy-spin-mfld-homogeneous}

Let us now turn directly to the question of constructing a background geometry on which a~filtered deformation of a Poincar\'e superalgebra can be realised. Let $(V,\eta)$ be a Lorentzian inner product space, $S$ a (possibly $N$-extended) spinor module of $\Spin(V)$, $\kappa\colon \Odot^2 S\to V$ a symmetric and causal Dirac current and $\fs$ the Poincar\'e superalgebra defined by this data (see Definition~\ref{def:flat-model-alg}). Our ultimate goal is to generalise \cite[Theorem~13]{Figueroa-OFarrill2017_1}, which says that geometrically realisable highly supersymmetric filtered subdeformations of the 11-dimensional Poincar\'e superalgebra are subalgebras of Killing superalgebras. We will find some potential obstructions to doing this in full generality, but the obstructions themselves will prove instructive.

\subsubsection{Super Harish-Chandra pairs}

We begin by introducing what can be viewed as a Lie supergroup associated to a Lie superalgebra. In a superspace formalism, the supersymmetric backgrounds we will construct could be considered as homogeneous spaces for such a supergroup. We will not fully use such a formalism, but the following definition will be useful for us. See \cite{Santi2010} and the references therein for more on the superspace perspective.

\begin{Definition}[super Harish-Chandra pair]\label{def:super-Harish-Chandra}
	A \emph{super Harish-Chandra pair} is a pair $(G_{\overline0},\fg)$ where $\fg$ is a finite-dimensional Lie superalgebra and $G_{\overline0}$ is a connected Lie group integrating $\fg_{\overline0}$ equipped with a representation $\Ad\colon G_{\overline0}\to\GL(\fg)$ (which we call the \emph{adjoint representation} of the pair) integrating the restriction $\ad\colon\fg_{\overline0}\to\fgl(\fg)$ of the adjoint representation of $\fg$ to $\fg_{\overline0}$.
\end{Definition}

Note that if $(G_{\overline0},\fg)$ is a super Harish-Chandra pair, the adjoint representation $\Ad\colon G_{\overline0}\to\GL(\fg)$ preserves the grading on $\fg$ since $\ad$ does. The induced subrepresentation on $\fg_{\overline0}$ is simply the adjoint representation of $G_{\overline0}$.

Let us now consider two motivating examples.

\begin{Example}
The spin cover $\Spin_0(V)\ltimes V$ of the connected Poincar\'e group $\SO_0(V)\ltimes V$ of course has the Poincar\'e algebra $\fs_{\overline 0}=\fiso(V)$ as its Lie algebra, and equipping it with the action $(a,v)\cdot s = a\cdot s$ on $\fs_{\overline 1}=S$ makes $(\Spin_0(V)\ltimes V,\fs)$ into a super Harish-Chandra pair.
\end{Example}

\begin{Example}
Let $(M,g,D)$ be a connected Lorentzian spin manifold with admissible connection and let $\fK_D$ be its Killing superalgebra (see Definition~\ref{def:killing-spinor}). Let us further assume that this is \emph{highly supersymmetric}, in the sense that $\dim\fS_D>\frac{1}{2}\dim S$ (i.e., $\fK_D$ is highly supersymmetric as a filtered subdeformation of $\fs$). Then, by Theorem~\ref{thm:homogeneity}, the values of the Killing vector fields arising as squares of Killing spinors span the tangent space at every point of $M$. Thus $M$ is \emph{locally homogeneous}; informally, it locally looks like a homogeneous space for a group with Lie algebra $\fV_D$. If we further assume that $M$ is geodesically complete, $\fV_D$ integrates to a connected subgroup $G$ of $\ISO(M,g)$ which acts transitively on $M$ \cite[Theorem~2.2]{Stehney1980}, whence $M$ is in fact homogeneous. By pulling back the action of $G$ to its universal covering group $\widetilde{G}$, it follows that $(M,g)$ is a homogeneous Lorentzian spin manifold for $\widetilde{G}$. Since $\widetilde{G}$ is simply connected, the restriction $\ad\colon\fV_D\to\fgl(\fK_D)$ of the adjoint representation of $\fK_D$ to $\fV_D$ integrates to a~representation $\Ad\colon\widetilde{G}\to\fK_D$. Thus $(\widetilde{G},\fK_D)$ is a super Harish-Chandra pair.
\end{Example}

The second example in particular shows us a way forward: where a Killing superalgebra is highly supersymmetric, the background which supports it is locally homogeneous, and if it is actually homogeneous, one can associate to it a super Harish-Chandra pair containing the Killing superalgebra.

Consider now a highly supersymmetric filtered subdeformation $\fg$ of the Poincar\'e superalgebra~$\fs$, say $\Gr\fg\cong\fa\subseteq\fs$ as graded Lie superalgebras so that $V=\fa_{-2}$, $S':=\fa_{-1}$ and $\fh:=\fa_0$ as usual. In particular, the filtration on $\fg$ is of the form $\fg=\fg^{-2}\supseteq\fg^{-1}\supseteq\fg^0\supseteq\fg^1=0$ and as Lie algebras, $\fh=\Gr_0\fg=\fg^0/\fg^1\cong\fg^0$, so there is a canonical embedding of $\fh$ as the zeroth part of the filtration of $\fg$.\footnote{This can also be seen by considering an explicit presentation of $\fg$ as a filtered deformation of $\fa$.}
Note also that the adjoint action of $\fa^0$ preserves each $\fg^k$ so that $\fg$ is also filtered as an $\fg^0$-module, whence $\Gr\fg$ is also naturally a graded $\fg^0$-module, and this action is the same one induced by the isomorphism $\fh\cong\fg^0$. We will thus identify $\fg^0$ and $\fh$, in particular treating $\fh$ as a subalgebra of $\fg$.

\begin{Lemma}\label{lemma:spin-lift-Ad}
	Let $(G_{\overline0},\fg)$ be a super Harish-Chandra pair where $\fg$ is a highly supersymmetric filtered subdeformation of $\fs$, and suppose that the connected subgroup $H$ of $G_{\overline0}$ generated by~${\fh=\fg^0}$ is closed. Then the inclusion $\fh\hookrightarrow\fso(V)$ integrates to a morphism $H\to\Spin(V)$. In particular, this lifts the isotropy representation and thus defines a homogeneous spin structure for $(G_{\overline0},H,\eta)$.
\end{Lemma}

\begin{proof}
	First note that $V\cong \fg^{-2}/\fg^{-1} \cong \fg_{\overline0}/\fh$ and $S'\cong\fg^{-1}/\fg^0\cong\fg_{\overline1}$ as $\fh$-modules. In particular, since $\eta$ is $\fh$-invariant, $(\fg_{\overline0},\fh,\eta)$ is a metric Lie pair, so $(G_{\overline0},H,\eta)$ is a metric Klein pair, and the inclusion $\fh\hookrightarrow\fso(V)$ is the isotropy representation of $\fh$. Similarly, at the level of groups, the adjoint representation induces a Lie group morphism $\varphi\colon H\to\SO(V)$ integrating $\fh\hookrightarrow\fso(V)$. This lifts to a morphism $\tilde\varphi\colon\Htilde \to \Spin(V)$ where $\pi_H\colon\Htilde\to H$ is the universal cover of $H$. We will show that there exists a diagonal arrow making the diagram
		\[
		\begin{tikzcd}
			& \Htilde \ar[d,"\pi_H"]\ar[r,"\tilde\varphi"]	& \Spin(V) \ar[d,"\pi"]	\\
			& H \ar[r,"\varphi"]\ar[ur,dotted]	& \SO(V)
		\end{tikzcd}
		\]
	commute, which proves the result. Recalling that there is a canonical short exact sequence of Lie groups
$ 1 \longrightarrow \pi_1(H) \longrightarrow\Htilde \longrightarrow H \longrightarrow1$,
	we will show that $\pi_1(H)\subseteq\ker\tilde\varphi$, hence $\tilde\varphi$ factors through $H$, giving the required arrow in the diagram.
	
	Note that the actions of $H$ on $V$ and $S'$ pull back to actions of $\Htilde$ (integrating the actions of~$\fh$) in which $\pi_1(H)$ acts trivially. On the other hand, $\Htilde$ also acts on $V$ and $S$ via the action of~$\Spin(V)$ on those spaces. The two actions of $\Htilde$ on $V$ agree by the commutation of the square diagram above. Since $\tilde\varphi$ integrates $\fh\hookrightarrow \fso(V)$ and the action of $\fh$ (via $\fso(V)$) on $S$ preserves $S'$, so too does the action of $\Htilde$ via $\Spin(V)$. The two actions of $\Htilde$ on $S'$ both integrate the same action of $\fh$ on $S'$, so they agree.
		
	The trivial action of $\pi_1(H)$ on $V$ gives us
$\pi_1(H)\subseteq\ker(\varphi\circ\pi_H)=\ker(\pi\circ\tilde\varphi)
$
	thus $\tilde\varphi(\pi_1(H))\subseteq\ker\pi=\{\pm\1\}$. On the other hand, the trivial action of $\pi_1(H)$ on $S'$ shows that $-\1\notin\tilde\varphi(\pi_1(H))$, thus $\pi_1(H)\subseteq\ker\tilde\varphi$. This shows that $\tilde\varphi$ factors through $H$ as claimed.
\end{proof}

\subsubsection{Reconstruction theorem}

We can now state and prove our final result which is a partial converse to Theorem~\ref{thm:killing-algebra-filtered} and generalisation of \cite[Theorem~13]{Figueroa-OFarrill2017_1}, with some caveats. Since this result essentially allows one to (locally) reconstruct a highly supersymmetric background from its Killing superalgebra, we follow \cite{Figueroa-OFarrill2017_1,Santi2022} and refer to it as the \emph{Reconstruction Theorem}.

\begin{Theorem}[reconstruction of highly supersymmetric background]\label{thm:homog-spin-mfld}	
	Let $\fg$ be a geometrically realisable highly supersymmetric filtered subdeformation of the Poincar\'e superalgebra $\fs$.
	Let $G_{\overline0}$ be the connected and simply connected Lie group corresponding to the Lie algebra $\fg_{\overline0}$ and suppose that the connected subgroup $H$ corresponding to $\fh$ is closed.
	Then there exists a homogeneous spin structure on the Lorentzian homogeneous space $M=G_{\overline0}/H$ with a Dirac current $\kappa$ and a~connection $D$ on the spinor bundle $\Sbundle=G_{\overline 0}\times_H S$ as well as an injective $\ZZ_2$-graded linear map~${\Psi\colon\fg\hookrightarrow\fV_D\oplus\fS_D}$ which restricts to a Lie algebra embedding $\fg_{\overline 0}\hookrightarrow\fV_D$ and which satisfies
	\begin{align}\label{eq:reconstruct-brackets}
		\Psi(\comm{X}{s}) = \eL_{\Psi(X)}\Psi(s),
		\qquad \Psi\qty(\comm{s}{s'}) = \kappa\qty(\Psi(s),\Psi(s'))
	\end{align}
	for all $X\in\fg_{\overline 0}$ and $s,s'\in \fg_{\overline 1}$.
	
	Furthermore, if $\Psi(\fg_{\overline 1})=\fS_D$ or if $D$ is admissible in the sense of Definition~{\rm\ref{def:killing-spinor}} otherwise, we find in particular that $\fg$ embeds into the Killing superalgebra $\fK_D$ of $(M,g,D)$.
\end{Theorem}

\begin{Remark}
	Let us make some comments about this result before giving the proof. First, it is necessary to hypothesise that $H$ is closed; it may be possible to adjust the choice of $\fh$ such that this is the case but we have not been able to show this. We note that this is also implicitly assumed in \cite[Theorem~13]{Figueroa-OFarrill2017_1}; although it is not in the statement of the theorem, it is stated in the preamble to the theorem.
	
	Moreover, the present theorem is very close to saying that \emph{all} geometrically realisable highly supersymmetric subdeformations of $\fs$ embed into Killing superalgebras, but there is an obstruction in that it does not seem to be guaranteed that the connection $D$ on $\Sbundle$ which we will construct in the proof is admissible. The algebraic conditions~(1) and~(2) of Definition~\ref{def:killing-spinor} are satisfied because we work with admissible cocycles which are expressed in terms of normalised cocycles of the whole Poincar\'e superalgebra $\fs$; the difficulty is in showing condition~(3), namely that $\eL_{\kappa_\epsilon}\beta=0$ for all $\epsilon\in\fS_D$; we can only show this for $\epsilon\in\Psi(\fg_{\overline 1})$. On the other hand, by a~dimension count we find that \emph{maximally} supersymmetric subdeformations do in fact embed in Killing superalgebras, since in this case we must have $\Psi(\fg_{\overline 1})=\fS_D$.
	
	To put a positive spin on this issue, it provides us with new motivation to study Killing superalgebras constructed from \emph{constrained} Killing spinors, which we briefly discussed in Section~\ref{sec:alg-kse}; while $D$ might not be guaranteed to satisfy the conditions for admissibility on sections of $\Sbundle$, it \emph{does} satisfy these conditions on sections of the sub-bundle $\Sbundle':=G_{\overline0}\times_H S'$ to which elements of~$\Psi(\fg_{\overline 1})$ belong and which we could describe as the kernel of some fibrewise linear operator on~$\Sbundle$.
\end{Remark}

For the proof, we will make use of the technology for working with connections on homogeneous spaces presented in Appendix~\ref{sec:nomizu}, namely Nomizu maps and Wang's theorem (see Theorem~\ref{thm:Wang}). If we have a spin-lift $\tilde\varphi$ of the isotropy representation, then~${P=G_{\overline0}\times_H\Spin(V)}$ is a $G_{\overline0}$-invariant principal $\Spin(V)$-bundle and $\tilde\varphi$ is a map from the isotropy group to the structure group of $P$ -- to connect with the notation in Theorem~\ref{thm:Wang}, we have $G=G_{\overline 0}$, $K=\Spin(V)$ and $\psi=\tilde\varphi$. Then one can easily check that the Nomizu map~${L\colon \fg\to\fso(V)}$ of the Levi-Civita connection (considered as a $G_{\overline0}$-invariant principal connection on $F_{\rm SO}$) satisfies the conditions of Wang's theorem so also corresponds with a principal connection on $P$, which is of course the spin-lift of the Levi-Civita connection.

\begin{proof}
	By hypothesis, we have a metric Klein pair $(G_{\overline0},H,\eta)$. Then by Lemma~\ref{lemma:spin-lift-Ad}, there exists a spin lift of the isotropy representation of $(G_{\overline0},H,\eta)$ with corresponding homogeneous spin structure $P=G_{\overline0}\times_H\Spin(V)$ and spinor bundle $\Sbundle=G_{\overline0}\times_H S$ on $M=G_{\overline0}/H$. We identify $V$ with $T_oM\cong\fg/\fh$. Since $H$ is connected, the image of the spin-lift of the isotropy representation lies in $\Spin_0(V)$, whence $(M,g)$ is time-orientable and so we can construct a bundle Dirac current~${\kappa\colon \Odot^2\Sbundle\to TM}$ by Lemma~\ref{lemma:bundle-dirac-current-exist}.
	
Now let $[\mu_-]$ be the admissible Spencer cohomology class for $\fa=\Gr\fg\subseteq \fs$ corresponding to~$\fg$, so that there exists \smash{$\beta+\gamma\in\bigl(\cH^{2,2}\bigr)^\fh$} satisfying $i_*[\mu_-]=i^*[\beta+\gamma]\in\ssH^{2,2}(\fa_-;\fs)$. Recall that~${\beta\in\Hom(V\otimes S,S)\cong V^*\otimes\End S}$ is $\fh$-invariant and thus $H$-invariant, so by Frobenius reciprocity (see Theorem~\ref{thm:Frobenius-recipr}) it defines a unique $G_{\overline0}$-invariant 1-form with values in $\End\Sbundle$, although we will actually adopt a sign convention where we work with $\boldsymbol{\beta}\in\Omega^1(M;\End\Sbundle)^{G_{\overline0}}$ corresponding to~$-\beta$. We define the connection $D:=\nabla-\boldsymbol{\beta}$.

Since $G_{\overline 0}$ acts by isometries by construction and $\boldsymbol{\beta}$ is $G_{\overline 0}$-equivariant, for each $X\in\fg_{\overline 0}$ the fundamental vector field $\xi^M_X\in\fX(M)$ is Killing and preserves $\boldsymbol{\beta}$, so it lies in $\fV_D$. Recalling that the assignment $X\to\xi^M_X$ is a Lie algebra \emph{anti}-homomorphism $\fg_{\overline 0}\to \fX(M)$, we define $\Psi(X)=-\xi^M_X$. This gives us a Lie algebra morphism $\Psi\colon\fg_{\overline 0}\to\fV_D$.
	
	To construct the spinor fields associated to elements of $\fg_{\overline 1}=S'$ and show that they are $D$-parallel, we will exploit the homogeneous bundle structure of $S$. Let us define a left action of $H$ on $G_{\overline 0}$ by $h\cdot g=gh^{-1}$ where we use group multiplication in $G_{\overline 0}$ on the right-hand side. Then, using some standard associated bundle theory, we have identifications between spaces of sections and $H$-equivariant maps
	\[
		\fX(M) \cong C^\infty\qty(G_{\overline 0};V)^H,\qquad
		\fS \cong C^\infty\qty(G_{\overline 0};S)^H.
	\]
	For $X\in\fg_{\overline 0}$, the fundamental vector field $\xi^M_X$ corresponds to the map $G_{\overline 0}\to V\cong\fg_{\overline 0}/\fh$ given by $X\mapsto \overline{\Ad_{g^{-1}}X}$ (where the overline denotes the image under the quotient by $\fh$). Similarly, for $s\in \fg_{\overline 1}$ we define $\Psi(s)\in\fS$ as the spinor field associated to the map $G_{\overline 0}\to S$ defined by~${g\mapsto \Ad_{g^{-1}}s}$. We now use Wang's theorem and the Levi-Civita Nomizu map to show that these spinor fields are $D$-parallel. First, let us take a presentation of the form \eqref{eq:admiss-cocycle-brackets} for $\fg$; in particular, we have a (non-reductive) split $\fg_{\overline0}=\fh\oplus V$ and a map $\lambda\colon V\to\fso(V)$. Then, writing elements $X\in\fg_{\overline 0}$ in the form $X=A+v$ for $A\in\fh$, $v\in V$, we can follow the construction of the Nomizu map $L\colon \fg_{\overline 0}\to\fso(V)$ corresponding to the Levi-Civita connection (see Appendix~\ref{sec:nomizu}) to compute
$L(A+v) = A + \lambda(v)$.
	In particular, for $s\in S'$ we have $L(A+v)\cdot s = \comm{A+v}{s} - \beta(v,s)$, where we use the definition of the bracket on $\fg$. Since this is also the map corresponding to the spin-lift of the Levi-Civita connection $\nabla$ under Wang's theorem, we can compute the $H$-equivariant map $G_{\overline 0}\to S$ corresponding to the spinor field $\nabla_{\xi^M_X}\Psi(s)$ as
	\[
		g \mapsto (L(\Ad_{g^{-1}}X) - \ad_{\Ad_{g^{-1}}X})\qty(\Ad_{g^{-1}}s) = - \beta\qty(\overline{\Ad_{g^{-1}}X},\Ad_{g^{-1}}s),
	\]
	where the equality is due to the expression for $L$ above. But this is nothing but the map corresponding to $\boldsymbol{\beta}\bigl(\xi^M_X,\Psi(s)\bigr)$ (recalling our sign convention for $\boldsymbol{\beta}$). By homogeneity, the values of the fundamental vector fields span every tangent space, so we have shown that $\nabla\Psi(s)=\boldsymbol{\beta}\Psi(s)$, whence $\Psi(s)\in\fS_D$ as required. We now have a $\ZZ_2$-graded map $\Psi\colon\fg\to\fV_D\oplus\fS_D$ which can show satisfies the equations~\eqref{eq:reconstruct-brackets} by examining the Killing transport data at the basepoint (see Section~\ref{sec:alg-structure-killing}, in particular the explicit presentation of the brackets of the Killing superalgebra as a filtered deformation \eqref{eq:KSA-deformed-brackets}). This will also show that $\Psi$ is injective since no non-zero element of $\fg$ has trivial transport data.
	
	Finally, we note that if $D$ is admissible, $\Psi\colon\fg\hookrightarrow\fK_D$ is a Lie superalgebra embedding. The Spencer cocycle conditions for $\beta+\gamma$ ensure $\boldsymbol{\beta}$ satisfies conditions~(1) and~(2) of Definition~\ref{def:killing-spinor}. Condition~(3) ($\eL_{\kappa_\epsilon}\boldsymbol{\beta}=0$ for all $\epsilon\in\fS_D$) is not automatically satisfied, but observe that we have~$\eL_{\xi^M_X}\boldsymbol{\beta}=0$ for all $X\in\fg_{\overline 0}$ by $G_{\overline 0}$-invariance of $\boldsymbol{\beta}$; in particular, $\eL_{\kappa_\epsilon}\boldsymbol{\beta}=0$ holds for~${\epsilon\in\Psi(\fg_{\overline 1})}$.
\end{proof}

We note that we have glossed over some details of the calculations above, in particular when constructing the spinorial part of the map $\Psi$; these details are routine exercises in principal connections but somewhat cumbersome to write out in full. A more elegant way of dealing with them is through the superspace formalism introduced in \cite{Santi2010} which is used, somewhat implicitly, in the proof of \cite[Theorem~13]{Figueroa-OFarrill2017_1}, but which is beyond the scope of this work.

\begin{Remark}
We have already noted some differences between the result above and the much cleaner statement in the 11-dimensional case. Another difference to be aware of is that a geometry $(M,g,D)$ does not seem to be determined by the odd-generated ideal $\comm{\fg_{\overline 1}}{\fg_{\overline 1}}\oplus\fg_{\overline 1}$ (or weak equivalence class of a maximal Lie pair $(S'=\Sbundle_p,\beta)$) in general since we may have $\cK^{2,2}(\fa_-)\neq 0$; the discussion of Section~\ref{sec:classification-odd-gen} shows that if this is the case then it is possible for two realisable subdeformations of $\fs$ to have the same odd-generated ideal but for neither to embed in the other. On the geometric side, this means that two simply connected backgrounds may be isometric to one another but have different admissible connections $D$ (corresponding to different flux configurations, in supergravity terminology). On the other hand, it is not hard to see that a \emph{strong} equivalence class of a maximal Lie pairs \emph{does} determine the connection uniquely.
\end{Remark}

\appendix

\section{Further background}
\label{app:further-background}

\subsection{Inner product spaces and spinors}
\label{sec:inner-products-spinors}

For the sake of brevity, we will not give a detailed account of the theory of Clifford algebras, spin groups and spinors here; standard references include \cite{Harvey1990,Lawson1989}; see also \cite{JMF-spin}. We will briefly establish some conventions, however. Throughout, let $(V,\eta)$ be a (pseudo-)inner product space.
We recall that using $\eta$, the vector spaces $V$ and $V^*$ can be identified via the \emph{musical isomorphisms}~${\flat\colon V\to V^*}$, $v\mapsto v^\flat$ where
$v^\flat(w) = \eta(v,w)
$
for all $v,w\in V$ and $\sharp\colon V^*\to V$, $\alpha\mapsto\alpha^\sharp$ which is defined as the inverse of $\flat$.

\subsubsection{Orthogonal groups and algebras}
\label{sec:orth-groups-algs}

We define the \emph{orthogonal group} of $(V,\eta)$ as the group of linear automorphisms of $V$ which preserve $\eta$,
\[
	\Orth(V,\eta) := \{A\in\GL(V) \mid \forall v,w\in V,\, \eta(Av,Aw) = \eta(v,w)\},
\]
and the \emph{special orthogonal group} as the subgroup of elements of $\Orth(V,\eta)$ which preserve orientations,
$\SO(V,\eta) := \{A\in \Orth(V) \mid \det A = 1\}$,
though we typically suppress the inner product $\eta$ in the notation where there is no ambiguity. The latter group is a normal subgroup of the former of index~2. Some basic topological properties depend on the signature of $\eta$ as follows.
\begin{itemize}\itemsep=0pt
	\item If $\eta$ is (positive- or negative-)definite, $\Orth(V)$ has two connected components, with the connected component of the identity being $\SO(V)$, and both groups are compact.
	\item If $\eta$ is not definite, $\Orth(V)$ has four connected components, two of which constitute $\SO(V)$. We denote the connected component of the identity by $\SO_0(V)$.\footnote{Following naming conventions used for the Lorentzian case in physics, this group is sometimes called the \emph{proper orthochronous orthogonal group}. It consists of invertible linear isometries of $(V,\eta)$ which preserve both orientations and time orientations on $V$.} None of these groups are compact.
\end{itemize}
Clearly, all of these groups share the same Lie algebra, which can be identified with the Lie algebra of $\eta$-skew-symmetric endomorphisms
\[
	\fso(V) := \{A\in\End(V) \mid \eta(Av,w) + \eta(v,Aw) = 0\}
\]
under the commutator bracket. We also define $\Orth(p,q):=\Orth(\RR^{p,q})$ etc., where $\RR^{p,q}$ is the standard inner product space with signature $(p,q)$; that is, with $p$ positive eigenvalues and $q$ negative eigenvalues.

\subsubsection{Clifford algebras}
\label{sec:clifford-algs}

The defining relation for the Clifford algebra $\Cl(V,\eta)$ is
\begin{equation}\label{eq:clifford-rel}
	v\cdot v = \pm\eta(v,v)\1
\end{equation}
for all $v\in V$. Note that there is a choice of sign convention here; the $+$ is more common in the physics literature while the $-$ is usually preferred in mathematics. We adopt the mathematical convention with the $-$ sign in defining $\Cl(p,q):=\Cl(\RR^{p,q})$. Any inner product space is (non-canonically) isomorphic to $\RR^{p,q}$ for some unique pair $(p,q)$, and such an isomorphism induces an isomorphism of associative algebras $\Cl(V,\eta)\cong \Cl(p,q)$. In previous works on supersymmetry and Spencer cohomology \cite{Beckett2021,deMedeiros2016,deMedeiros2018,Figueroa-OFarrill2017} and in the upcoming work treating 2-dimensional examples mentioned in the introduction, a $+$ sign is adopted for calculational convenience. We also typically suppress $\eta$ in the notation for the Clifford algebra (and the related groups defined below), writing $\Cl(V)=\Cl(V,\eta)$.

There is an isomorphism of vector spaces $\Wedge^\bullet V\simeq \Cl(V,\eta)$ which is most easily described in terms of an orthonormal basis $\{\be_\mu\}_{\mu=1}^{\dim V}$, although it is actually independent of this choice; to each exterior product of the form $\be_{\mu_1}\wedge\be_{\mu_2}\wedge\cdots\wedge\be_{\mu_r}$, we associate the Clifford product~${\be_{\mu_1}\cdot\be_{\mu_2}\cdot\cdots\cdot\be_{\mu_r}}$. We thus identify these two spaces and view Clifford multiplication $\cdot$ as a~new algebra structure on the exterior algebra $\Wedge^\bullet V$.

The exterior algebra $\Wedge^\bullet V$ has a natural $\ZZ$-grading by form degree. Clifford multiplication does not respect this grading, but it does respect the $\ZZ_2$-grading given by taking the $\ZZ$-grading modulo 2. We denote the even subalgebra of the Clifford algebra \big(which is generated by the subspace corresponding to $\Wedge^2 V$\big) by $\Cl_{\overline0}$.

There is a vector space embedding $\omega\colon\fso(V)\hookrightarrow\Wedge^\bullet V$ whose image is $\Wedge^2 V$ given by $A\mapsto \pm\frac{1}{2}\omega_A$ where $\omega_A$ is the 2-vector defined by the equation
$\imath_{v^\flat}\omega_A = -Av
$
where $\imath_{v^\flat}$ is contraction with the musical image $v^\flat\in V^*$ of an arbitrary vector $v\in V$.
This an embedding of Lie algebras if $\Wedge^\bullet V$ is equipped with the commutator Lie bracket with respect to Clifford multiplication.

\subsubsection{Pinor and spinor representations}
\label{sec:pinor-spinor}

The Clifford algebra $\Cl(V)$ has either a single isomorphism class of irreducible real finite-dimensional modules or two such classes, depending on $p-q\mod 4$ where $(p,q)$ is the signature of $(V,\eta)$.
\begin{itemize}\itemsep=0pt
	\item For $p-q\neq 3\mod 4$, there is a unique module which we denote by $\ssS$.
	\item For $p-q = 3 \mod 4$, there are two such modules which can be distinguished by the fact that a unit volume element in $\Wedge^{p+q}V$ acts as $+\1$ (the identity) on one of these modules and $-\1$ on the other; given an orientation on $V$, we denote by $\ssS$ the module on which the canonical unit volume element acts as $+\1$.
\end{itemize}
All of these are referred to as (irreducible, real) \emph{pinor} modules. Under the induced action of the even subalgebra $\Cl_{\overline 0}(V)$, or equivalently the spin group (see below), the pinor modules may be reducible, depending on $p-q\mod 8$.
\begin{itemize}\itemsep=0pt
	\item For $p-q=3,5,7 \mod 8$, $\ssS$ is irreducible; for later convenience we define $\ssS_1:=\ssS$ in this case.
	\item For $p-q=1,2,6$, $\ssS$ has a decomposition into irreducible $\Cl_{\overline0}(V)$-modules $\ssS=\ssS_+\oplus\ssS_-$ which happen to be isomorphic (as real representations); let us set $\ssS_1:=\ssS_+\cong\ssS_-$ here.
	\item For $p-q=0,4\mod 8$, there is a similar decomposition $\ssS=\ssS_+\oplus\ssS_-$, except here $\ssS_+\not\cong\ssS_-$. In any case, $\ssS_1$ or $\ssS_\pm$ are called the (irreducible, real) \emph{spinor} modules.
\end{itemize}

\subsubsection{Pin and spin groups}
\label{sec:pin-spin-group}

We define the \emph{pin group} and \emph{spin group} as the subgroups
\begin{gather*}
	\Pin(V) = \{v_1\cdot v_2\cdot\dots\cdot v_{r} \mid v_i\in V, \, \eta(v_i,v_i) = \pm 1 \},	\\
	\Spin(V) = \{v_1\cdot v_2\cdot\dots\cdot v_{2r} \mid v_i\in V, \, \eta(v_i,v_i) = \pm 1 \} = \Pin(V)\cap \Cl_{\overline0},
\end{gather*}
of the group of units $\Cl(V)^\times$. The spin group is a normal subgroup of the pin group of index 2, and similarly to the orthogonal groups, we have the following:
\begin{itemize}\itemsep=0pt
	\item If $\eta$ is definite, $\Pin(V)$ has two connected components, with the connected component of the identity being $\Spin(V)$, and both groups are compact.
	\item If $\eta$ is not definite, $\Pin(V)$ has four connected components, two of which constitute $\Spin(V)$. We denote the connected component of the identity by $\Spin_0(V)$. None of these groups are compact.
\end{itemize}
There is a canonical double-covering of Lie groups $\Pin(V)\to\Orth(V)$ which restricts to double-coverings $\Spin(V)\to \SO(V)$ and (where the groups exist) $\Spin_0(V)\to\SO_0(V)$. These covering maps differentiate to a canonical isomorphism of Lie algebras $\mathfrak{spin}_{\overline0}(V)=\mathfrak{spin}(V)\to\fso(V)$, thus we will identify these algebras.

The pin and spin groups act on the pinor modules by restriction of the action of the Clifford algebra. The irreducible pinor module $\ssS$ is irreducible under the action of $\Pin(V)$. It is irreducible under the action of $\Spin(V)$ or $\Spin_0(V)$ if and only if it is irreducible under $\Cl_{\overline 0}(V)$ and, in the appropriate signatures, $\ssS_1$ or $\ssS_\pm$ are irreducible under the action of these groups. The induced action of $\fso(V)$ on the (s)pinor modules is given by
\begin{equation}\label{eq:sov-spinor-action}
	A\cdot s = \pm\frac{1}{2}\omega_A\cdot s
\end{equation}
for $A\in\fso(V)$ and $s\in \ssS$, where the sign is the same as that in the Clifford relation \eqref{eq:clifford-rel} and the dot on the left-hand side is the Clifford action.

\subsection{Lie derivatives}
\label{sec:lie-der}

A standard concept in differential geometry, Lie derivatives play a central role in the definition of the Killing (super)algebra. In this appendix, we will first introduce some useful formulae for Lie derivatives of tensor fields and then show how these generalise to allow us to take Lie derivatives of spinor fields along Killing vectors.

\subsubsection{Lie derivative of vectors and tensors}
\label{sec:lei-der-vector-tensor}

Let $(M,g)$ be a connected pseudo-Riemannian manifold. We denote its Lie algebra of Killing vector fields\footnote{The notation should not be taken to imply that this is the Lie algebra of the isometry group of $(M,g)$ which is in general only a subalgebra of $\fiso(M,g)$; Killing vector fields may not integrate to a global action by a one-parameter group of isometries. We will not need to discuss the isometry group itself in this appendix or in the main text so this should not cause confusion.}
by $\fiso(M,g)$. The Levi-Civita connection will be denoted by $\nabla$. The bundle of skew-symmetric endomorphisms (with respect to $g$) of $TM$ will appear when we discuss Killing vectors. We denote the space of sections of this bundle by $\fso(M,g)$ and note that it is an (infinite-dimensional) Lie algebra under the commutator Lie bracket. The perspective presented here is due to Kostant \cite{Kostant1955}.

Given a vector field $X\in\fX(M)$, we define an endomorphism of $TM$ by
$A_X(Y) := -\nabla_Y X
$
for all $Y\in\fX(M)$. This gives us a useful formula for the Lie derivative along $X$: for $Y\in\fX(M)$,
\begin{equation}\label{eq:lie-der-vec-formula}
	\eL_X Y = \comm{X}{Y} = \nabla_X Y - \nabla_Y X = \nabla_X Y + A_X Y,
\end{equation}
where we have used the torsion-freeness of $\nabla$. We will also use
\begin{equation}\label{eq:comm-formula}
	\comm{X}{Y} = A_X Y - A_Y X.
\end{equation}
A similar formula holds a tensor field $T$ of any type
\begin{equation}\label{eq:lie-der-tensor-formula}
	\eL_X T = \nabla_X T + A_X \cdot T,
\end{equation}
where the action of $A_X$ on $T$ is defined via a Leibniz formula. Note that we have a different Leibniz rule of the form
\[
	\nabla_X (A_Y\cdot T) = \nabla_X(A_Y)\cdot T + A_Y\cdot\nabla_X T
\]
for all $X,Y\in\fX(M)$, or, abstracting $T$, $\comm{\nabla_X}{A_Y}=\nabla_X A_Y$ as endomorphisms of the space of tensor fields.

We can use the formula \eqref{eq:lie-der-tensor-formula} to show that $X$ is a Killing vector if and only if $A_X\in\fso(M,g)$: using the metric-compatibility of $\nabla$, we have $\eL_X g = A_X\cdot g$, so
\begin{equation}\label{eq:killing-condition}
	\eL_X g = 0 \iff g(A_X Y, Z) + g(Y,A_X Z) = 0,\qquad\forall Y,Z\in\fX(M).
\end{equation}
One can also show \cite[Lemma~2.2]{Kostant1955} that, for $X\in\fiso(M,g)$ and $Y\in\fX(M)$,
\begin{equation}\label{eq:nabla-A_X}
	\nabla_Y A_X = -R(Y,X),
\end{equation}
where $R$ is the Riemann curvature, so in particular $R(Y,X)\in\fso(M,g)$. Then for two Killing vectors $X,Y\in\fiso(M,g)$ and a vector field $Z\in\fX(M)$, using equations~\eqref{eq:comm-formula} and \eqref{eq:nabla-A_X} along with the algebraic Bianchi identity,
\begin{align*}
	A_{\comm{X}{Y}}Z
		&= -\nabla_Z\qty(A_X Y - A_Y X)= -(\nabla_ZA_X)Y - A_X\nabla_ZY + (\nabla_ZA_Y)X + A_Y\nabla_ZY\\
		&= R(Z,X)Y + A_XA_YZ - R(Z,Y)X- A_YA_XZ= \comm{A_X}{A_Y}Z - R(X,Y)Z.
\end{align*}
We thus have
\begin{equation}\label{eq:A-comm}
	A_{\comm{X}{Y}} = \comm{A_X}{A_Y} - R(X,Y)
\end{equation}
for all Killing vectors $X$ and $Y$.

\subsubsection{Spin structures and spinor bundles in arbitrary signature}
\label{sec:spin-structures}

Now let $(M,g)$ be a pseudo-Riemannian spin manifold of signature $(p,q)$. We recall that this means that $M$ is oriented and that there exists a $\Spin(p,q)$-principal bundle $P\to M$ over $M$ and a $\Spin(p,q)$-equivariant bundle map $\varpi\colon P\to F_{\rm SO}$, where $F_{\rm SO}\to M$ is the special orthonormal frame bundle of $(M,g)$ on which the spin group acts via the natural map $\Spin(p,q)\to\SO(p,q)$.

In indefinite signature, the structure group of the bundle $P\to M$ can be reduced to $\Spin_0(p,q)$ if and only if $(M,g)$ is time-orientable, and if this is the case then $\varpi$ restricts to a $\Spin_0(p,q)$-equivariant mapping $\varpi_0\colon P_0\to F_{SO_0}$ of the reduced bundles. Note that many authors assume time-orientability and refer to this reduced structure as the spin structure.\footnote{The reduced structure has been called a \emph{strong spin structure}, and manifolds admitting them \emph{strongly spin} \cite{Cortes2021,Shahbazi2024_2,Shahbazi2024_1}.}
	
Let $S$ be a (possibly $N$-extended) spinor module of $\Spin(p,q)$ in the sense discussed in Section~\ref{sec:isom-alg-spinors} and let $\Sbundle:=P\times_{\Spin(p,q)} S\to M$ be the associated vector bundle over $M$, which we will call the \emph{spinor bundle}. Sections of $\Sbundle$ will be called \emph{spinor fields}, and we denote the space of such fields by $\fS=\Gamma(\Sbundle)$. We also let $\Cl(M,g)=P\times_{\Spin(p,q)}\Cl(p,q)\cong F_{\rm SO}\times_{\SO(p,q)}\Cl(p,q)$ be the associated Clifford bundle (where the action is via the twisted adjoint, or the algebra automorphism induced by the natural action on $\RR^{p,q}$) and recall that we can identify it with the exterior bundle $\Wedge^\bullet T^*M$. Sections of this bundle can be Clifford-multiplied and Clifford-act on spinor fields via pointwise Clifford multiplication and action. The (spin-lift of) the Levi-Civita connection extends to all of these bundles, and there is a Leibniz rule for $\nabla$ with respect to Clifford multiplication and action.

The adjoint bundle $\ad F_{\rm SO}\cong\ad P$ can be identified with the space of $g$-skew-symmetric endomorphisms of $TM$, and it embeds in $\Cl(M,g)$ as $\Wedge^2 TM$. As a Lie algebra, the space of sections $\fso(M,g)$ acts on spinor fields via the Clifford action, i.e., we apply equation \eqref{eq:sov-spinor-action} pointwise.

\subsubsection{Lie derivative of spinors}
\label{sec:lie-der-spinors}

The formulae~\eqref{eq:lie-der-vec-formula} and \eqref{eq:lie-der-tensor-formula} suggest the following definition of the Lie derivative to spinor fields.

\begin{Definition}[spinorial Lie derivative \cite{Kosmann1971}]\label{def:lie-der}
	The Lie derivative of a spinor field $\epsilon\in\fS=\Gamma(\Sbundle)$ along a Killing vector $X\in\fiso(M,g)$ is given by
$
		\eL_X \epsilon := \nabla_X \epsilon + A_X\cdot \epsilon$.
\end{Definition}

Note, however, that while the formulae for the Lie derivatives of vector and tensor fields hold for the Lie derivative along any vector field $X$, the right-hand side of the spinor formula is only defined for Killing vectors\footnote{The spinorial Lie derivative is actually defined for \emph{conformal Killing} vector fields in \cite{Kosmann1971}, but we will only use it for Killing vectors here.}
$X$, since sections of $\End(TM)$ do not act on $\fS$, while~$\fso(M,g)$ does. Note that $\fso(M,g)$ is an (infinite-dimensional) Lie algebra and the action on $\fS$ defines a~representation. This action obeys a Leibniz rule of the form $\comm{\nabla_X}{A}=\nabla_X A$ as endomorphisms of $\fS$ for all $A\in\fso(M,g)$ and $X\in\fX(M)$; more explicitly,
$\nabla_X (A\cdot\epsilon) = (\nabla_X A)\cdot\epsilon + A\cdot(\nabla_X \epsilon)
$
for all $\epsilon\in\fS$. It follows immediately from the above that we have yet another Leibniz rule~${\eL_X (A\cdot\epsilon) = (\eL_X A)\cdot\epsilon + A\cdot(\eL_X \epsilon)
}$
for all $X\in\fiso(M,g)$, $A\in\fso(M,g)$, or equivalently~${\comm{\eL_X}{A}=\eL_X A}$ as endomorphisms of $\fS$. One can also show that where $X$ is a Killing vector field and $Y$ any vector field, the following compatibility condition between the spinorial Lie derivative and the Levi-Civita connection on spinors holds
\begin{equation}\label{eq:lie-der-compat}
	\comm{\eL_X}{\nabla_Y} = \nabla_{\comm{X}{Y}},
\end{equation}
while if both $X$ and $Y$ are Killing vector fields,
$\eL_{\comm{X}{Y}} = \eL_X\eL_Y - \eL_Y\eL_X$.
In particular, the spinorial Lie derivative gives $\fS$ the structure of an $\fiso(M,g)$ module.

We can also use the Leibniz rule to define an action of $\fso(M,g)$, and thus a Lie derivative along Killing vectors, on sections of $\End\Sbundle$ and on spaces of forms with values in this bundle~$\Omega^p(M;\End\Sbundle)$. More generally, the space of sections of any bundle constructed by taking duals and tensor products of $TM$ and $\Sbundle$ is equipped with such a structure, and all the Leibniz rules and compatibility conditions discussed above also hold when acting on these spaces.


\subsection{Grading and filtration of superspaces and superalgebras}
\label{sec:grading-filtration-superspace}

We assume that the reader is familiar with ordinary $\ZZ$-graded vector spaces and with vector superspaces but recall and set some terminology for $\ZZ$-graded superspaces and algebras, and we also introduce filtrations.

\subsubsection{Graded vector superspaces and superalgebras}
\label{sec:graded-superalg}

A vector superspace $V$ is said to be \emph{$\ZZ$-graded} with grading $V=\bigoplus_{k\in\ZZ}V_k$ if the latter is a~$\ZZ$-grading of $V$ as a vector space and the $\ZZ$-grading is compatible with the $\ZZ_2$ grading on $V$ in the sense that an element of degree $d$ in the $\ZZ$-grading has degree (parity) $\overline{d}=d\mod 2$ in the $\ZZ_2$ grading. Note that we can always make a $\ZZ$-graded \emph{vector space} into a $\ZZ$-graded vector \emph{superspace} by taking the $\ZZ$-grading modulo $2$. We say that a $\ZZ$-graded vector (super)space~${V=\bigoplus_{k\in\ZZ}V_k}$ has \emph{finite depth} $h\in\ZZ_{\geq 0}$ if $V_{-h}\neq 0$ but $V_k=0$ for all $k<-h$.

A Lie superalgebra $\fg$ is said to be \emph{$\ZZ$-graded} with grading $\fg=\bigoplus_{k\in\ZZ}\fg_k$ if the latter is a~$\ZZ$-grading of $\fg$ as a vector superspace and the superbracket respects this grading in the sense that~${\comm{\fg_k}{\fg_l}\subseteq\fg_{k+l}}$.

If $V$ is a vector superspace, we define the exterior algebra $\Wedge^\bullet V:=\mathcal{T}(V)/I$ where $\mathcal{T}(V)$ is the tensor algebra and I is the ideal generated by elements of the form $v\otimes w + (-1)^{\abs{v}\abs{w}}w\otimes v$, where~$v$,~$w$ are homogeneous elements with parities $\abs{v}$, $\abs{w}$, respectively. This is a $\ZZ$-graded associative algebra of depth 0: the grading is $\Wedge^\bullet V = \bigoplus_{k=0}^\infty \Wedge^k V$ where $\Wedge^k V$ is the image of~${\mathcal{T}^k(V)= V^{\Otimes k}}$ under the quotient map. One can show that, for $V=V_{\overline{0}}\oplus V_{\overline{1}}$,
\begin{equation}\label{eq:super-wedge}
	\Wedge^k V = \bigoplus_{i+j=k}\Wedge^i V_{\overline{0}}\otimes\Odot^j V_{\overline{1}}.
\end{equation}
It will occasionally be useful to remind the reader of this definition of $\Wedge^k$ of a vector superspace; when we do so, we will say that we mean it in the ``super-sense''. There is of course an analogous ``super'' version of the symmetric algebra $\Odot^\bullet V$ but we will not require it in this work.

\subsubsection{Filtered vector spaces and algebras}
\label{sec:filtered-alg}

We say that a vector space $V$ is \emph{$\ZZ$-filtered} if there exists a decreasing sequence of vector subspaces~${V \supseteq \dots \supseteq V^{k-1} \supseteq V^k \supseteq V^{k+1} \supseteq \dots \supseteq 0}$.
To emphasise that we are taking $V$ with its filtration, we will denote it by~$V^\bullet$. We say that~$V^\bullet$ is \emph{finite} if there exist some $k_1$, $k_2$ such that~${V^{k_1}=V}$ and $V^{k_2}=0$. Of course, if $V$ is finite-dimensional then $V^\bullet$ must be finite. We say that $V^\bullet$ (not necessarily finite) has \emph{finite depth}~${h\in\ZZ_{\geq 0}}$ if $V_{-h} = V$ but $V_k\neq V$ for all $k>-h$.

Given a $\ZZ$-filtered vector space, there is an \emph{associated graded vector space}
\[
	\Gr V^\bullet = \bigoplus_{k\in\ZZ} \Gr_k V,\qquad \text{where}\quad \Gr_k = \faktor{V^k}{V^{k+1}}.
\]
If the filtration $V^\bullet$ is finite then so is the grading on $\Gr V^\bullet$, and if $V^\bullet$ has finite depth $h$, so does $\Gr V^\bullet$. The associated graded vector space $\Gr V^\bullet$ is (non-canonically) isomorphic to~$V^\bullet$, which can be seen as follows. For each $V^k$, we choose a complement $W_k$ to $V^{k+1}$ so that~${V^k=W_k\oplus V^{k+1}}$. Then we have a grading $V=\bigoplus_{k\in\ZZ} W_k$, and the quotient map $V^k\to\Gr_k V^\bullet$ restricts to a linear isomorphism $W_k\to\Gr_kV^\bullet$. Taking the sum of these maps gives us a graded vector space isomorphism $V\to\Gr V^\bullet $.

If $V^\bullet$ is equipped with a Lie algebra structure, then it is \emph{$\ZZ$-filtered as an algebra} if $\comm{V^k}{V^l}\subseteq V^{k+l}$. Its associated graded vector space is naturally a $\ZZ$-graded Lie algebra; for $X\in \Gr_k V^\bullet, Y\in \Gr_l V^\bullet$, we choose representatives $\tilde{X}\in V^k,\tilde{Y}\in V^l$ and then define a bracket $\comm{-}{-}_{\mathrm{Gr}}$,
\[
	\comm{X}{Y}_{\mathrm{Gr}} := \overline{\big[\tilde{X},\tilde{Y}\big]},
\]
where the overline denotes the natural projection $V^{k+l}\to\Gr_{k+l}V^\bullet$. One can check that this bracket is well-defined precisely because of the condition that the algebra operation on $V^\bullet$ respects the filtration, it satisfies the appropriate axioms, and it is graded by construction.

We note that if $V=\bigoplus_{k\in\ZZ}V_k$ is a graded Lie algebra then it possesses a natural filtration~${V^k=\bigoplus_{l=k}^\infty V_l}$ which has finite depth $h$ if and only if the grading does. Note that we have~${V^k=V_k\oplus V^{k+1}}$, and there is a canonical isomorphism of graded vector spaces $V\to\Gr V^\bullet$ defined as for a general choice of complements $W_k$ as described above. In fact, it is easy to see that this is an isomorphism of graded algebras.

A \emph{filtered morphism} $f\colon V^\bullet\to W^\bullet$ of filtered vector spaces is a vector space morphism~$f\colon V\to W$ such that $f\bigl(V^i\bigr)\subseteq W^i$. We analogously define a \emph{filtered morphism} of filtered Lie superalgebras. A filtered morphism induces an associated graded morphism $\Gr f\colon \Gr V^\bullet \to W^\bullet$ by~${\Gr f([v])=[f(v)]}$. One can easily check that this is well-defined and respects gradings. We say that a filtered morphism is \emph{strict} if $f\bigl(V^k\bigr)=f(V)\cap W^k$.

If $f$ is an injective (respectively surjective) strict filtered morphism, then $\Gr f$ is injective (respectively surjective). Moreover, a linear isomorphism which is filtered is strict if and only if its inverse is filtered, and the inverse of a strict filtered isomorphism is also strict.\footnote{Analogous results do not hold for non-strict filtered morphisms, hence the need to introduce this technicality. See, e.g., \cite[\href{https://stacks.math.columbia.edu/tag/0108}{Example 0108}]{stacks-project}.}
We say that two filtered vector spaces (resp.\ Lie superalgebras) are \emph{isomorphic} if there exists a \emph{strict} filtered isomorphism between them.

\subsubsection{Filtered vector superspaces and superalgebras}
\label{sec:filtered-superspace-superalg}

If $V=V_{\overline{0}}\oplus V_{\overline{1}}$ is a vector superspace, we say that a vector space filtration $V \supseteq \dots \supseteq V^{k-1} \supseteq V^k \supseteq V^{k+1} \supseteq \dots \supseteq 0$ (or $V^\bullet$ for short) is a filtration of the superspace if it is compatible with the $\ZZ_2$-grading in the following sense. We demand that the filtered components are graded subalgebras
\begin{equation}\label{eq:filter-super-1}
	V^k = V^k_{\overline{0}} \oplus V^k_{\overline{1}},
\end{equation}
(note that $V^\bullet_{\overline{0}}$ and $V^\bullet_{\overline{1}}$ are filtrations of the even and odd subspace respectively) such that
\begin{equation}\label{eq:filter-super-2}
	V^{2k-1}_{\overline{0}} = V^{2k}_{\overline{0}}
	\qquad	\text{ and } \qquad
	V^{2k}_{\overline{1}} = V^{2k+1}_{\overline{1}}.
\end{equation}
Note that this means $\Gr_{2k+1}V^\bullet_{\overline 0}=\Gr_{2k}V^\bullet_{\overline 1}=0$, while $\Gr_{2k} V^\bullet=\Gr_{2k} V^\bullet_{\overline0}$ and $\Gr_{2k+1} V^\bullet = \Gr_{2k+1} V^\bullet_{\overline 1}$. This ensures that one can define a natural $\ZZ_2$-grading on $\Gr V^\bullet$ compatible with the $\ZZ$-grading.\footnote{These compatibility conditions seem to be somewhat under-emphasised in the literature; the need for such conditions is mentioned in the classic reference \cite{Kac1977}, but sufficient conditions are only given for a particular class of filtrations there, not for a general filtration. Our main reference \cite{Cheng1998} also mentions that a compatibility condition is required but does not spell it out. The conditions given here appear in \cite{Scheunert1979}.}

If $V^\bullet$ is a Lie superalgebra which is $\ZZ$-filtered as a vector superspace, then it is \emph{$\ZZ$-filtered as a Lie superalgebra} if $\comm{V^k}{V^l}\subseteq V^{k+l}$. One then defines a Lie superalgebra bracket on the associated graded superspace $\Gr V^\bullet$ in the obvious way.

To avoid overloading notation, in the main text we often drop the $\bullet$ superscript from filtered vector spaces and the $\mathrm{Gr}$ subscript from the operation on graded algebras; the filtration structure should be implicitly understood whenever we discuss such spaces. We note that a bullet superscript is also used for cohomological grading elsewhere in this work.
	
\subsection{Homogeneous spaces}
\label{sec:homog-spaces}

Here we establish some terminology, conventions and important results for working with homogeneous spaces. Such spaces play an important role in this work due to the homogeneity theorem (see Theorem~\ref{thm:homogeneity}): a highly supersymmetric supergravity background is (at least locally) homogeneous.

We give our own presentation of standard material here, but some authoritative references for material with no explicit citation in the text are the books of Kobayashi--Nomizu \cite{Nomizu1969} and Sharpe \cite{Sharpe1997}. Wise's work \cite{Wise2010} on the Cartan-geometric formulation of gravity also contains a~brief pedagogical introduction to (metric) Klein pairs.
	
\subsubsection{Homogeneous spaces and Klein pairs}

\begin{Definition}
	Let $G$ be a Lie group. A \emph{homogeneous $G$-space} is a manifold $M$ on which $G$ acts smoothly and transitively. A morphism of such spaces is a $G$-equivariant smooth map, and an isomorphism is a $G$-equivariant diffeomorphism.
	The \emph{isotropy group} of a point $p$ in such a~space is the stabiliser subgroup $G_p=\{g\in G \mid {g\cdot p = p} \}$.
\end{Definition}

Let $G$ be a Lie group and $M$ a homogeneous $G$-space. Fixing a point $p\in M$, we note that the \emph{orbit map} $\alpha_p\colon G\to M$ given by $\alpha_p(g) = g\cdot p$ is surjective (by transitivity), and $G_p=\alpha_p^{-1}(p)$, so~$G_p$ is a closed subgroup of $G$. If $G$ is connected, $M$ must also be connected since it is the image of $G$ under $\alpha_p$. We will often take $G$ to be connected when discussing homogeneous $G$-spaces. Note that all of the isotropy groups of the action of $G$ on $M$ are conjugate; for all~${g\in G}$ we have $G_{g\cdot p} = gG_p g^{-1}$.


The derivative at $p\in M$ of the diffeomorphism $\phi_g\colon M\to M$ associated to an element $g\in G$ is an invertible linear map $d_p\phi_g\colon T_pM\to T_{g\cdot p}M$; in particular, if $h\in G_p$ then $d_p\phi_h\in \GL(T_pM)$. This defines a representation of $G_p$ on $T_pM$ known as the \emph{isotropy representation}.

The derivative of $\alpha_p$ at the identity of $G$ is a linear map $d_e\alpha_p\colon \fg=T_eG\to T_pM$ with kernel~${\fg_p=\Lie G_p}$. One can show that $d_e\alpha_p\circ\Ad_h=d_p\phi_h\circ d_e\alpha_p$ for all $h\in G_p$; in other words, $d_e\alpha_p$ is $G_p$-equivariant where $G_p$ acts on $\fg$ via the restriction of the adjoint representation of $G$. Thus there exists a short exact sequence of $G_p$-modules
\[
\begin{tikzcd}
	&0 \ar[r] &\fg_p \ar[r] & \fg \ar[r,"d_e\alpha_p"] & T_pM \ar[r] & 0.
\end{tikzcd}
\]

The isotropy representation is of fundamental importance in the geometry of homogeneous spaces due to a particular manifestation of \emph{Frobenius reciprocity}: $G$-invariant tensor fields on~$M$ are in one-to-one correspondence with $G_p$-invariant tensors on $T_p M$. This correspondence is easy to describe: if $\tau$ is a $G$-invariant tensor field then its value $\tau_p$ at $p$ is clearly a $G_p$-invariant tensor on $T_p M$; conversely, if $t$ is a $G_p$-invariant tensor on $T_p M$ then we can define a tensor field $\tau$ on $M$ by defining $\tau_{g\cdot p}:=(d_p\phi_{g})(t)$. This tensor field $\tau$ is well-defined and $G$-invariant precisely because $t$ is $G_p$-invariant and $d\phi_{gg'}=d\phi_g\circ d\phi_{g'}$. We will see a more general version of Frobenius reciprocity in due course (see Proposition~\ref{thm:Frobenius-recipr}).

\begin{Definition}\label{def:lie-pair-algs}
	A \emph{Lie pair} is a pair $(\fg,\fh)$ where $\fg$ is a~finite-dimensional Lie algebra and $\fh$ a~subalgebra of $\fg$.
	A \emph{Klein pair} is a pair $(G,H)$ where $G$ is a connected Lie group and $H$ is a~closed subgroup of~$G$.
	The \emph{homogeneous $G$-space associated to the Klein pair $(G,H)$} is the left coset space $G/H$.
\end{Definition}

If $(G,H)$ is a Lie pair, the coset space $G/H$ is indeed a homogeneous $G$-space; clearly $G$ acts on it transitively, and it is a smooth manifold because $H$ is closed. It is connected because $G$ is connected, and $H$ is the isotropy group for $o:=H$ (the coset of the identity). Note that we assume that $G$ is connected but $H$ may not be connected. We will discuss the effect of this on the topology of the associated homogeneous space in Appendix~\ref{sec:homog-spaces-principal-bundles}.

Conversely, if $M$ is a homogeneous space for a connected group $G$ and $p\in M$ then $(G,G_p)$ is a Klein pair, and $M$ is of course the associated homogeneous space up to isomorphism of $G$-spaces; $M\cong G/G_p$. Fixing a connected Lie group $G$, one can show that there is a one-to-one correspondence between conjugacy classes of closed subgroups of $G$ and isomorphism classes of homogeneous $G$-spaces.

The discussion of the isotropy representation above shows that given a Lie pair $(G,H)$, we have $T_o(G/H)\cong \fg/\fh$ as $H$-modules, and that $H$-invariant tensors on $\fg/\fh$ induce $G$-invariant tensor fields on $G/H$. We will be particularly interested in $H$-invariant (pseudo-)inner products~$\eta$ on $\fg/\fh$, which induce $G$-invariant pseudo-Riemannian metrics $g$ on $G/H$.

\begin{Definition}A \emph{metric Lie pair} is a triple $(\fg,\fh,\eta)$ where $(\fg,\fh)$ is a Lie pair and $\eta$ is an $\fh$-invariant inner product on $\fg/\fh$.
A \emph{metric Klein pair} is a triple $(G,H,\eta)$ where $(G,H)$ is a~Klein pair and $\eta$ is an $H$-invariant inner product on $\fg/\fh$.
The \emph{homogeneous pseudo-Riemannian $G$-space associated to $(G,H,\eta)$} is the pair $(G/H,g)$, where $g$ is the pseudo-Riemannian metric induced by $\eta$.
\end{Definition}

Given a metric Klein pair $(G,H,\eta)$, set $M=G/H$, $o=H$ and note that since $\eta$ is $H$-invariant, the image of the isotropy representation, which we now denote by $\varphi\colon H\to GL(T_oM)$, is contained in $\Orth(T_oM)\subseteq\GL(T_oM)$. If $M$ is oriented, then since $G$ is connected, its action is orientation-preserving, hence the action of $H$ on $T_oM$ must also be orientation-preserving, so~${\Im\varphi\subseteq\SO(T_oM)}$. If the signature of $\eta$ is indefinite and $M$ is both oriented and time-oriented, we similarly find that $\Im\varphi\subseteq\SO_0(T_oM)$.

\subsubsection{Homogeneous spaces as principal bundles and their cohomology}
\label{sec:homog-spaces-principal-bundles}


Given a Klein pair $(G,H)$, right multiplication by elements of $H$ on $G$ gives the coset map $G\to M=G/H$ ($g\mapsto gH$) the structure of a (right) $H$-principal bundle over $M$. Note that we \emph{only} consider the \emph{left} action of $G$ on itself, and $G$ does not act on the right of $M$ in a natural way, while $H$ acts on both $G$ and $M$ from the left and right -- the left action on $M$ is non-trivial, while the right action is trivial. The coset map is equivariant with respect to both the left action of $G$ and the right action of $H$.

We have a fibration $H\to G\to M$, hence a long exact sequence of groups in homotopy ending~in
\[
\begin{tikzcd}
	\cdots \ar[r] &\pi_1(G) \ar[r] &\pi_1(M) \ar[r] &\pi_0(H) \ar[r] &\pi_0(G)=1,
\end{tikzcd}
\]
where we take the basepoints of $G$ and $H$ to be the identity element and that of $M$ to be $o=H$. The group structure on $\pi_0(H)$ is induced by that of $H$, and $\pi_0(G)$ is trivial since we assume $G$ to be connected. Thus if $M$ is simply connected, $H$ must be connected. The converse is not true in general, but if we assume that $G$ is simply connected, $\pi_1(M)\cong\pi_0(H)$, so $M$ is simply connected if and only if $H$ is connected.

When working with a homogeneous $G$-space $M=G/H$, it is often useful to pass to the universal cover $\widetilde{G}$ of $G$ by pulling back the action of $G$ on $M$ to $\widetilde{G}$. We can then consider $M$ as a homogeneous space corresponding to the Klein pair \smash{$\bigl(\widetilde{G},H'\bigr)$} where $H'$ is the preimage of $H$ under the covering map \smash{$\widetilde{G}\to G$}. Note that the universal cover $\widetilde{M}$ of $M$ is then a homogeneous space corresponding to $\bigl(\widetilde{G},H'_0\bigr)$, where $H'_0$ is the connected component of $H'$.

\subsubsection{Associated bundles and Frobenius reciprocity}

For the remainder of this section, we let $V=T_o M$ for our convenience. Many natural bundles over~$M$ can be described as associated bundles with respect to representations of $H$. For example, there is a $G$-equivariant isomorphism of vector bundles
$
	TM \cong G\times_H V$,
where $H$ acts on $V$ via the isotropy representation and $G$ acts on $TM$ via the push-forward of the action on $M$ and on the associated bundle via $g\cdot[g',v]=[gg',v]$. The isomorphism in the direction~${G\times_H V\to TM}$ is given by $[g,v]\mapsto g_* v$. The isotropy representation also induces representations of $H$ on~${V^*=T^*_o M}$ and its exterior powers, giving us further $G$-equivariant isomorphisms
$T^*M \cong G\times_H V^*$,	$
	\Wedge^\bullet T^*M \cong G\times_H \Wedge^\bullet V^*$.
The action of $H$ on $\End(V)$ by conjugation by elements of $\GL(V)$ similarly gives us
$\End(TM) \cong G\times_H \End(V)$.

Considering the isotropy representation as a Lie group morphism $H\to\GL(V)$ gives various actions of $H$ on $\GL(V)$ -- by conjugation or by left or right multiplication by the images of elements in $H$, for example -- but we will concentrate on the left action by left multiplication. This allows us form a (right) $\GL(V)$-principal bundle $G\times_H \GL(V)$ equipped with a compatible (left) action of $G$, where the two groups act by right multiplication on the right factor and left multiplication on the left factor respectively. Now, fix a frame $f=(f_1,\dots,f_n)$ at $o$ ($n=\dim M$), that is, $f$ is an ordered basis for $V$. This choice induces a group isomorphism $\GL(V)\cong\GL(\RR^n)$ by representing elements of $\GL(V)$ in the frame $f$; in Einstein notation, for $A\in\GL(V)$ we define the components of the matrix $\underline{A}\in \GL(\RR^n)$ by $A(f_i)=A\indices{^j_i}f_j$. This allows us to treat the frame bundle $FM\to M$ as a $\GL(V)$-principal bundle; if $f'\in FM$ is another frame (at any point), we have
\smash{$(f'\cdot A)_i = f'_jA\indices{^j_i}$}.
This choice of frame $f$ also allows us to define a $G$-equivariant isomorphism of $\GL(V)$-principal bundles (that is, a $(G,\GL(V))$-equivariant bundle morphism)
$FM \cong G\times_H \GL(V)$,
where the left action of $G$ on frames is induced by its action on tangent vectors,
\[
	g\cdot (f'_1,\dots,f'_n) = (g\cdot f'_1,\dots,g\cdot f'_n) = (g_*f'_1,\dots,g_*f'_n);
\]
the isomorphism is given by mapping $G\times_H \GL(V)\ni [g,A]\mapsto g\cdot f \cdot A$.

Given a \emph{metric} Klein pair $(G,H,\eta)$ with orientation, the isotropy representation $H\to\SO(V)$ allows us to construct a $G$-equivariant isomorphism of $\SO(V)$-bundles
$F_{\rm SO} \cong G\times_H \SO(V)$,
where~${F_{\rm SO}\to M}$ is the special orthonormal frame bundle of $(M,g)$, in a completely analogous way -- the construction requires a choice of \emph{oriented orthonormal} frame $f$. Similarly, without an orientation we have $F_O\cong G\times_H\Orth(V)$, and with time orientation we have $F_{SO_0}\cong G\times_H\SO_0(V)$, where the notation should be self-explanatory.

We can make particularly fruitful use of these isomorphisms with Frobenius reciprocity, which we have already seen an application of. This a powerful result which appears in a number of different guises in the literature but in general relates representations of a group $G$ restricted to a subgroup $H$ to representations of $G$ induced by representations of $H$. The most useful version to us here is due to Bott.

\begin{Proposition}[Frobenius reciprocity theorem \cite{Bott1965}]\label{thm:Frobenius-recipr}
	Let $(G,H)$ be a Klein pair, let $W$ be a~$G$-module and let $W'$ be an $H$-module. Then there is a an isomorphism of vector spaces
	\[
		\Hom_G(W,\Gamma(G\times_H W'\to G/H)) \simeq \Hom_H(W,W').
	\]
\end{Proposition}

In the statement above, we treat $W$ as an $H$-module by restricting the action of $G$ on $W$ and we treat the space of sections $\Gamma(G\times_H W'\to G/H)$ as a $G$-module by acting on the left factor of each fibre. We recover our previous statement about invariant tensor fields by noting, for instance, that, setting $W=\RR$ the trivial module and $W'=V=T_oM$, the $G$-equivariant bundle isomorphism $TM \cong G\times_H V$ gives
\smash{$\fX(M)^G \simeq V^H$},
so $G$-invariant vector fields on $M$ are in correspondence with $H$-invariant vectors at $o$. Similarly, using the $G$-equivariant bundle morphism $G\times_H \Wedge^\bullet V^* \cong \Wedge^\bullet T^*M$,
we have
\smash{$\Omega^\bullet(M)^G \simeq (\Wedge^\bullet V^*)^H$},
and similarly $\End(TM)^G\simeq \End(V)^H$ etc.

\subsubsection{Connections, Wang's theorem and Nomizu maps}
\label{sec:nomizu}

Just like tensor fields, invariant connections on homogeneous spaces have a simpler description in terms of equivariant linear maps. This formalism was introduced for affine connections on the homogeneous spaces of reductive Lie pairs by Nomizu \cite{Nomizu1954}. We will discuss a generalisation due to Wang \cite{Wang1958}. A good reference for all statements is \cite{Nomizu1969}.

For a Lie pair $(G,H)$, we once again denote $M=G/H$, $V=T_oM\cong\fg/\fh$ where $o=H$. Note that since $H$ acts on $V$ via the isotropy representation which we denote $\varphi\colon H\to\GL(V)$, it acts on $\fgl(V)$ by conjugation.


Suppose that $\pi\colon P\to M$ is a right principal $K$-bundle for some Lie group $K$, and suppose that the (left) action of $G$ on $M$ lifts to a (left) action on $P$ which is compatible with the action of $K$; that is, we have an action of $G$ on $P$ such that $\pi$ is $G$-equivariant and $(g\cdot p)\cdot k=g\cdot (p\cdot k)$ for all $g\in G$ and $k\in K$. Then note that $H$ preserves the fibre $P_o$, and by fixing a~basepoint~${p\in P_0}$ we can define a~group homomorphism $\psi\colon H\to K$ by declaring $h\cdot p = p\cdot\psi(h)$. A~different choice of basepoint changes $\psi$ by conjugation by an element of $k$. We can now state the following.\looseness=1

\begin{Theorem}[\cite{Wang1958}]\label{thm:Wang}
	Let $(G,H)$ be a Lie pair and let $P\to M=G/H$ be a right principal $K$-bundle for some Lie group $K$. Suppose that the left action of $G$ lifts to $P$ compatibly with the action of $K$. Fix a basepoint $p\in P_o$ and let $\psi\colon H\to K$ be the induced group morphism. Then there is a one-to-one correspondence between~$G$-invariant principal connections $\omega\in\Omega^1(P;\fk)$ on~$P$ and linear maps $\Psi\colon \fg\to\fk$ such that
	\begin{enumerate}\itemsep=0pt
		\item[$(1)$] $\Psi\circ\Ad^G_h = \Ad^K_{\psi(h)}\circ\Psi$ for all $h\in H$,
		\item[$(2)$] $\Psi(X)=(d_e\psi)(X)$ for all $X\in\fh$.
	\end{enumerate}
	The correspondence is given by
$\Psi(X) = \omega_p(\xi_X)
$
	where $X\in \fg$ and $\xi_X$ is the fundamental vector field on $P$ associated to $X$. Moreover, the curvature $F_\omega\in\Omega^2(P;\fk)^G$ is given by
	\[
(F_\omega)_p(\xi_X,\xi_Y) = \comm{\Psi(X)}{\Psi(Y)}_\fk - \Psi(\comm{X}{Y}_\fg).
	\]
\end{Theorem}

The final equation uniquely determines $\omega_p$ since $G$-transitivity on the base implies that any fibre is carried to $P_o$ by the action of some element of $G$, $K$ acts freely and transitively on each fibre, and $\omega$ is a $G$-invariant principal connection.

We will refer to the map $\Psi$ as the \emph{Nomizu map}, since it plays the same role in Wang's theorem that the more well-known notion of Nomizu map plays in the reductive case due to Nomizu.

If $P=FM$ is the frame bundle, recall from Appendix~\ref{sec:homog-spaces-principal-bundles} that a choice of frame~$f$ in~$F_oM$ gives~$FM$ the structure of a principal $K=\GL(V)$-bundle and that we have a compatible left action of~$G$. Then we immediately have $\psi=\varphi\colon H\to\GL(V)$, the isotropy representation. Moreover, if $(G,H,\eta)$ is an orientable metric Lie pair, we can similarly take $P=F_{\rm SO}$, ${K=\SO(V)}$ and $\psi=\varphi\colon H\to\SO(V)$. Now recall that an affine connection is equivalent to a principal connection on $FM$, so Wang's theorem tells us that such a connection on a homogeneous space is characterised by a map $\Psi\colon \fg\to\fgl(V)$, and similarly, metric affine connections on a homogeneous pseudo-Riemannian space correspond to maps $\Psi\colon \fg\to\fso(V)$ (subject to the conditions of the theorem).

We finish this section by describing the Nomizu map associated to the Levi-Civita connection for a metric Lie pair $(G,H,\eta)$. This is a generalisation of the description in the reductive case given in \cite{Nomizu1969} to the non-reductive case that appears, e.g., in~\cite{Figueroa-OFarrill2017_1}. We first define a (degenerate) symmetric bilinear form on $\fg$ by $\pair{X}{Y}=\eta\bigl(\overline{X},\overline{Y}\bigr)$ where $\overline{X}$, $\overline{Y}$ denote the images in $V=\fg/\fh$ of $X,Y\in\fg$; note that $\pair{-}{-}$ is $H$-invariant since $\eta$ is. We then define $U\colon \Odot^2\fg\to V$ by the equation
$2\eta(U(X,Y),v) = \pair{X}{\comm{Z}{Y}} + \pair{\comm{Z}{X}}{Y}
$
for all $X,Y\in\fg$, $v\in V$ and $Z\in\fg$ an arbitrary element such that $\overline{Z}=v$. This is well-defined by $\fh$-invariance of $\pair{-}{-}$. We can then define $\widetilde{L}\colon \fg\to\Hom(\fg,V)$ by
\[
	\widetilde{L}(X)Y = \frac{1}{2}\overline{\comm{X}{Y}} + U(X,Y).
\]
One can check that $\widetilde{L}(X)(Y)=0$ for all $X\in\fg$ and $Y\in\fh$, so $\widetilde{L}$ factors through a map $L\colon \fg\to\fgl(V)$ which is a Nomizu map and also satisfies
\begin{enumerate}\itemsep=0pt
	\item $\Im L\in\fso(V)$,
	\item $\widetilde{L}(X)(Y) - \widetilde{L}(X)(Y) - \overline{\comm{X}{Y}} = 0$ for all $X,Y\in \fg$.
\end{enumerate}
These two properties precisely say that the affine connection corresponding to $L$ is metric-compatible and torsion-free, hence it is the Levi-Civita connection.

\subsection*{Acknowledgements}

The author would like to thank Jos\'e Figueroa-O'Farrill, under whose supervision this work was done, for his guidance and patience. Thanks also to Andrea Santi, James Lucietti, C.S.~Shahbazi and the sorely missed Paul de Medeiros for many enlightening conversations and helpful comments. Finally, thanks to the anonymous referees, whose comments and suggestions greatly improved the quality and clarity of this work.

This work was carried out with scholarship funding from the Science and Technologies Facilities Council (STFC) and the School of Mathematics at the University of Edinburgh.

This work previously appeared as part of the author's Ph.D.~Thesis \cite{Beckett2024} but has not been published elsewhere. Some changes in notation and terminology have been made for this version to aid comparison with other work, and there are also some minor changes in exposition and corrections; in particular, Lemma~\ref{lemma:bundle-dirac-current-exist} and its proof have been updated.

\LastPageEnding

\end{document}